\documentclass{amsart}
\usepackage{amssymb,amsmath,amscd}

    \def\nz{\hbox{\text{0}{ \raise 3pt \hbox{\kern -10pt
    \vrule width 8pt height 0.4 pt }\kern 2pt}}}

\begin{document}
\baselineskip 15 pt
\title{TRANSLATES OF POLYNOMIALS}
\thanks{
Mathematics Department, Purdue University, 
West Lafayette, IN 47907, USA,\\
\hspace*{5mm}e-mail: ram@cs.purdue.edu\\
\hspace*{4mm}
Mathematics Department, Purdue University, 
West Lafayette, IN 47907, USA,\\
\hspace*{4mm}
e-mail: heinzer@math.purdue.edu\\
\hspace*{4mm}
Mathematics Department, University of Kentucky, 
Lexington, KY 40506, USA,\\
\hspace*{4mm}
e-mail: sohum@ms.uky.edu\\
2000 Mathematical Subject Classification: 12F10, 14H30, 20D06, 20E22. 
Abhyankar's work was partly supported by NSF Grant DMS 99-88166 and
NSA grant MDA 904-99-1-0019.}
\author{By Shreeram S. Abhyankar, William J. Heinzer, and Avinash Sathaye}
\date{}
\begin{abstract}
We undertook this study of affine pencils especially to celebrate the
70th birthday of our good friend C. S. Seshadri. The first named author met
Seshadri in Paris in 1958 and had the pleasure of seeing him frequently ever
since. We are very happy to say to him: JEEVEMA SHARADAH SHATAM.
\end{abstract}
\maketitle

\centerline{\bf Section 1: Introduction}

Let $f=0$ be a hypersurface in the $n$-dimensional affine space over a field
$k$ with $n>1$, i.e., $f\in R\setminus k$ where $R$ is the polynomial ring
$k[X_1,\dots,X_n]$. We want to consider the pencil of hypersurfaces
$f-c=0$ with $c$ varying over $k$, and wish to consider the sets
singset$(f)=\{c\in k:f-c\text{ is singular}\}$ and
redset$(f)=\{c\in k:f-c\text{ is reducible}\}$.
In the first case 
$f-c$ is singular means the local ring $R_P/((f-c)R_P)$ is nonregular 
for some $P\in\text{spec}(R)$ with $f-c\in P$, 
and in the second case
$f-c$ is reducible means $f-c=gh$ with $g,h$ in $R\setminus k$.

As consequences of the two famous theorems of Bertini, which may be
called Bertini I or Singular Bertini, and Bertini II or Reducible Bertini,
it can be shown that, under suitable conditions, singset$(f)$ and 
redset$(f)$ are finite. One of our aims is to give short direct proofs of 
these consequences. We shall do this in Sections 3 and 4. In Section 2
we shall recall the Theorems of Bertini and also the Theorem of L\"uroth. 
Before outlining the contents of the rest of the paper, let us fix
some notation. 

By $|S|$ we denote the cardinality of a set $S$.
By $U(S)$ we denote the multiplicative group of all units in a ring $S$,
and by $S^{\times}$ we denote the set of nonzero elements in it.
By QF$(S)$ we denote the quotient field of a domain $S$.
By $A=R/(fR)$ we denote the affine coordinate ring of $f=0$, and we identify
$k$ with its image under the
residue class epimorphism $\phi:R\to A$;
if $f$ is irreducible in $R$ then by $L=\text{QF}(A)$ we denote the 
function field of $f=0$.
By an affine domain over a field $k'$ we mean an overdomain of $k'$ which 
is a finitely generated ring extension of $k'$.
By a DVR we mean a real discrete valuation ring; if the said ring has
quotient field $L'$ then we call it a DVR of $L'$; if it also contains a 
subfield $k'$ then we call it a DVR of $L'/k'$.
Note that a finitely generated free abelian group is isomorphic to
$\mathbb Z^r$ for a unique nonnegative integer $r$ which is called its rank;
we shall apply this to the multiplicative group $U(A')/U(k')$ for
an overdomain $A'$ of a field $k'$.

Given polynomials $g_1,\dots,g_m$ in one or more variables with 
coefficients in some field, we write gcd$(g_1,\dots,g_m)=1$ or $\ne 1$
to mean that they do not or do have a nonconstant common factor.
In particular we shall apply this to the partial derivatives
$f_{X_1},\dots,f_{X_n}$ of $f$.

Let $R^*=k^*[X_1,\dots,X_n]$ where $k^*$ is an algebraic closure of $k$. Let 
singset$(f)^*=\{c\in k^*:R^*_P/((f-c)R^*_P)\text{ is nonregular for some }
P\in\text{spec}(R^*)\text{ with }f-c\in P\}$, and
redset$(f)^*
=\{c\in k^*:f-c=gh\text{ for some }g,h\text{ in }R^*\setminus k^*\}$. 
We shall also consider the multiple set and the primary set of $f$ defined
by putting
multset$(f)^*=\{c\in k^*:f-c=gh^2\text{ for some }g\in R^*\setminus\{0\}
\text{ and }h\in R^*\setminus k^*\}$, and
primset$(f)=\{c\in k:f-c=gh^\mu\text{ for some }g\in k^{\times}
\text{ and }h\in R\setminus k\text{ and integer }\mu>1\}$.

In Section 5 we shall discuss the notion of composite pencils, and
we shall find some bounds for $|\text{redset}(f)|$.
A series of Examples, to be outlined in Remark 5, will illustrate the
spread of the various values which $|\text{redset}(f)|$ can take.

In Section 6 we shall extend our study to more general pencils $f-cw$ where
$w\in R^{\times}$ with gcd$(f,w)=1$.
In Remark 8 we shall link-up the redset of a general pencil to Klein's
parametrization of a special rational surface, and in Question 4 we shall
pose a related problem.

In Section 7 we shall employ a refined version of redset$(f)$ to give
necessary conditions for the ring $R[1/f]$ to be isomorphic to the ring
$R[1/f']$ where $f'=0$ is another hypersurface. This is when $f$ and $f'$
are irreducible. In the reducible case we shall give necessary conditions
in terms of the quotient group $U(R[1/f])/U(k)$ and the function fields of the
irreducible components of $f=0$. 
Now an isomorphism of the rings $R[1/f]$ and $R[1/f']$ can be geometrically
paraphrased as biregular equivalence of the complements of the hypersurfaces
$f=0$ and $f'=0$ in affine $n$-space; moreover, an automorphism of $R$
sending $f$ to $f'$ induces such an isomorphism.
In Questions 5 to 8 of Section 7, these facts provide a link-up of the results
of that section to: Abhyankar's theorem on exceptional nonruled 
varieties,
the birational invariance of the arithmetic genus of a nonsingular projective
variety via the domination part of Abhyankar's desingularization 
theory,
the epimorphism theorems and problems discussed by Abhyankar in his Kyoto 
Notes, and the work of Zariski, Fan, Teicher, and others on the topology of 
complements.
We were originally motivated to consider the said isomorphisms by a
question which Roger Wiegand asked us in 1988.

In Section 8 we shall deduce the finiteness of the redset of a  hypersurface
from that of a plane curve by the intervention of Zariski's famous Lemma 5.
It is by means of this Lemma that Abhyankar proved the Galois case of the
Jacobian Problem. More precisely he deduced the Galois case from the
birational case which was itself proved by using Zariski's Main Theorem.

In Section 9, we shall find a bound for the singset of plane curve $f$ 
in terms of its deficiency set defset$(f)$ which is the set of all constants
$c$ for which the algebraic rank $\rho_a(f-c)$ is different from the 
pencil-rank $\rho_\pi(f)$.
In the complex case, $\rho_a(f)$ of an irreducible $f$ coincides with its
first homology rank. In the general case, $\rho_a(f)$ is defined in terms
of the genera and numbers of branches of the various irreducible components
of $f$, and the pencil-rank $\rho_\pi(f)$ is the general value of $\rho_a(f-c)$
taken over all constants $c$. We shall show how the pencil rank is related
to the Zeuthen-Segre invariant.
We shall also discuss Jung's formula which relates the rank to the
Zeuthen-Segre invariant. 

In Section 10 we shall extend our study of the deficiency set to that
of a general pencil. 

To put things in proper perspective, 
in Section 11 we shall briefly talk about linear systems and pencils on normal
varieties, and say a few words about the Zeuthen-Segre invariant of a
nonsingular projective algebraic surface.

\centerline{}

\centerline{\bf Section 2: Theorems of Bertini and L\"uroth}

Considering a linear systems of codimension one subvarieties of an algebraic
variety, and calling it irreducible if its generic member is irreducible,
in Mantra form, Bertini's Theorems may be stated thus:

\centerline{}

BERTINI I OR SINGULAR BERTINI. Outside the singularities of the variety and 
outside its base points, members of an irreducible
linear system do not have variable singularities. 

\centerline{}

BERTINI II OR REDUCIBLE BERTINI. If a linear system, without fixed components,
is not composite with a pencil, then it is irreducible.

\centerline{}

These were obtained by Bertini in his 1882 paper \cite{Ber}. 
They were revisited by Zariski in \cite{Za1} and \cite{Za2}.

\centerline{}

We shall also use the equally hoary:

\centerline{}

THEOREM Of L\"UROTH. If a curve has a rational parametrization then it has
a faithful rational parametrization.

\centerline{}

This is in his 1875 paper \cite{Lur}. 
We need the refined version given by Abhyankar-Eakin-Heinzer in \cite{AEH}.
Also see Igusa \cite{Igu} and Nagata \cite{Na2}.

\centerline{}

\centerline{\bf Section 3: Singset}

Let us now prove our:

\centerline{}

SINGSET THEOREM. If $k$ is of characteristic zero then singset$(f)$
is finite.

\centerline{}

PROOF. Let $I$ be the ideal in $R^*$
generated by $f_{X_1},\dots,f_{X_n}$. For any
$P\in\text{spec}(R^*)$ with $I\subset P$, consider the residue class map
$\Phi_P:R^*\to R^*/P$. Since all the partials of
$f$ belong to $P$, it follows that $D(\Phi_P(f))=0$ for every 
$\Phi_P(k^*)$-derivation of QF$(R^*/P)$. Therefore,
since $k^*$ is of characteristic zero, we have 
$\Phi_P(f)=\Phi_P(\kappa(P))$ for a unique $\kappa(P)\in k^*$.
Clearly $P+(f-c)R^*=R^*$ whenever $\kappa(P)\ne c\in k^*$.
Let $P_1,\dots,P_s$ be the minimal primes of $I$ in $R^*$, where 
$s=0\Leftrightarrow I=R^*$. 
Then for all $c\in k^*\setminus\{\kappa(P_1),\dots,\kappa(P_s)\}$  we have
$I+(f-c)R^*=R^*$. Since $k$ is of characteristic zero, it follows that
singset$(f)\subset\{\kappa(P_1),\dots,\kappa(P_s)\}$, and hence 
singset$(f)$ is finite.

\centerline{}

\centerline{\bf Section 4: Redset}

\centerline{}

Next we prove our:

\centerline{}

REDSET THEOREM. If $f$ is irreducible in $R$ and $k$ is relatively
algebraically closed in $L$, then redset$(f)$ is finite. 

\centerline{}

PROOF. By the following Lemma we can find a finite number of DVRs 
$V_1,\dots,V_t$ of $L/k$ such that $A\cap V_1\cap\dots\cap V_t=k$.
For every $z\in L^{\times}$ let 
$W_i(z)=\text{ord}_{V_i}(z)$, and let $W:L^{\times}\to\mathbb Z^t$ be the map
given by putting $W(z)=(W_1(z),\dots,W_t(z))$.
Let $G$ be the set of all $g\in R\setminus k$ such that $gh=f-c$ for some
$h\in R\setminus k$ and $c\in k^{\times}$. Since the degree
of $g$ is clearly smaller than the degree of $f$, the set $G$ is contained
in a finite dimensional $k$-vector-subspace of $R$. Therefore for every
$i\in\{1,\dots,t\}$, the set $W_i(\phi(g))_{g\in G}$ is bounded from below.
Since $h$ also belongs to $G$ and clearly $W_i(\phi(g))=-W_i(\phi(h))$,
it follows that the set $W_i(\phi(g))_{g\in G}$ is also bounded from above.
Since, for $1\le i\le t$, 
the set $W_i(\phi(g))_{g\in G}$ is bounded from both sides, it follows that
$W(\phi(G))$ is a finite set. Also clearly $\phi(G)\subset U(A)$.
Let $g_1h_1=f-c_1$ and $g_2h_2=f-c_2$ with $g_1,h_1,g_2,h_2$ in 
$R\setminus k$ and $c_1,c_2$ in $k^{\times}$ be such that
$W(\phi(g_1))=W(\phi(g_2))$.
Then $\phi(g_1)/\phi(g_2)\in A\cap V_1\cap\dots\cap V_t=k$ and  hence
$\phi(g_2)=c\phi(g_1)$ for some $c\in k^{\times}$. 
Consequently $g_2-cg_1$ is divisible by $f$ in $R$ and hence,
because $\text{deg}(g_2-cg_1)<\text{deg}(f)$, we must have $g_2=cg_1$.
Therefore, by
subtracting the equation $g_2h_2=f-c_2$ from the
equation $g_1h_1=f-c_1$ we get
$c_2-c_1=g_1h_1-g_2h_2=g_1(h_1-ch_2)$ which implies that $c_2-c_1\in k$ 
is divisible in $R$ by the positive degree
polynomial $g_1$. Consequently we must have $c_2=c_1$. Therefore, because
the set $W(\phi(G))$ is finite, we conclude that redset$(f)$ is finite.

\centerline{}

LEMMA. Given any affine domain $A'$ over a field $k'$, let $k''$ be the
algebraic closure of $k'$ in $L'=\text{QF}(A')$, and let $A''$ be the integral
closure of $A'$ in $L'$. Then
there exists a finite number of DVRs $V_1,\dots,V_t$ of $L'/k'$ such that
$A''\cap V_1\cap\dots\cap V_t=k''$. Moreover, if $k''=k'$ then
$A'\cap V_1\cap\dots\cap V_t=k'$ and $U(A')/U(k')$ is a finitely
generated free abelian group of rank $r$ with $r\le\text{max}(0,t-1)$.

\centerline{}

PROOF. Let $m$ be the transcendence degree of $L'$ over $k'$.
If $m=0$ then we can take $t=0$, and if also $k'=k''$ then
clearly $U(A')/U(k')=1$ and hence $r=0$.
So assume $m>0$. Then by Noether 
Normalization $A'$ is integral over a polynomial ring
$\widehat A=k'[Y_1,\dots,Y_m]\subset L'$.
Let $\widehat L=k'(Y_1,\dots,Y_m)$ and
$\widehat V=\{g/h:g,h\in\widehat A\text{ with }h\ne 0
\text{ and deg}(g)\le\text{deg}(h)\}$. 
Then $\widehat V$ is a DVR of $\widehat L/k'$
with $\widehat A\cap\widehat V=k'$.  	
Let $V''$ be the integral closure of $\widehat V$ in $L'$. Then
$V''=V_1\cap\dots\cap V_t$ where $V_1,\dots,V_t$ are DVRs of $L'/k'$ with
$t>0$, and by the following Sublemma we have $A''\cap V''=k''$. 
It follows that $A''\cap V_1\cap\dots\cap V_t=k''$, and if $k''=k'$ then
$A'\cap V_1\cap\dots\cap V_t=k'$.
Now assuming $k''=k'$, let $W:U(L')\to\mathbb Z^t$ be the homomorphism
(from a multiplicative group to an additive group) given by
$W(z)=(W_1(z),\dots,W_t(z))$ with $W_i(z)=\text{ord}_{V_i}(z)$. 
Then $U(A')\cap\text{ker}(W)=U(k')$ and hence we get a monomorphism
$\overline W:U(A')/U(k')\to\mathbb Z^t$. Therefore 
$U(A')/U(k')$ is a finitely generated free abelian group of rank $r$ with 
$r\le t$. Suppose if possible that $r=t$.
Then we can find $z_1,\dots,z_t$ in $U(A')$ such that the $t\times t$
matrix $W_i(z_j)$ has a nonzero determinant. Now the column vectors of
this matrix are $\mathbb Q$-linearly independent vectors in $\mathbb Q^t$ 
and hence the column vector $(1,0,\dots,0)$ can be expressed as a
$\mathbb Q$-linear combination of them, i.e., we can find rational
numbers $a_1,\dots,a_t$ such that $\sum_{1\le j\le t}W_i(z_j)a_j=1$ or $0$
according as $i=1$ or $2\le i\le t$. Next we can find integers
$a,b_1,\dots,b_t$ with $a>0$ such that $aa_j=b_j$ for $1\le j\le t$.
Let $z=\prod_{1\le j\le t}z_j^{b_j}$. Then $z\in U(A')$ with
$W_1(z)=a>0=W_i(z)$ for $2\le i\le t$.
Consequently $z\in A'\cap V_1\dots\cap V_t$ but $z\not\in k'$ which is a
contradiction. Therefore we must have $r\le t-1$. 

\centerline{}

SUBLEMMA. Let $\widehat A$ and $\widehat V$ be normal domains with a common
quotient field $\widehat L$. Let $A''$ and $V''$ be the respective integral
closures of $\widehat A$ and $\widehat V$ in an algebraic field extension
$L'$ of $\widehat L$. Then $A''\cap V''$ is the integral closure of
$\widehat A\cap \widehat V$ in $L'$.

\centerline{}

PROOF. Clearly the integral closure of $\widehat A\cap\widehat V$ in $L'$
is contained in $A''\cap V''$. Conversely, given any $z\in L'$ let $H(Z)$
be the minimal monic polynomial of $z$ over $\widehat L$. By
Kronecker's Theorem (or obviously), 
if $z\in A''$ then $H(Z)\in\widehat A[Z]$, and
if $z\in V''$ then $H(Z)\in\widehat V[Z]$. Consequently,
if $z\in A''\cap V''$ then $H(Z)\in(\widehat A\cap\widehat V)[Z]$.
Thus $A''\cap V''$ is integral over $\widehat A\cap\widehat V$.
Therefore $A''\cap V''$ is the integral closure of
$\widehat A\cap \widehat V$ in $L$.

\centerline{}

REMARK 1. With notation as in the proof of the Redset Theorem, 
in the $n=2$ case, the Lemma also follows by taking 
$V_1,\dots,V_t$ to be the valuation rings of the places at infinity of the 
plane curve $f=0$. Moreover, in this case, in the proof of the Redset
Theorem, the finiteness of
$W(\phi(G))$ can be deduced from the boundedness from below of the
sets $W_1(\phi(G)),\dots,W_t(\phi(G))$ by invoking the fact that
the number of zeros of a rational function on the curve $f=0$ equals the
number of its poles.

\centerline{}

REMARK 2. Assuming $k$ to be of characteristic zero, the $n>2$ case of
the Redset Theorem can be deduced from the $n=2$ case by invoking the
famous Lemma 5 of Zariski's paper \cite{Za1} thus.
By applying a linear $k$-automorphism to $R$ and multiplying $f$ by an
element in $k^\times$, we can arrange $f$ to be a monic polynomial of 
positive degree in $X_1$ with coefficients
in $S=k[X_2,\dots,X_n]$. Now $\phi(X_2),\dots,\phi(X_n)$ is a 
transcendence basis of $L/k$ and hence, invoking Zariski's Lemma 5 and 
applying a $k$-linear automorphism to $S$, we can arrange the field
$k(\phi(X_3),\dots,\phi(X_n))$ to be relatively algebraically closed in
$L$. Let $\widetilde R=\widetilde k[X_3,\dots,X_n]$ with
$\widetilde k=k(X_3,\dots,X_n)$, and let 
$\widetilde L=\text{QF}(\widetilde A)$ with
$\widetilde A=\widetilde R/(f\widetilde R)$.
By Gauss Lemma, $f$ is irreducible in $\widetilde R$ and hence, by the $n=2$
case of the Redset Theorem, the set
$\{c\in\widetilde k:f-c\text{ is reducible in }\widetilde R\}$ is finite.
Therefore, again by Gauss Lemma, the set
$\{c\in k:f-c\text{ is reducible in }R\}$ is finite.

\centerline{}

REMARK 3. Assume $f$ is irreducible in $R$.
Now the Redset Theorem says that if the ground field $k$ is relatively
algebraically closed in the function field $L$ then redset$(f)$, i.e.,
the set of all constants $c$ for which the polynomial $f-c$ factors, is
finite. We want to observe that this set, although determined by $f$,
is not determined by the affine coordinate ring $A$. Indeed if redset$(f)$
is nonempty, say it contains a constant $c$ (which is necessarily nonzero
because $f$ is irreducible), then every constant $\gamma$ in $A$ factors. 
Namely, since $c$ is in redset$(f)$, we can write $f-c=gh$ with $g,h$ in 
$R\setminus k$, and multiplying both sides by $-\gamma/c$ we get
$(-\gamma/c)f+\gamma=g'h'$ with $g'=g$ and $h'=(-\gamma/c)h$.
Thus we have factored $\gamma$ not only in $A$, i.e., modulo the ideal
$fR$, but also ``modulo'' the constant multiples $kf$ of $f$.
Note that if $\gamma\ne 0$ then both $g'$ and $h'$ belong to
$R\setminus k$, but if $\gamma=0$ then $h'$ does not. To take care of
this extra ``desire'' and to make sure that the ``multiplier $e=-(\gamma/c)$''
of $f$ is also nonzero, by taking $e=ug-(\gamma/c)$ with any 
$u\in R^{\times}$ we get $ef+\gamma=g'h'$ with $g'=g\in R\setminus k$ and
$h'=uf-(\gamma/c)h\in R\setminus k$.
As example, $f=X_1X_2+1$ is irreducible with
redset$(f)=\{1\}$, and for any $\gamma\in k$
we have $(X_1X_2-\gamma)f+\gamma=X_1[X_1X_2^2+(1-\gamma)X_2]$.

\centerline{}

\centerline{\bf Section 5: Generic Members and Composite Pencils}

Let $R^{\sharp}=k(Z)[X_1,\dots,X_n]$
where $Z$ is an indeterminate over $R$.
By the generic member of the pencil $(f-c)_{c\in k}$ we mean the
hypersurface $f^{\sharp}=0$ with $f^{\sharp}=f-Z\in R^{\sharp}$.
Let singset$(f^{\sharp})$ be the set of all $c\in k(Z)$ such that
$R^{\sharp}_P/((f^{\sharp}-c)R^{\sharp}_P)$ is nonregular for some
$P\in\text{spec}(R^{\sharp})$ with $f^{\sharp}-c\in P$, and let
redset$(f^{\sharp})$ be the set of all $c\in k(Z)$ such that $f^{\sharp}-c$
is reducible in $R^{\sharp}$.

By Gauss Lemma $f^{\sharp}$ is irreducible in $R^{\sharp}$, i.e.,
$0\not\in\text{redset}(f^{\sharp})$.
Let $\phi^{\sharp}: R^{\sharp}\to R^{\sharp}/(f^{\sharp}R^{\sharp})$
be the residue class epimorphism. Clearly 
$R\cap\text{ker}(\phi^{\sharp})=\{0\}$ and 
$\phi^{\sharp}(f)=\phi^{\sharp}(Z)$,
and hence there exits a unique isomorphism
$\psi^{\sharp}:k(X_1,\dots,X_n)\to\text{QF}(R^{\sharp}/(f^{\sharp}R^{\sharp}))$
such that for all $r\in R$ we have $\psi^{\sharp}(r)=\phi^{\sharp}(r)$. 
Thus the triple
$\phi^{\sharp}(k(Z))\subset\phi^{\sharp}(R^{\sharp})
\subset\text{QF}(\phi^{\sharp}(R^{\sharp}))$ is isomorphic to the triple
$$
k^{\sharp}=k(f)\subset A^{\sharp}=k(f)[X_1,\dots,X_n]
\subset L^{\sharp}=k(X_1,\dots,X_n)
$$
and hence we may regard the above three displayed sets as
the ground field, the affine coordinate ring, and the function field
of $f^{\sharp}=0$. 

\centerline{}

The ring $A^{\sharp}$ is regular because it is the localization
of the regular ring $R$ at the multiplicative set $k[f]^{\times}$.
Thus we have the:

\centerline{}

GENERIC SINGSET THEOREM. $0\not\in\text{singset}(f^{\sharp})$.

\centerline{}

In view of the above isomorphism of triples, by the Redset Theorem
we get the:

\centerline{}

GENERIC REDSET THEOREM. 
If $k^{\sharp}$ is relatively algebraically closed in $L^{\sharp}$ then 
redset$(f^{\sharp})$ is finite.

\centerline{}

In \cite{AEH} we proved the following:

\centerline{}

REFINED L\"UROTH THEOREM. Assume that
$k^{\sharp}$ is not relatively algebraically closed in $L^{\sharp}$. Let 
$\widehat B^{\sharp}$ be the integral closure of $k[f]$ in $L^{\sharp}$,
let $\widehat k^{\sharp}$ be the algebraic closure of $k^{\sharp}$ in 
$L^{\sharp}$, and let $\nu=[\widehat k^{\sharp}:k^{\sharp}]$.
Then $\nu$ is an integer with $\nu>1$, and there exist
$\widehat f\in R\setminus k$ and $\Lambda\in k[Z]\setminus k$ with
$\widehat B^{\sharp}=k[F]$ and $\widehat k^{\sharp}=k(F)$ such that 
$f=\Lambda(\widehat f)$ and deg$_Z\Lambda=\nu$.

\centerline{}

REMARK 4. In the above situation, every member of the pencil
$(f-c)_{c\in k^*}$ consists of $\nu$ members (counted properly) of the
pencil $(\widehat f-\widehat c)_{\widehat c\in k^*}$, and so we say that 
{\bf the pencil $(f-c)_{c\in k^*}$ is composite with the pencil 
$(\widehat f-\widehat c)_{\widehat c\in k^*}$}.  
In greater detail, for any $c\in k^*$ we have 
$\Lambda(Z)-c=\widehat c_0\prod_{1\le i\le\nu}(Z-\widehat c_i)$ where 
$\widehat c_0,\widehat c_1,\dots,\widehat c_{\nu}$ in $k^*$ with 
$\widehat c_0\ne 0$, and by substituting $\widehat f$ for $Z$ we get
$f-c=\widehat c_0\prod_{1\le i\le\nu}(\widehat f-\widehat c_i)$ and so the 
hypersurface $f=c$ is the union of the hypersurfaces 
$\widehat f=\widehat c_1,\dots,\widehat f=\widehat c_{\nu}$. Since 
deg$_Z\Lambda=\nu>1$, if $k$ is of characteristic $0$ then by taking
a root $\zeta$ of the $Z$-derivative of $\Lambda(Z)$ in $k^*$ we see that
$\Lambda(Z)-\Lambda(\zeta)$ has a multiple
root in $k^*$ and hence $f-\Lambda(\zeta)$ has a nonconstant multiple factor 
in $R^*$, and therefore multset$(f)^*\ne\emptyset$.
Note that if $k$ is of characteristic $p>0$ then this does not work as can
be seen by taking $\Lambda(Z)=Z^p+Z$.
Without assuming any condition on the pair $(k^{\sharp},L^{\sharp})$ we
see that, for any $c\in k^*$, every multiple factor of $f-c$ in $R^*$
divides $f_{X_j}$ for $1\le j\le n$, and hence:
multset$(f)^*\ne\emptyset\Rightarrow\text{gcd}(f_{X_1},\dots,f_{X_n})\ne 1$.
Again without assuming any condition on the pair $(k^{\sharp},L^{\sharp})$ we
see that $c\mapsto c-Z$ gives an injection of redset$(f)$ into
redset$(f^{\sharp})$, and hence
$|\text{redset}(f)|\le|\text{redset}(f^{\sharp})|$. Thus, in view of the
above two Theorems, we get the: 

\centerline{}

COMPOSITE PENCIL THEOREM. We have the following:

(I) 
gcd$(f_{X_1},\dots,f_{X_n})=1\Rightarrow\text{multset}(f)^*=\emptyset$.

(II) 
characteristic of $k$ is $0$ and 
multset$(f)^*=\emptyset\Rightarrow k^{\sharp}$ is relatively algebraically 
closed in $L^{\sharp}$.

(III)
$k^{\sharp}$ is relatively algebraically closed in 
$L^{\sharp}\Rightarrow$ redset$(f)$ and redset$(f^{\sharp})$
are finite with $|\text{redset}(f)|\le|\text{redset}(f^{\sharp})|$.

(IV)
$k^{\sharp}\text{ is not relatively algebraically closed in }
L^{\sharp}\Rightarrow\text{redset}(f)^*=k^*$.

\centerline{}

Next we prove the:

\centerline{}

MIXED PRIMSET THEOREM. If $k^{\sharp}$ is relatively algebraically closed 
in $L^{\sharp}$ then primset$(f)=\emptyset$.

\centerline{}

PROOF. If primset$(f)\ne\emptyset$ then for some $c\in k$ and integer
$\mu>1$ we have $f-c=gh^\mu$ with $g\in k^{\times}$ and $h\in R\setminus k$,
and this implies $[k(h):k(f)]=\mu$ and hence $k^{\sharp}$ is not relatively
algebraically closed in $L^{\sharp}$.

\centerline{}

Let us now prove the following theorem which is some kind of a mixture of
the Redset Theorem and the Generic Redset Theorem.

\centerline{}

MIXED REDSET THEOREM. 
If $k^{\sharp}$ is relatively algebraically closed in $L^{\sharp}$ then 
$U(A^{\sharp})/U(k^{\sharp})$ is a finitely generated free abelian group of
rank $r$ with $|\text{redset}(f)|\le r$.

\centerline{}

PROOF. Assuming that $k^{\sharp}$ is relatively algebraically closed in 
$L^{\sharp}$, by the Lemma we see that 
$U(A^{\sharp})/U(k^{\sharp})$ is a finitely generated free abelian group of
some rank $r$. Suppose if possible that 
$|\text{redset}(f)|>r$ and take distinct elements $c_1,\dots,c_s$ in
redset$(f)$ with integer $s>r$. By the Mixed Primset Theorem, 
for $1\le i\le s$ we have 
\begin{equation*}
f-c_i=g_ih_i\quad\text{ where }\quad
g_i,h_i\text{ in }R\setminus k
\quad\text{ with }\quad\text{gcd}(g_i,h_i)=1.
\end{equation*}
Since $s>r$, there exist integers
$a_1,\dots,a_s$ at least one of which is nonzero, say $a_{\iota}\ne 0$,
such that
\begin{equation*}
\prod_{1\le i\le s}g_i^{a_i}\in U(k(f)).
\tag{2}
\end{equation*}
For the rest of the proof we shall give two alternative arguments.

\centerline{}

FIRST ARGUMENT.
Clearly the images of $g_i$ and $h_i$ in $U(A^{\sharp})/U(k^{\sharp})$ are 
inverses of each other and hence replacing $g_i$ by $h_i$ for those $i$
for which $a_i<0$, we can arrange matters so that $a_i\ge 0$ for all $i$,
and hence in particular $a_{\iota}>0$.
Again since the images of $g_i$ and $h_i$ in $U(A^{\sharp})/U(k^{\sharp})$ 
are inverses of each other, by (2) we get
\begin{equation*}
\prod_{1\le i\le s}h_i^{a_i}\in U(k(f)).
\tag{$3'$}
\end{equation*}
Any element in $k(f)^{\times}$ can be written as $u(f)/v(f)$ where
$u(Z),v(Z)$ in $k[Z]^{\times}$ with $u(Z)\widehat u(Z)+v(Z)\widehat v(Z)=1$
for some $\widehat u(Z),\widehat v(Z)$ in $k[Z]$, and substituting $Z=f$
we get $u(f)\widehat u(f)+v(f)\widehat v(f)=1$; it follows that if
$u(f)/v(f)\in R$ then $v(f)\in k^{\times}$ and hence $u(f)/v(f)\in k[f]$.
Thus $R\cap k(f)=k[f]$ and therefore by (2) and $(3')$ we get
\begin{equation*}
\prod_{1\le i\le s}g_i^{a_i}=g(f)
\quad\text{ and }\quad
\prod_{1\le i\le s}h_i^{a_i}=h(f)
\tag{$4'$}
\end{equation*}
with $g(Z),h(Z)$ in $k[Z]^{\times}$.
Multiplying the two equations in $(4')$ we obtain
\begin{equation*}
\prod_{1\le i\le s}(f-c_i)^{a_i}=g(f)h(f)
\end{equation*}
and hence by the $k$-isomorphism $k[Z]\to k[f]$ with $Z\mapsto f$ we get
\begin{equation*}
\prod_{1\le i\le s}(Z-c_i)^{a_i}=g(Z)h(Z)
\end{equation*}
and therefore we have
\begin{equation*}
g(Z)=\overline g\prod_{1\le i\le s}(Z-c_i)^{\alpha_i}
\quad\text{ and }\quad
h(Z)=\overline h\prod_{1\le i\le s}(Z-c_i)^{\beta_i}
\tag{$5'$}
\end{equation*}
with $\overline g,\overline h$ in $k^{\times}$
and nonnegative integers $\alpha_i,\beta_i$ such that 
$\alpha_i+\beta_i=a_i$. 
Substituting $Z=f$ and $f-c_i=g_ih_i$ in $(5')$ 
and then comparing it with $(4')$ we get
\begin{equation*}
\prod_{1\le i\le s}g_i^{a_i}=
\overline g\prod_{1\le i\le s}(g_ih_i)^{\alpha_i}
\quad\text{ and }\quad
\prod_{1\le i\le s}h_i^{a_i}=
\overline h\prod_{1\le i\le s}(g_ih_i)^{\beta_i}.
\tag{$6'$}
\end{equation*}
Clearly gcd$(f-c_i,f-c_j)=1$ for all $i\ne j$, and hence by (1) and $(6')$
we get
\begin{equation*}
a_i=\alpha_i\text{ and }a_i=\beta_i\text{ for }1\le i\le s.
\tag{$7'$}
\end{equation*}
Since $\alpha_{\iota}+\beta_{\iota}=a_{\iota}>0$ with $\alpha_{\iota}\ge 0$ 
and $\beta_{\iota}\ge 0$, we must have either $\alpha_{\iota}<a_{\iota}$ or 
$\beta_{\iota}<a_{\iota}$ which contradicts $(7')$.
Therefore $|\text{redset}(f)|\le r$.

\centerline{}

SECOND ARGUMENT. We shall make this more condensed.  By (2) we get
\begin{equation*}
\prod_{1\le i\le s}g_i^{a_i}=\overline u
\prod_{1\le i\le \sigma}u_i(f)^{\alpha_i}
\quad\text{ with }\quad\overline u\in k^{\times}
\tag{$3''$}
\end{equation*}
where $\alpha_1,\dots,\alpha_\sigma$ are nonzero integers and
$u_1(Z),\dots,u_\sigma(Z)$ are pairwise coprime monic
members of $k[Z]\setminus k$.
As in the First Argument, 
for any $u,v$ in $k[Z]\setminus k$ with gcd$(u,v)=1$ 
we have $u(f),v(f)$ in $R\setminus k$ with gcd$(u(f),v(f))=1$. 
Therefore, by looking at divisibility by a nonconstant irreducible factor
of $g_i$ in $R$, in view of (1) and $(3'')$, we see that for each 
$i\in\{1,\dots,s\}$ with $a_i\ne 0$ there is some 
$\theta(i)\in\{1,\dots,\sigma\}$ with $u_{\theta(i)}(Z)=Z-c_i$ and
$\alpha_{\theta(i)}=a_i$. In particular we get
$u_{\theta(\iota)}(Z)=Z-c_\iota$ and $\alpha_{\theta(\iota)}=a_\iota$.
Again in view of (1) and $(3'')$, 
by looking at divisibility by a nonconstant irreducible factor of $h_{\iota}$ 
in $R$, we get a contradiction. 
Therefore $|\text{redset}(f)|\le r$.

\centerline{}

REMARK 5. Assume $n=2$. Let $t$ be the number of places at infinity 
of the irreducible plane curve $f^{\sharp}=0$ and let
$V_1,\dots,V_t$ be their valuation rings. Then $t$ is a 
positive integer and by the Lemma and the Mixed Redset Theorem we see that if 
$k^{\sharp}$ is relatively algebraically closed in $L^{\sharp}$, then
$R^{\sharp}\cap V_1\cap\dots\cap V_t=k^{\sharp}$ and 
$|\text{redset}(f)|\le t-1$. We shall show that this bound is the best possible.
In fact, in the following Examples 1 to 5, 
assuming $k$ to contain sufficiently many elements,
we shall show that, 
for any $t>0$, $|\text{redset}(f)|$ can be equal to any integer $\mu$ with 
$0\le\mu\le t-1$; 
Example 1 works for $0\le \mu = t-1$ and $t>0$;
Example 2 works for $1\le \mu < t-1$ and $t>2$;
Example 3 works for $ \mu=0$ and $t>3$;
Example 4 works for $ \mu=0$ and $t>0$;
Example 5 works for $ \mu=1$ and $t>1$.
As common notation for these examples,
let $n=2$, let $m\ge 0$ be an integer, and
let $a(X_1)=(X_1-a_1)\dots(X_1-a_m)$
where $a_1,\dots,a_m$ are pairwise distinct elements in $k$.
Moreover, for the irreducible $f\in R\setminus k$ to be constructed, let 
$\tau(f)$ denote the number of places at infinity of the plane curve $f=0$, and 
let $\tau(f^{\sharp})$ denote the number of places at infinity of the plane 
curve $f^{\sharp}=0$.

\centerline{}

EXAMPLE 1. Let
$f=a(X_1)X_2+X_1-z$ with $z\in k\setminus\{a_1,\dots,a_m\}$. 
Then $f$ is irreducible in $R$ and
redset$(f)=\text{redset}(f)^*=\{a_1-z,\dots,a_m-z\}$. 
Moreover, $k^{\sharp}$ is relatively algebraically closed in $L^{\sharp}$
and we have $|\text{redset}(f)|=|\text{redset}(f)^*|=m$ with 
$\tau(f)=\tau(f^{\sharp})=m+1$.

\centerline{}

PROOF. Clearly $(X_1-a_i)$ divides $f-(a_i-z)$ in $R$ for $1\le i\le m$, and
hence $\{a_1-z,\dots,a_m-z\}\subset\text{redset}(f)$. By Gauss Lemma we 
also see that
redset$(f)^*\subset\{a_1-z,\dots,a_m-z\}$. Therefore $f$ is irreducible 
in $R$, and we have redset$(f)=\text{redset}(f)^*=\{a_1-z,\dots,a_m-z\}$. 
In particular $|\text{redset}(f)|=|\text{redset}(f)^*|=m<\infty$, and hence 
by the Composite Pencil Theorem we see 
that $k^{\sharp}$ is relatively algebraically closed in $L^{\sharp}$.
For the degree form (which consists of the highest degree terms) of $f$
we have defo$(f)=X_1^mX_2$, and hence the points at infinity of the curve
$f=0$ are $(X_0,X_1,X_2)=(0,1,0)$ and $(X_0,X_1,X_2)=(0,0,1)$.
Homogenizing $f$ we obtain
$$
(X_1-a_1X_0)\dots(X_1-a_mX_0)X_2+X_0^mX_1-zX_0^{m+1}
$$ 
Putting $X_1=1$ we get 
$$
f_1=(1-a_1X_0)\dots(1-a_mX_0)X_2+X_0^m-zX_0^{m+1}
$$
and hence $(0,1,0)$ is a simple point. Putting $X_2=1$ we get 
$$
f_2=(X_1-a_1X_0)\dots(X_1-a_mX_0)+X_0^mX_1-zX_0^{m+1}
$$
and hence $(0,0,1)$ is an $m$-fold point with $m$ distinct tangents.
Therefore $\tau(f)=m+1$, and hence by taking $(z+Z,k^{\sharp})$ for
$(z,k)$ we get $\tau(f^{\sharp})=m+1$.

\centerline{}

EXAMPLE 2. Assume $m\ge 1$ and let 
$\mu$ be any integer with $1\le\mu\le m$. 
Let $b(X_1)=(X_1-b_1)\dots(X_1-b_m)$ where
$b_1,\dots,b_m$ be pairwise distinct elements in $k$ such that 
$b_i=a_i$ for $1\le i\le\mu$, and
$b_i\not\in\{a_1,\dots,a_m\}$ for $\mu<i\le m$. 
If $m\ne 1$ then let $\gamma$ be any nonzero element of $k$, and
if $m=1$ then let $\gamma$ be any nonzero element of $k$ such that
$Z^2+\gamma Z+1=(Z-\gamma_1)(Z-\gamma_2)$ with $\gamma_1\ne\gamma_2$ in $k$.
Let $f=a(X_1)X_2^2+\gamma b(X_1)X_2+X_1-z$ with 
$z\in k\setminus\{a_1,\dots,a_{\mu}\}$. 
Then $f$ is irreducible in $R$ and
redset$(f)=\text{redset}(f)^*=\{a_1-z,\dots,a_{\mu}-z\}$. 
Moreover, $k^{\sharp}$ is relatively algebraically closed in $L^{\sharp}$
and we have $|\text{redset}(f)|=|\text{redset}(f)^*|=\mu$ with 
$\tau(f)=\tau(f^{\sharp})=m+2$.

\centerline{}

PROOF. Clearly $(X_1-a_i)$ divides $f-(a_i-z)$ in $R$ for $1\le i\le\mu$, and
hence $\{a_1-z,\dots,a_{\mu}-z\}\subset\text{redset}(f)$. 
Again by Gauss Lemma we see that if for some 
$c\in k^*\setminus\{a_1-z,\dots,a_{\mu}-z\}$ we have $c\in\text{redset}(f)$
then there exists $\eta\in(k^*)^{\times}$ and a disjoint partition
$\{1,\dots,m\}=U\coprod V$ such that upon letting
$$
u(X_1)=\prod_{i\in U}(X_1-a_i)
\quad\text{ and }\quad
v(X_1)=\prod_{i\in V}(X_1-a_i)
$$ 
we have
$$
f-z-c=\left[u(X_1)X_2+\eta\right]\left[v(X_1)X_2+(1/\eta)(X_1-z-c\right].
$$
Equating the coefficients of $X_2$ we get
\begin{equation*}
\gamma b(X_1)=(1/\eta)(X_1-z-c)u(X_1)+\eta v(X_1).
\tag{*}
\end{equation*}
This gives a contradiction because the LHS is divisible by $(X_1-a_1)$,
but exactly one term in the RHS is so divisible. 
Therefore 
redset$(f)^*\subset\{a_1-z,\dots,a_m-z\}$. Consequently $f$ is irreducible 
in $R$, and we have redset$(f)=\text{redset}(f)^*=\{a_1-z,\dots,a_m-z\}$. 
In particular $|\text{redset}(f)|=|\text{redset}(f)^*|=m<\infty$, and hence 
by the Composite Pencil Theorem we see 
that $k^{\sharp}$ is relatively algebraically closed in $L^{\sharp}$.
Now defo$(f)=X_1^mX_2^2$, and hence the points at infinity of the curve
$f=0$ are $(X_0,X_1,X_2)=(0,1,0)$ and $(X_0,X_1,X_2)=(0,0,1)$.
Homogenizing $f$ we obtain
\begin{align*}
&(X_1-a_1X_0)\dots(X_1-a_mX_0)X_2^2\\
&+\gamma (X_1-b_1X_0)\dots(X_1-b_mX_0)X_2X_0+X_1X_0^{m+1}-zX_0^{m+2}.
\end{align*}
Putting $X_1=1$ we get 
\begin{align*}
f_1=&(1-a_1X_0)\dots(1-a_mX_0)X_2^2\\
&+\gamma (1-b_1X_0)\dots(1-b_mX_0)X_2X_0+X_0^{m+1}-zX_0^{m+2}
\end{align*}
and hence $(0,1,0)$ is a double point with two distinct tangents.
Putting $X_2=1$ we get 
\begin{align*}
f_2=&(X_1-a_1X_0)\dots(X_1-a_mX_0)\\
&+\gamma (X_1-b_1X_0)\dots(X_1-b_mX_0)X_0+X_1X_0^{m+1}-zX_0^{m+2}
\end{align*}
and hence $(0,0,1)$ is an $m$-fold point with $m$ distinct tangents.
Therefore $\tau(f)=m+2$, and hence by taking $(z+Z,k^{\sharp})$ for
$(z,k)$ we get $\tau(f^{\sharp})=m+2$.

\centerline{}

EXAMPLE 3. Assume $m\ge 2$. Let $\alpha_i$ and $\alpha_{ij}$
be the elements in $k$ such that
$$
a(X_1)=X_1^m+\sum_{1\le j\le m}\alpha_{j}X_1^{m-j}
$$ 
and
$$
a(X_1)/(X_1-a_i)=X_1^{m-1}+\sum_{1\le j\le m-1}\alpha_{ij}X_1^{m-1-j}
\quad\text{ for }\quad 1\le i\le m.
$$
Let $b(X_1)=X_1^m+\sum_{i\in\{1,2,m\}}\beta_i X_1^{m-i}$
where $(\beta_i)_{i\in\{1,2,m\}}$ are elements $k$ such that: 

(i) $b(a_i)\ne 0$ for $1\le i\le m$;

(ii) if $m>3$ then $\beta_1\ne\alpha_1$, and 
$\beta_2\ne\beta_1\alpha_{i1}-\alpha_{i1}^2+\alpha_{i2}$ for $1\le i\le m$;

(iii) if $m=3$ then $\beta_1\ne\alpha_1$, and
$\beta_2\ne\beta_1\alpha_{i1}-\alpha_{i1}^2+\alpha_{i2}+1$ for $1\le i\le m$;

(iv) if $m=2$ then $\beta_1\ne 1+\alpha_1$, and
$\beta_2\ne\beta_1\alpha_{i1}-\alpha_{i1}^2-\alpha_{i1}-a_i$.

\noindent
Let $f=a(X_1)X_2^2+b(X_1)X_2+X_1-z$ with $z\in k$. 
Then $f$ is irreducible in $R$ and
redset$(f)=\text{redset}(f)^*=\emptyset$. 
Moreover, $k^{\sharp}$ is relatively algebraically closed in $L^{\sharp}$
and we have $|\text{redset}(f)|=|\text{redset}(f)^*|=0$ with 
$\tau(f)=\tau(f^{\sharp})=m+2$.

\centerline{}

PROOF. Let $\mu=0$. By (i) we have gcd$(a(X_1),b(X_1))=1$ and hence
the proof is identical with the proof of Example 1
except we have to get a contradiction to (*) where $\gamma=1$ and
$c$ is any element of $k^*$.
So assume (*) and let $y=z+c$.
Then comparing degrees and coefficients of $X_1^m$ in the
LHS and RHS of (*) we see that $\eta=1$ and either
$(u(X_1),v(X_1))=(1,a(X_1))$ or
$(u(X_1),v(X_1))=(a(X_1)/a_i(X_1),a_i(X_1))$ for some $i\in\{1,\dots,m\}$.
Let us expand the said LHS and RHS and compare the coefficients of $X_1^{m-1}$
and $X_1^{m-2}$ in them. Then in case of $m>3$ we have
\begin{align*}
&(u(X_1),v(X_1))=(1,a(X_1))\\
&\Rightarrow
X_1^m+\sum_{j\in\{1,2,m\}}\beta_j X_1^{m-j}
=(X_1-y)+\left[X_1^m+\sum_{1\le j\le m}\alpha_j X_1^{m-j}\right]\\
&\Rightarrow
\beta_1=\alpha_1\\	
\end{align*}
and
\begin{align*}
&(u(X_1),v(X_1))=(a(X_1)/a_i(X_1),a_i(X_1))\\
&\Rightarrow
X_1^m+\sum_{j\in\{1,2,m\}}\beta_j X_1^{m-j}
=(X_1-y)\left[X_1^{m-1}+\sum_{1\le j\le m}\alpha_{ij} X_1^{m-1-j}\right]
+(X_1-a_i)\\
&\Rightarrow
\beta_1=\alpha_{i1}-y\text{ and }\beta_2=-y\alpha_{i1}+\alpha_{i2}\\	
&\Rightarrow
\beta_2=\beta_1\alpha_{i1}-\alpha_{i1}^2+\alpha_{i2}\\	
\end{align*}
and hence we get a contradiction by (ii).
Likewise in case of $m=3$ we have
\begin{align*}
&(u(X_1),v(X_1))=(1,a(X_1))\\
&\Rightarrow
X_1^3+\beta_1 X_1^2+\beta_2 X_1+\beta_3
=(X_1-y)+(X_1^3+\alpha_1 X_1^2+\alpha_2 X_1+\alpha_3)\\
&\Rightarrow
\beta_1=\alpha_1\\	
\end{align*}
and
\begin{align*}
&(u(X_1),v(X_1))=(a(X_1)/a_i(X_1),a_i(X_1))\\
&\Rightarrow
X_1^3+\beta_1 X_1^2+\beta_2 X_1+\beta_3
=(X_1-y)(X_1^2+\alpha_{i1} X_1+\alpha_{12})+(X_1-a_i)\\
&\Rightarrow
\beta_1=\alpha_{i1}-y\text{ and }\beta_2=-y\alpha_{i1}+\alpha_{i2}+1\\	
&\Rightarrow
\beta_2=\beta_1\alpha_{i1}-\alpha_{i1}^2+\alpha_{i2}+1\\	
\end{align*}
and hence we get a contradiction by (iii).
Finally in case of $m=2$ we have
\begin{align*}
&(u(X_1),v(X_1))=(1,a(X_1))\\
&\Rightarrow
X_1^2+\beta_1 X_1+\beta_2
=(X_1-y)+(X_1^2+\alpha_1 X_1+\alpha_2)\\
&\Rightarrow
\beta_1=1+\alpha_1\\	
\end{align*}
and
\begin{align*}
&(u(X_1),v(X_1))=(a(X_1)/a_i(X_1),a_i(X_1))\\
&\Rightarrow
X_1^2+\beta_1 X_1+\beta_2
=(X_1-y)(X_1+\alpha_{i1})+(X_1-a_i)\\
&\Rightarrow
\beta_1=\alpha_{i1}-y+1\text{ and }\beta_2=-y\alpha_{i1}-a_i\\	
&\Rightarrow
\beta_2=\beta_1\alpha_{i1}-\alpha_{i1}^2-\alpha_{i1}-a_i\\	
\end{align*}
and hence we get a contradiction by (iv).

\centerline{}

EXAMPLE 4. Let
$f=a(X_1)X_2^2+X_2+X_1^{m+1}-z$ with $z\in k$.
Then $f$ is irreducible in $R$
and redset$(f)=\text{redset}(f)^*=\emptyset$.
Moreover, $k^{\sharp}$ is relatively algebraically closed in $L^{\sharp}$
and we have $|\text{redset}(f)|=|\text{redset}(f)^*|=0$ with 
$\tau(f)=\tau(f^{\sharp})=m+1$.

\centerline{}

PROOF.  By Gauss Lemma we see that if $c\in\text{redset}(f)^*$ then 
for some $u(X_1)$, $v(X_1)$, $p(X_1)$, $q(X_1)$ in $R^{\times}$ we have
$$
f-c=[u(X_1)X_2+p(X_1)][v(X_1)X_2+q(X_1)]
$$
and equating coefficients of powers of $X_2$ we get
\begin{equation*}
u(X_1)v(X_1)=a(X_1)
\tag{$4_1$}
\end{equation*}
and
\begin{equation*}
p(X_1)q(X_1)=r(X_1)
\quad\text{ with }\quad
r(X_1)=X_1^{m+1}-z-c
\tag{$4_2$}
\end{equation*}
and
\begin{equation*}
u(X_1)q(X_1)+v(X_1)p(X_1)=1.
\tag{$4_3$}
\end{equation*}
For a while suppress the variable $X_1$. 
Note that then deg$(a)=m$ and deg$(r)=m+1$.
First suppose $m$ is even.  Then $m+1$ is odd. 
Multiplying $(4_3)$ by $u$ and using $(4_1)$ we get
\begin{equation*}
u^2q+ap=u.
\tag{$4'_4$}
\end{equation*}
By $(4_2)$ we know that deg$(pq)$ is odd and hence 
$(\text{deg}(q),\text{deg}(p))=$ (even, odd) or (odd, even). Therefore
$(\text{deg}(u^2q),\text{deg}(ap))=$ (even, odd) or (odd, even), and hence
deg$(u^2q+ap)=\text{max}(\text{deg}(u^2q),\text{deg}(ap))>\text{deg}(u)$
which contradicts $(4'_4)$.
Next suppose $m$ is odd.  Then $m+1$ is even. 
Multiplying $(4_3)$ by $p$ and using $(4_2)$ we get
\begin{equation*}
ur+vp^2=p.
\tag{$4''_4$}
\end{equation*}
By $(4_1)$ we know that deg$(uv)$ is odd and hence 
$(\text{deg}(u),\text{deg}(v))=$ (even, odd) or (odd, even). Therefore
$(\text{deg}(ur),\text{deg}(vp^2))=$ (even, odd) or (odd, even), and hence
deg$(ur+vp^2)=\text{max}(\text{deg}(ur),\text{deg}(vp^2))>\text{deg}(p)$
which contradicts $(4''_4)$.
Thus we have shown that redset$(f)^*=\emptyset$. 
Consequently $f$ is irreducible in $R$ and redset$(f)=\emptyset$. 
In particular $|\text{redset}(f)|=|\text{redset}(f)^*|=0<\infty$, and hence 
by the Composite Pencil Theorem we see 
that $k^{\sharp}$ is relatively algebraically closed in $L^{\sharp}$.
Now defo$(f)=X_1^mX_2^2$, and hence the points at infinity of the curve
$f=0$ are $(X_0,X_1,X_2)=(0,1,0)$ and $(X_0,X_1,X_2)=(0,0,1)$.
Homogenizing $f$ we obtain
$$
(X_1-a_1X_0)\dots(X_1-a_mX_0)X_2^2+X_2X_0^{m+1}+X_1^{m+1}X_0-zX_0^{m+2}.
$$
Putting $X_1=1$ we get 
$$
f_1=(1-a_1X_0)\dots(1-a_mX_0)X_2^2+X_2X_0^{m+1}+X_1^{m+1}X_0-zX_0^{m+2}
$$
and hence $(0,1,0)$ is a simple point. Putting $X_2=1$ we get 
$$
f_2=(X_1-a_1X_0)\dots(X_1-a_mX_0)+X_0+X_1^{m+1}X_0-zX_0^{m+2}
$$
and hence $(0,0,1)$ is an $m$-fold point with $m$ distinct tangents.
Therefore $\tau(f)=m+1$, and hence by taking $(z+Z,k^{\sharp})$ for
$(z,k)$ we get $\tau(f^{\sharp})=m+1$.

\centerline{}

EXAMPLE 5. Assume $m\ge 2$. Let
$f=X_2^m a(X_1/X_2)+z$ with $z\in k^{\times}$.
Then $f$ is irreducible in $R$ 
and redset$(f)=\text{redsset}(f)^*=\{z\}$.
Moreover, $k^{\sharp}$ is relatively algebraically closed in $L^{\sharp}$
and we have $|\text{redset}(f)|=|\text{redset}(f)^*|=1$ with 
$\tau(f)=\tau(f^{\sharp})=m$.

\centerline{}

PROOF. Obviously $z\in\text{redset}(f)$. Moreover, 
for any $c\ne z$, homogenizing $f-c$ and then putting $X_2=1$ we
get the polynomial $f_2=a(X_1)+(z-c)X_0^m$ which is clearly irreducible
and hence so in $f-c$. Also obviously $f-c=0$ has $m$ places at infinity.

\centerline{}

\centerline{\bf Section 6: More General Pencils}

Let us now study more general pencils of hypersurfaces $f-cw=0$ with $c$ 
varying over $k$, where $w\in R^{\times}$ is such that gcd$(f,w)=1$. 
Let $R^{\flat}=R^{\sharp}=k(Z)[X_1,\dots,X_m]$.

The assumption gcd$(f,w)=1$ says that the pencil
$(f-cw)_{c\in k}$ is without fixed components.
However, it may have base loci of dimension $< n-1$, i.e., 
there may be primes of height 
$>1$ in $R$ which contain $f$ and $w$ both, and nothing can be said about
the singularities of $f-cw=0$ at a base point.
To indicate these primes, 
for any $g_1,\dots,g_m$ in $R$ we let 
$\mathcal V(g_1,\dots,g_m)$ denote the variety defined by
$g_1=\dots=g_m=0$, i.e. we put
$\mathcal V(g_1,\dots,g_m)=\{P\in\text{spec}(R)\text{ with }
(g_1,\dots,g_m)R\subset P\}$; 
now $\mathcal V(f,w)$ is the set of all base loci of our pencil; 
for $n=2$ the set $\mathcal V(f,w)$ is finite and its members
are the base points of the pencil.
Likewise for any $g_1,\dots,g_m$ in $R^*$ we put 
$\mathcal V(g_1,\dots,g_m)^*=\{P\in\text{spec}(R^*)\text{ with }
(g_1,\dots,g_m)R^*\subset P\}$, 
and for any $g_1,\dots,g_m$ in $R^{\flat}$ we put 
$\mathcal V(g_1,\dots,g_m)^{\flat}=\{P\in\text{spec}(R^{\flat})\text{ with }
(g_1,\dots,g_m)R^{\flat}\subset P\}$. 

Let singset$(f,w)=\{c\in k:R_P/((f-cw)R_P)\text{ is nonregular for some }
P\in\text{spec}(R)\setminus\mathcal V(f,w)\text{ with }f-cw\in P\}$, and let
singset$(f,w)^*=\{c\in k^*:R^*_P/((f-cw)R^*_P)$
is nonregular for some
$P\in\text{spec}(R^*)\setminus\mathcal V(f,w)^*\text{ with }f-cw\in P\}$. 
Also let
redset$(f,w)=\{c\in k:f-cw=gh \text{ for some }g,h\text{ in }R\setminus k\}$,
and let redset$(f,w)^*=\{c\in k^*:f-cw=gh$
for some $g,h\text{ in }R^*\setminus k^*\}$.
Finally let
multset$(f,w)^*=\{c\in k^*:f-cw=gh^2\text{ for some }g\in (R^*)^{\times}
\text{ and }h\in R^*\setminus k^*\}$, and let
primset$(f,w)=\{c\in k:f-cw=gh^\mu$ for some 
$g\in k^{\times}\text{ and }h\in R\setminus k\text{ with integer }\mu>1\}$.

Sometimes we need to projectivise the pencil $f-cw=0$ by allowing $c$ to vary 
over $k\cup\{\infty\}$ and declaring that $f-\infty w$ means $w$. To take 
care of this we put redset$(f,w)_+=$ redset$(f,w)$ or 
redset$(f,w)\cup\{\infty\}$ 
according as we cannot or can write $w=gh$ with $g,h$ in $R\setminus k$.
Likewise we put redset$(f,w)^*_+=$ redset$(f,w)^*$ or 
redset$(f,w)^*\cup\{\infty\}$ according as we cannot or can write $w=gh$ 
with $g,h$ in $R^*\setminus k^*$. 
Similarly  we put multset$(f,w)^*_+=$ multset$(f,w)^*$ or 
multset$(f,w)^*\cup\{\infty\}$ according as we cannot or can write 
$w=gh^2$ for some $g\in (R^*)^{\times}$ and $h\in R^*\setminus k$.
Finally  we put primset$(f,w)_+=$ primset$(f,w)$ or 
primset$(f,w)\cup\{\infty\}$ according as we cannot or can write 
$w=gh^{\mu}$ with $g\in k^{\times}$ and $h\in R\setminus k$ with
integer $\mu>1$. Let
$$
d=\text{max}(\text{deg}(f),\text{deg}(w))
$$
and observe that, without assuming deg$(f)>0$ but assuming $d>0$,
the condition gcd$(f,w)=1$ implies that
\begin{equation*}
\text{deg}(f-cw)<d\text{ for at most one }c\in k^*\cup\{\infty\}
\tag{$'$}
\end{equation*}
and
\begin{equation*}
\text{gcd}(f-c_1w,f-c_2w)=1\text{ for all }c_1\ne c_2\text{ in } 
k^*\cup\{\infty\}.
\tag{$''$}
\end{equation*}

By the generic member of the pencil $(f-cw)_{c\in k}$ we mean the
hypersurface $(f,w)^{\flat}=0$ with $(f,w)^{\flat}=f-Zw\in R^{\flat}$.
Let singset$(f,w)^{\flat}=\{c\in k(Z):
R^{\flat}_P/(f-(Z+c)w)R^{\flat}_P)$ is nonregular for some
$P\in\text{spec}(R^{\flat})\setminus\mathcal V(f,w)^{\flat}$ with 
$f-(Z+c)w \in P\}$. 
Let redset$(f,w)^{\flat}=\{c\in k(Z)$ such that $f-(Z+c)w=gh$ for some
in $g,h$ in $R^{\flat}\setminus k(Z)\}$.

By Gauss Lemma $(f,w)^{\flat}$ is irreducible in $R^{\flat}$, i.e.,
$0\not\in\text{redset}(f,w)^{\flat}$.
Let $\phi^{\flat}: R^{\flat}\to R^{\flat}/((f,w)^{\flat}R^{\flat})$
be the residue class epimorphism. Clearly
$R\cap\text{ker}(\phi^{\flat})=\{0\}$ and
$\phi^{\flat}(f)=\phi^{\flat}(Z)\phi^{\flat}(w)$,
and hence there exists a unique isomorphism
$\psi^{\flat}:k(X_1,\dots,X_n)\to\text{QF}(R^{\flat}/((f,w)^{\flat}R^{\flat}))$
such that for all $r\in R$ we have $\psi^{\flat}(r)=\phi^{\flat}(r)$.
Thus the triple
$\phi^{\flat}(k(Z))\subset\phi^{\flat}(R^{\flat})
\subset\text{QF}(\phi^{\flat}(R^{\flat}))$ is isomorphic to the triple
$$
k^{\flat}=k(f/w)\subset A^{\flat}=k(f/w)[X_1,\dots,X_n]
\subset L^{\flat}=k(X_1,\dots,X_n)
$$
and hence we may regard the above three displayed sets as
the ground field, the affine coordinate ring, and the function field
of $(f,w)^{\flat}=0$.

\centerline{}

Let us reiterate that singset$(f,w)$, singset$(f,w)^*$, and
singset$(f,w)^{\flat}$ denote singularities outside the base points.
For any $P\in\text{spec}(A^{\flat})$ with $f-(Z+c)w \in P$ and
$(f,w)A^{\flat}\not\subset P$, 
upon letting $Q=R\cap P$, we see that the local ring $A^{\flat}_P$ is 
regular because it equals the localization of the regular local
ring $R_Q$ at its multiplicative subset $k[f/w]^{\times}$ or
$k[w/f]^{\times}$ depending on whether $w\not\in P$ or $f\not\in P$.
Thus we have the:

\centerline{}

GENERAL GENERIC SINGSET THEOREM.
$0\not\in\text{singset}(f,w)^{\flat}$.

\centerline{}

In view of the above isomorphism of triples, by the General Redset Theorem 
which we shall state and prove in a moment, we get the:

\centerline{}

GENERAL GENERIC REDSET THEOREM. If $k^{\flat}$ is relatively
algebraically closed in $L^{\flat}$, then redset$(f,w)^{\flat}$ is
finite.

\centerline{}

GENERAL REDSET THEOREM. If $f$ is irreducible in $R$ with
$\text{deg}(f)\ge\text{deg}(w)$ and $k$ is relatively
algebraically closed in $L$, then redset$(f,w)$ is finite. 

\centerline{}

PROOF. By the Lemma we can find a finite number of DVRs 
$V_1,\dots,V_t$ of $L/k$ such that $A[1/\phi(w)]\cap V_1\cap\dots\cap V_t=k$.
For every $z\in L^{\times}$ let 
$W_i(z)=\text{ord}_{V_i}(z)$, and let $W:L^{\times}\to\mathbb Z^t$ be the map
given by putting $W(z)=(W_1(z),\dots,W_t(z))$.
Let $G$ be the set of all $g\in R\setminus k$ such that $gh=f-cw$ for some
$h\in R\setminus k$ and $c\in k^{\times}$. Since the degree
of $g$ is clearly smaller than the degree of $f$, the set $G$ is contained
in a finite dimensional $k$-vector-subspace of $R$. Therefore for every
$i\in\{1,\dots,t\}$, the set $W_i(\phi(g))_{g\in G}$ is bounded from below.
Since $h$ also belongs to $G$ and clearly 
$W_i(\phi(g))=W_i(\phi(w))-W_i(\phi(h))$ with integer 
$W_i(\phi(w))$ which depends only on $i$ and $w$,
it follows that the set $W_i(\phi(g))_{g\in G}$ is also bounded from above.
Since, for $1\le i\le t$, 
the set $W_i(\phi(g))_{g\in G}$ is bounded from both sides, it follows that
$W(\phi(G))$ is a finite set. Also clearly $\phi(G)\subset U(A[1/\phi(w)])$.
Let $g_1h_1=f-c_1w$ and $g_2h_2=f-c_2w$ with $g_1,h_1,g_2,h_2$ in 
$R\setminus k$ and $c_1,c_2$ in $k^{\times}$ be such that
$W(\phi(g_1))=W(\phi(g_2))$.
Then $\phi(g_1)/\phi(g_2)\in A[1/\phi(w)]\cap V_1\cap\dots\cap V_t=k$ and  
hence $\phi(g_2)=c\phi(g_1)$ for some $c\in k^{\times}$. 
Consequently $g_2-cg_1$ is divisible by $f$ in $R$ and hence,
because $\text{deg}(g_2-cg_1)<\text{deg}(f)$, we must have $g_2=cg_1$.
Therefore, by subtracting the equation $g_2h_2=f-c_2w$ from the
equation $g_1h_1=f-c_1w$ we get $(c_2-c_1)w=g_1h_1-g_2h_2=g_1(h_1-ch_2)$ 
which implies that $(c_2-c_1)w$ 
is divisible in $R$ by the positive degree polynomial $g_1$ with
gcd$(w,g_1)=1$. Consequently we must have $c_2=c_1$. Therefore, because
the set $W(\phi(G))$ is finite, we conclude that redset$(f,w)$ is finite.

\centerline{}

The following version of L\"uroth's Theorem was first proved by
Igusa \cite{Igu}, and then it was deduced by Nagata \cite{Na2} as a
consequence of Abhyankar's paper \cite{A01}. 

\centerline{}

GENERAL REFINED L\"UROTH THEOREM. Assume that
$k^{\flat}$ is not relatively algebraically closed in $L^{\flat}$. Let 
$\widehat k^{\flat}$ be the algebraic closure of $k^{\flat}$ in 
$L^{\flat}$, and let $\nu=[\widehat k^{\flat}:k^{\flat}]$.
Then $\nu$ is an integer with $\nu>1$, and there exist
$\Gamma\in k[Z]^{\times}$, $\Omega\in k[Z]^{\times}$,
$\widehat f\in R\setminus k$, and $\widehat w\in R\setminus k$, with
gcd$(\Gamma,\Omega)=1$, max$(\text{deg}(\Gamma),\text{deg}(\Omega))=\nu$,
gcd$(\widehat f,\widehat w)=1$, 
max$(\text{deg}(\widehat f),\text{deg}(\widehat w))>0$, and
deg$(\widehat w)=\text{min}\{\text{deg}(\widehat f-\widehat c\widehat w):
\widehat c\in k^*\cup{\infty}\}$,
such that $\widehat k^{\flat}=k(\widehat f/\widehat w)$ and 
$f/w=\Gamma(\widehat f/\widehat w)/\Omega(\widehat f/\widehat w)$. 

\centerline{}

REMARK 6. In the above situation, every member of the pencil
$(f-cw)_{c\in k^*}$ consists of $\nu$ members (counted properly) of the
pencil $(\widehat f-\widehat c\widehat w)_{\widehat c\in k^*}$, and so we say 
that {\bf the pencil $(f-cw)_{c\in k^*}$ is composite with the pencil 
$(\widehat f-\widehat c\widehat w)_{\widehat c\in k^*}$}.
In greater detail, upon letting $\Lambda(Y,Z)=\Gamma(Z)-Y\Omega(Z)$
we get irreducible $\Lambda(Y,Z)$ in $k[Y,Z]$
of $Y$-degree $1$ and $Z$-degree $\nu$
such that $\Lambda(f/w,\widehat f/\widehat w)=0$.
For most $c\in k^*$ we have 
$\Lambda(c,Z)=\widehat c_0\prod_{1\le i\le\nu}(Z-\widehat c_i)$
where $\widehat c_0,\widehat c_1,\dots,\widehat c_{\nu}$ in $k^*$ with
$\widehat c_0\ne 0$ and so the ``object'' $f/w=c$ is the  union of
the ``objects'' $\widehat f/\widehat w=\widehat c_1,\dots, 
\widehat f/\widehat w=\widehat c_{\nu}$.
Alternatively, by letting
$\gamma(Z,T)=T^{\nu}\Gamma(Z/T)$,
$\omega(Z,T)=T^{\nu}\Omega(Z/T)$, and
$\lambda(Y,Z,T)=T^{\nu}\Lambda(Y,Z/T)$, we get polynomials which are 
homogeneous of degree $\nu$ in $(Z,T)$, and for which we have
$f/w=\gamma(\widehat f,\widehat w)/\omega(\widehat f,\widehat w)$ and
$\lambda(Y,Z,T)/\omega(Z,T)=(\gamma(Z,T)/\omega(Z,T))-Y$.
Now 
$\lambda(c,Z,T)=\widehat c_0\prod_{1\le i\le\nu}(Z-\widehat c_iT)$,
and so the hypersurface $f-cw=0$ is the union of the hypersurfaces
$\widehat f-\widehat c_1\widehat w=0,
\dots,\widehat f-\widehat c_{\nu}\widehat w=0$.
In the above phrase ``most $c\in k^*$'' we were referring to the tacit
assumption that $\lambda(c,Z,T)$ is not divisible by $T$; 
in the contrary case, if $T^{\mu}$ is the highest power of $T$ which
divides $\lambda(c,Z,T)$ then $\mu$ of the roots, say $c_1,\dots,c_{\mu}$,
have ``gone to infinity'' and the hypersurface $f-cw=0$ is composed of the
hypersurfaces 
$\widehat w=0,\dots,\widehat w=0,\widehat f-\widehat c_{\mu+1}\widehat w=0,
\dots,\widehat f-\widehat c_{\nu}\widehat w=0$ with
$\widehat w=0$ occurring $\mu$ times.
By factoring $\omega(Z,T)$ we get
$\omega(Z,T)=\widetilde c_0 T^{\epsilon}\prod_{\epsilon+1\le i\le \nu}
(Z-\widetilde c_i T)$ where
$\widetilde c_0,\widetilde c_{\epsilon+1}\dots,\widetilde c_{\nu}$
are elements in $k^*$ with $\widetilde c_0\ne 0$. The assumption 
gcd$(\Gamma,\Lambda)=1$ implies that $\widetilde c_i\ne \widehat c_j$
for $\epsilon+1\le i\ne \nu$ and $\mu+1\le j\le\nu$, 
and if $\epsilon>0$ then the phrase 
``most $c\in k^*$'' can be changed to the phrase ``all $c\in k^*$.''
In other words, for two different members of the ``old pencil''
$(f-cw)_{c\in k^*\cup\{\infty\}}$ the corresponding $\nu$-tuples in the 
``new pencil'' 
$(\widehat f-\widehat c\widehat w)_{\widehat c\in k^*\cup\{\infty\}}$
are disjoint.
Upon multiplying the polynomials $\Gamma$ and $\Omega$ by suitable constants
it can be arranged that $f=\gamma(\widehat f,\widehat w)$ and
$w=\omega(\widehat f,\widehat w)$; now for all $c\in k^*$ we
have $f-cw=\lambda(c,\widehat f,\widehat w)$.
The picture is completed by noting that every member of the projectivised 
pencil $(f-cw)_{c\in k^*\cup\{\infty\}}$ consists of $\nu$ members of the
projectivised pencil 
$(\widehat f-\widehat c\widehat w)_{\widehat c\in k^*\cup\{\infty\}}$, 
counted with multiplicities.

The above correspondence $(\widehat c_1,\dots,\widehat c_{\nu})\mapsto c$
can be elucidated by noting that the 
$(Z-\widehat c_1)\text{-adic},\dots,(Z-\widehat c_{\nu})\text{-adic}$
valuations of $k^*(Z)$ are the extensions of the $(Y-c)$-adic valuation of
$k^*(Y)$ with $Y=\Gamma(Z)/\Omega(Z)$. Here the $(Z-\infty)$-adic and
the $(Y-\infty)$-adic valuations are the negative degree functions, and the 
multiplicity of the root $\widehat c_i$ is the ramification exponent of
the $(Z-\widehat c_i)$-adic valuation. 
As said above, since gcd$(\Gamma,\Omega)=1$, for any two different $c$'s in
$k^*\cup\{\infty\}$ the two corresponding sets
$(\widehat c_1,\dots,\widehat c_{\nu})$ are disjoint.
By $(')$ we see that deg$(\widehat f-\widehat c\widehat w)<\widehat d$ for
at most one $\widehat c\in k^*\cup\{\infty\}$ and this can equal a
$\widehat c_i$ for at most one $c\in k^*\cup\{\infty\}$.

Since $[k^*(Z):k^*(Y)]=\nu>1$,
if the characteristic of $k$ is $0$ then, for some $c\in k^*$,
there exits $\widehat c\in k^*$ such that $(Z-\widehat c)$-adic 
valuation lies above the $(Y-c)$-adic valuation and has ramification
exponent greater than one, and hence $f-cw$ has a
nonconstant multiple factor in $R$, and therefore multset$(f,w)^*\ne\emptyset$.
To find $c$ explicitly, upon letting $\Theta(Y)$ be the $Z$-resultant of
$\Lambda(Y,Z)$ and its $Z$-derivative $\lambda_Z(Y,Z)$ we can show that
$\Theta(Y)$ has a root in $k^*$ which is not a root of the
coefficient of $Z^{\nu}$ in $\Lambda(Y,Z)$; now for $c$ we can take such
a root.

Without assuming any condition on the pair $(k^{\flat},L^{\flat})$, if 
$c\in k^*$ is such that $f-cw=gh^2$ with $g\in R^*\setminus\{0\}$ and
$h\in R^*\setminus k^*$ then $(f/w)-c=(g/w)h^2$ and hence taking the partials
of both sides relative to $X_i$ and multiplying by $w^2$ we get
$fw_{X_i}-wf_{X_i}=w^2[(g/w)_{X_i}h^2+2hh_{X_i}(g/w)]
=h[(gw_{X_i}-wg_{X_i})h+2gwh_{X_i}]$; therefore multset$(f,w)^*\ne\emptyset
\Rightarrow\text{gcd}(fw_{X_1}-wf_{X_1},\dots,fw_{X_n}-wf_{X_n})\ne 1$.

Again without assuming any condition on the pair $(k^{\flat},L^{\flat})$ we
see that $c\mapsto c-Z$ gives an injection of redset$(f,w)$ into
redset$(f,w)^{\flat}$, and hence
$|\text{redset}(f,w)|\le|\text{redset}(f,w)^{\flat}|$. 

Thus, in view of the General Generic Redset Theorem, the General
Refined L\"uroth Theorem, and the above two observations $(')$ and $('')$, 
we get the: 

\centerline{}

GENERAL COMPOSITE PENCIL THEOREM. We have the following:

(I) gcd$(fw_{X_1}-wf_{X_1},\dots,fw_{X_n}-wf_{X_n})=1
\Rightarrow\text{multset}(f,w)^*=\emptyset$.

(II) 
characteristic of $k$ is $0$ and 
multset$(f,w)^*=\emptyset\Rightarrow k^{\flat}$ is relatively algebraically 
closed in $L^{\flat}$.

(III)
$k^{\flat}$ is relatively algebraically closed in 
$L^{\flat}\Rightarrow$ redset$(f,w)$ and redset$(f,w)^{\flat}$
are finite with $|\text{redset}(f,w)|\le|\text{redset}(f,w)^{\flat}|$.

(IV)
$k^{\flat}\text{ is not relatively algebraically closed in }
L^{\flat}\Rightarrow|(k^*\cup\{\infty\})\setminus\text{redset}(f,w)^*_+|$
$\le 1$ and deg$(f-cw)<d$ for any
$c\in(k^*\cup\{\infty\})\setminus\text{redset}(f,w)^*_+$.

\centerline{}

By slightly changing the proof of the Singset Theorem we shall now
prove the:

\centerline{}

GENERAL SINGSET THEOREM.
If $k$ is of characteristic zero then singset$(f,w)$ is finite.

\centerline{}

PROOF. Let $I$ be the ideal in $R^*$ generated by the $n$ elements 
$fw_{X_i}-wf_{X_i}$ with $1\le i\le n$, and note that these elements when
divided by $w^2$ give us the partials $(f/w)_{X_i}$. For any
$P\in\text{spec}(R^*)$ with $I\subset P$ and $w\not\in P$, 
consider the residue class map
$\Phi_P:R^*\to R^*/P$. Since all the partials of $f/w$ when multiplied by
$w^2$ belong to $P$, it follows that
$D(\Phi_P(f)/\Phi_P(w))=0$ for every $\Phi_P(k^*)$-derivation of QF$(R^*/P)$. 
Therefore, since $k^*$ is of characteristic zero, we have 
$\Phi_P(f)/\Phi_P(w)=\Phi_P(\kappa(P))$ for a unique $\kappa(P)\in k^*$.
For any $c\in k^*$ we clearly have: $f-cw\in P\Leftrightarrow c=\kappa(P)$.
Let $P_1,\dots,P_s$ be the minimal primes of $I$ in $R^*$ which do not
contain $w$. 
Then for any $c\in k^*\setminus\{\kappa(P_1),\dots,\kappa(P_s)\}$ and
$Q\in\text{spec}(R^*)$ with $f-cw\in Q$ and $(f,w)R^*\not\subset Q$, we must
have $I\not\subset Q$. Since $k$ is of characteristic zero, it follows that
singset$(f,w)\subset\{\kappa(P_1),\dots,\kappa(P_s)\}$, and hence 
singset$(f,w)$ is finite.

\centerline{}

REMARK 7. To explain the ideas behind the proofs of the Singset and General
Singset Theorems, suppose $k$ to be of characteristic zero. If the equations
$f_{X_1}=\dots=f_{X_n}=0$ have only a finite number of common solutions 
$(a_{i1},\dots,a_{in})_{1\le i\le s}$ in $(k^*)^n$ then upon letting
$c_i=f(a_{i1},\dots,a_{in})$ we clearly get singset$(f)=\{c_1,\dots,c_s\}$
provided $k=k^*$ and hence singset$(f)\subset\{c_1,\dots,c_s\}$ without that
proviso.  Without assuming the common solutions to be finite,
this is generalized by letting $P_1,\dots,P_s$
to be the minimal primes of $(f_{X_1},\dots,f_{X_n})R^*$, i.e., upon letting
$\mathcal V(P_1)^*,\dots,\mathcal V(P_s)^*$ to be the irreducible components
of the affine variety $f_{X_1}=\dots=f_{X_n}=0$, then showing that $f$ is
constant on $\mathcal V(P_i)^*$ and letting $c_i$ be that constant value.
To apply this to get the singset of the more general pencil $f-cw$ outside
its base points we take only those irreducible components of the 
variety $fw_{X_1}-wf_{X_1}=\dots=fw_{X_n}-wf_{X_n}=0$
which are not contained in the hypersurface $w=0$.
Note that, since in these two theorems the characteristic is assumed to be 
zero, the ideal $I$ is a nonzero ideal. Observe that here
the characteristic assumption is essential as is shown by taking
$f=X_1^p+\dots+X_n^p$ with characteristic $p>0$ and $k=k^*$ and noting that
then singset$(f)=k$. However, no assumption on the characteristic is required 
in the Generic Singset Theorem and the General Generic Singset Theorem.

\centerline{}

In the Mixed Primset Theorem we showed that if the special pencil $f-c$ is 
noncomposite, i.e., if $k^{\sharp}$ is relatively algebraically closed in
$L^{\sharp}$, then primset$(f)$ is empty. However, if the general pencil
$f-cw$ is noncomposite, i.e., if $k^{\flat}$ is relatively algebraically 
closed in $L^{\flat}$, then primset$(f,w)$ may be nonempty.
More precisely, we prove the: 

\centerline{}

GENERAL MIXED PRIMSET THEOREM
Assume that $k$ is perfect and $k^{\flat}$ is relatively algebraically 
closed in $L^{\flat}$. 
For any $c\in k\cup\{\infty\}$ let $d(c)=\text{deg}(f-cw)$. For any 
$c\in\text{primset}(f,w)_+$ by definition we have
$f-cw=g(c)h(c)^{\mu(c)}$ with $g(c)\in k^{\times}$, $h(c)\in R\setminus k$,
and integer $\mu(c)>1$; we shall assume that the triple $(g(c),h(c),\mu(c))$
is so chosen that $\mu(c)$ is maximal. Call $c'\in k\cup\{\infty\}$
inseparable if $(f-c'w)_{X_i}=0$ for $1\le i\le n$; otherwise call $c'$
separable.  Then we have the following:

(I) primset$(f,w)_+$ has at most one inseparable member. 

(II) If primset$(f,w)_+$ has an inseparable $c'$ then it has
at most one $c\ne c'$ with $d(c)=d$.
If primset$(f,w)_+$ has an inseparable member then
$|\text{primset}(f,w)_+|\le 3$. 

(III) primset$(f,w)_+$ has at most three members $c$ with $d(c)=d$.

(IV) If primset$(f,w)_+$ has three distinct members $c_i$ with
$d(c_i)=d$ for $1\le i\le 3$ then 
$\{\mu(c_1),\mu(c_2),\mu(c_3)\}=\{2,3,5\}$.

(V) $|\text{primset}(f,w)|\le |\text{primset}(f,w)_+|\le 4$. 

\centerline{}

PROOF.  To prove (I) and (II), assuming primset$(f,w)_+$ to have an inseparable
member, and replacing $f$ and $w$ by suitable $k$-linear combinations of them,
we may suppose $\infty$ to be that inseparable member. Note that now
$k$ must be of characteristic $p>0$ and $w=(w')^p$ with $w'\in R\setminus k$.
If primset$(f,w)_+$ has another inseparable member then replacing it by a 
suitable $k$-linear combination of it and $w$ we may suppose the other 
inseparable member to be $0$; now $f=(f')^p$ with $f'\in R\setminus k$
and $k(f/w)\subset k(f'/w')$ with $[k(f'/w'):k(f/w)]=p$ contradicting the
noncompositness of our pencil. This proves (I). In particular, we must have
$f_{X_j}\ne 0$ for some $j\in\{1,\dots,n\}$. If primset$(f,w)_+$ has two 
separable members $c_1\ne c_2$ with $d(c_1)=d=d(c_2)$ then for $1\le i\le 2$
we get $f_{X_j}=(f-c_iw)_{X_j}=g(c_i)h(c_i)_{X_j}h(c_i)^{\mu(c_i)-1}$ and hence,
because clearly gcd$(h(c_1)^{\mu(c_1)-1},h(c_2)^{\mu(c_2)-1})=1$, we see that
$h(c_1)^{\mu(c_1)-1}h(c_2)^{\mu(c_2)-1}$ divides $f_{X_j}$ in $R$, and 
therefore 
\begin{align*}
d-1&\ge\text{deg}(f_{X_j})\\
&\ge\text{deg}(h(c_1)^{\mu(c_1)-1})+\text{deg}(h(c_2)^{\mu(c_2)-1})\\
&=d-(d/\mu(c_1))+d-(d/\mu(c_2))
\end{align*}
and dividing the first and the last expressions by $d$ and rearranging terms
suitably we obtain
$$
(1/\mu(c_1))+(1/\mu(c_2))\ge 1 + (1/d).
$$
Since our pencil is noncomposite, we must have gcd$(\mu(c_1),\mu(c_2))=1$,
and hence the left hand side is at most $5/6$ which is less
than the right hand side, giving a contradiction. This proves (II). 

To prove (III) to (V), let if possible 
$c_i$ in primset$(f,w)_+$ with $d(c_i)=d$ for $1\le i\le u$ be distinct
members where $u=3$ or $4$.
If $(f/w)_{X_j}=0$ for $1\le j\le n$ then $k$ must be of characteristic $p>0$
and we must have $f/w=(f'/w')^p$ for some $f',w'$ in $R^{\times}$ 
contradicting the noncompositness of our pencil. Therefore
$(f/w)_{X_j}\ne 0$ for some $j\in\{1,\dots,n\}$.
By the equation $f-c_iw=g(c_i)h(c_i)^{\mu(c_i)}$ we see that $f-c_iw$
and $(f-c_iw)_{X_j}$ are both divisible by $h(c_i)^{\mu(c_i)-1}$ in $R$.
Clearly $w^2(f/w)_{X_j}=fw_{X_j}-wf_{X_j}$ and hence if $c_i=\infty$ then
$w^2(f/w)_{X_j}$ is divisible by $h(c_i)^{\mu(c_i)-1}$ in $R$.
If $c_i\in k$ then
$w^2(f/w)_{X_j}=w^2((f-c_iw)/w)_{X_j}=(f-c_iw)w_{X_j}-w(f-c_iw)_{X_j}$ and 
hence again $w^2(f/w)_{X_j}$ is divisible by $h(c_i)^{\mu(c_i)-1}$ in $R$.
Clearly gcd$(h(c_i),h(c_{i'})=1$ for all $i\ne i'$, and hence
$w^2(f/w)_{X_j}$ is divisible by $\prod_{1\le i\le u}h(c_i)^{\mu(c_i)-1}$ 
in $R$. Therefore
\begin{align*}
2d-1&\ge\text{deg}(w^2(f/w)_{X_j})\\
&\ge\sum_{1\le u}\text{deg}(h(c_i)^{\mu(c_i)-1})\\
&=\sum_{1\le i\le u}[d-(d/\mu(c_i))]
\end{align*}
and dividing the first and the last expressions by $d$ and rearranging terms
suitably we obtain
$$
\sum_{1\le i\le u}(1/\mu(c_1))\ge (u-2) + (1/d).
$$
Again as above, two different $\mu_i$ must have gcd $1$, for otherwise, the
system is easily seen to be composite. In case of $u=4$, the LHS of the
above inequality has each term at most $1/2$ and at most one term equal
to $1/2$,  leading to a maximum value less than 
$2=(u-2)$, which is a contradiction in view of the RHS.
This proves (III) and hence also (V).
In case of $u=3$, we see that
$2,3,5$ are the only choices for mutually coprime $\mu(c_i)$ giving a sum 
bigger than $1=(u-2)$ for the LHS, which proves (IV).

\centerline{}

REMARK 8. Referring to (IV) above, for $n=2$ and $k$ any field of 
characteristic zero, there is indeed an example of a pencil with three such 
members in the primset. Namely we use the fact that the surface 
$X^2+Y^3+Z^5=0$ is rational. The specific parametrization is in Klein's 
Lectures on the Icosahedron \cite{Kle}.
To quote it explicitly, from Article 13 of Chapter I of \cite{Kle}, let
$$
f=H_1^2
\quad\text{ and }\quad
w=H_2^3
\quad\text{ and }\quad
v=1728H_3^5
$$
where
$$
H_1=X_1^{30}+X_2^{30}-10005X_1^{10}X_2^{10}(X_1^{10}+X_2^{10})
+522X_1^5X_2^5(X_1^{20}-X_2^{20})
$$
and
$$
H_2=-[X_1^{20}+X_2^{20}+494X_1^{10}X_2^{10}]
+228X_1^5X_2^5(X_1^{10}-X_2^{10})
$$
and
$$
H_3=X_1X_2[X_1^{10}-X_2^{10}+11X_1^5X_2^5].
$$
Let us put
$$
P=X_1^{10}-X_2^{10}\quad\text{ and }\quad Q=X_1^5X_2^5.
$$
Then
\begin{align*}
H_1
&=(X_1^{10}+X_2^{10})^3-10008X_1^{10}X_2^{10}(X_1^{10}+X_2^{10})
+522X_1^5X_2^5(X_1^{20}-X_2^{20})\\
&=(X_1^{10}+X_2^{10})(P^2-10004Q^2+522PQ)
\end{align*}
and hence
\begin{align*}
f&=(P^2+4Q^2)(P^2-10004Q^2+522PQ)^2\\
&=(P^2+4Q^2)[P^4+1044P^3Q+((522)^2-20008)P^2Q^2\\
&\quad\quad\qquad\qquad\qquad-(1044)(10004)PQ^3+(10004)^2Q^4]\\
&=P^6+1044P^5Q+[(522)^2-20004]P^4Q^2\\
&\quad-(1044)(10000)P^3Q^3+[(10004)^2+(1044)^2-4(20008)]P^2Q^4\\
&\quad-4(1044)(10004)PQ^5+(20008)^2Q^6.
\end{align*}
Also
$$
H_2=-P^2-496Q^2+228PQ
$$
and hence
\begin{align*}
w&=-[P^6+3(496)P^4Q^2+3(496)^2P^2Q^4+(496)^3Q^6]\\
&\quad+3(228)PQ[P^4+2(496)P^2Q^2+(496)^2Q^4]\\
&\quad-3(228)^2P^2Q^2[P^2+496Q^2]+(228)^3P^3Q^3\\
&=-P^6+684P^5Q-3[496+(228)^2]P^4Q^2\\
&\quad+[6(228)(496)+(228)^3]P^3Q^3-3[(496)^2+(228)^2(496)]P^2Q^4\\
&\quad+3(228)(496)^2PQ^5-(496)^3Q^6.
\end{align*}
Moreover
$$
H_3=X_1X_2(P+11Q)
$$
and hence
\begin{align*}
H_3^5
&=Q(P+11Q)^5\\
&=P^5Q+5(11)P^4Q^2+10(11)^2P^3Q^3+10(11)^3P^2Q^4+5(11)^4PQ^5+(11)^5Q^6.
\end{align*}
Adding the coefficients of $P^6,P^5Q,\dots,Q^6$ in $f$ and $w$, and
comparing the sums to the corresponding coefficients in $H_3^5$, we get 
$f+w=1728H_3^5=v$.

\centerline{}

Now we shall prove the:

\centerline{}

GENERAL MIXED REDSET THEOREM. Assume that $k$ is perfect and
$k^{\flat}$ is relatively algebraically closed in $L^{\flat}$. 
Note that now by the Lemma there exists a finite number
of DVRs $V_1,\dots,V_t$ of $L^{\flat}/k^{\flat}$ with
$A^{\flat}[1/w]\cap V_1 \cap \dots \cap V_t= k^{\flat}$ and, for any such $t$,
we have that $t$ is a positive integer and 
$U(A^{\flat}[1/w])/U(k^{\flat})$ is a finitely generated free
abelian group of rank $r\le t-1$.
By the General Mixed Primset Theorem, upon letting 
$\rho =|\text{primset}(f,w)|$, we get an integer $\rho$ with 
$0\le \rho\le 4$. We claim that $|\text{redset}(f,w)|\le r+\rho$. 

\centerline{}

PROOF. 
As in the proof of the Mixed Redset Theorem, suppose if possible that we have 
distinct elements $c_1,\dots,c_s$ in $k$ with integer
$s > r+\rho$ such that 
\begin{equation*}
f-c_iw =g_ih_i\quad\text{ where }\quad
g_i,h_i\text{ in }R\setminus k.
\tag{1*}
\end{equation*}
Let $\tau=s-\rho$. Then $\tau$ is an integer with $\tau>r$,
and upon relabelling $c_1,\dots,c_s$ we can arrange matters so that
\begin{equation*}
\quad\text{gcd}(g_i,h_i)=1 \quad \text{for}\quad 1 \le i \le\tau.
\tag{2*}
\end{equation*}
Since $\tau>r$, there exist integers $a_1,\dots,a_{\tau}$ at least one 
of which is nonzero, say $a_{\iota}\ne 0$, such that
\begin{equation*}
\prod_{1\le i\le \tau}g_i^{a_i}\in U(k(f/w)).
\tag{3*}
\end{equation*}
As in the Second Argument of the proof of the Mixed Redset Theorem, by (3*) we
get
\begin{equation*}
\prod_{1\le i\le\tau}g_i^{a_i}=\overline u
\prod_{1\le i\le \sigma}u_i(f/w)^{\alpha_i}
\quad\text{ with }\quad\overline u\in k^{\times}
\tag{4*}
\end{equation*}
where $\alpha_1,\dots,\alpha_\sigma$ are nonzero integers and
$u_1(Z),\dots,u_\sigma(Z)$ are pairwise coprime monic
members of $k[Z]\setminus k$.
For any $u,v$ in $k[Z]\setminus k$ with gcd$(u,v)=1$ 
we have $w^{\text{deg}(u)}u(f/w),w^{\text{deg}(v)}v(f/w)$ in $R\setminus k$ 
with gcd$(w^{\text{deg}(u)}u(f/w),w^{\text{deg}(v)}v(f/w))=1$. 
Therefore, by looking at divisibility by a nonconstant irreducible factor
of $g_i$ in $R$, in view of (1*), (2*), and (4*), we see that for each 
$i\in\{1,\dots,\tau\}$ with $a_i\ne 0$ there is some 
$\theta(i)\in\{1,\dots,\sigma\}$ with $u_{\theta(i)}(Z)=Z-c_i$ and
$\alpha_{\theta(i)}=a_i$. In particular we get
$u_{\theta(\iota)}(Z)=Z-c_\iota$ and $\alpha_{\theta(\iota)}=a_\iota$.
Again in view of (1*), (2*), and (4*), 
by looking at divisibility by a nonconstant irreducible factor of $h_{\iota}$ 
in $R$, we get a contradiction. 

\centerline{}

REMARK 9. In the $n=2$ case of the above theorem we can take $V_1,\dots,V_t$
to be valuation rings of the places at infinity of the curve $f^{\flat}=0$
together with the valuation rings of its places at finite distance centered at
points where it meets the curve $w=0$.

\centerline{}

REMARK 10. To draw a deduction chart for the various incarnations of the
Redset Theorem, as already observed, we have:

\begin{align*}
\text{Lemma }\Rightarrow \text{Redset Theorem }
&\Rightarrow\text{ Generic Redset Theorem}\\
&\Rightarrow\text{ Composite Pencil Theorem}
\end{align*}
and
\begin{align*}
\text{Lemma }\Rightarrow \text{General Redset Theorem }
&\Rightarrow\text{ General Generic Redset Theorem}\\
&\Rightarrow\text{ General Composite Pencil Theorem}
\end{align*}
and we have:
$$
\text{Lemma }\Rightarrow\text{ Mixed Redset Theorem}
$$
and
$$
\text{Lemma }\Rightarrow\text{General Mixed Redset Theorem}
$$
where A $\Rightarrow$ B means B can be deduced from A.
Also clearly:
$$
\text{Composite Pencil Theorem }
\Rightarrow \text{Redset Theorem when }k=k^*
$$
and
$$
\text{General Composite Pencil Theorem }
\Rightarrow \text{General Redset Theorem when }k=k^*.
$$

\centerline{}

REMARK 11.
The condition $\text{deg}(f)\ge\text{deg}(w)$ in the General Redset Theorem
is necessary as well as reasonable. It is reasonable because if 
$\text{deg}(f)<\text{deg}(w)$ and we pass to the projective $n$-space
by homogenizing $f$ and $w$, then $f$ becomes reducible by acquiring the
hyperplane at infinity counted $\text{deg}(w)-\text{deg}(f)$ times, and we 
have not yet obtained any control over the rest of the pencil.
To see that it is necessary, take $w$ to be a polynomial in $f$ of degree
$\nu>1$ and note that then, for every $c\in k^*$, $f-cw$ is clearly a product
of $\nu$ members of $R^*\setminus k^*$.

\centerline{}

\centerline{}

It is time to whet the appetite of the reader by raising a few questions. 

\centerline{}

QUESTION 1.
As the bounds found in the Mixed Redset Theorem were sharpened in
Remark 5 and Examples 1 to 5, can you do similar things to sharpen the bounds
found in the General Mixed Redset Theorem?

\centerline{}

QUESTION 2.
In the Mixed Redset Theorem it was shown that
$k^{\sharp}$ relatively algebraically closed in $L^{\sharp}\Rightarrow
|\text{redset}(f)|\le\text{rank}(U(A^{\sharp})/U(k^{\sharp}))$.
Can you similarly obtain a good bound in the context of the Redset Theorem
by showing that, assuming $f$ to be irreducible in $R$,
$k$ relatively algebraically closed in 
$L\Rightarrow |\text{redset}(f)|\le\text{rank}(U(A)/U(k))$?
Can you relate the conditions 
``$k$ relatively algebraically closed in $L$'' and
``$k^{\sharp}$ relatively algebraically closed in $L^{\sharp}$''? 
In case of $n=2$, how far can you relate the numbers of places at infinity of 
the curves $f=0$ and $f^{\sharp}=0$?

\centerline{}

QUESTION 3.
Can you, in a manner similar to Question 2, relate the hypotheses and conclusions
of the General Redset Theorem and the General Mixed Redset Theorem?

\centerline{}

QUESTION 4. Let the unique component set of $(f,w)$ be defined by putting
uniset$(f,w)=\{c\in k:f-cw=gh^\mu$ for some $g\in k^{\times}$ and irreducible
$h\in R\setminus k$ with integer $\mu>1\}$, and let us put
uniset$(f,w)_+=$ uniset$(f,w)$ or 
uniset$(f,w)\cup\{\infty\}$ according as we cannot or can write 
$w=gh^{\mu}$ with $g\in k^{\times}$ and irreducible $h\in R\setminus k$ with
integer $\mu>1$. 
Note that clearly uniset$(f,w)\subset\text{primset}(f,w)$ and
uniset$(f,w)_+\subset\text{primset}(f,w)_+$. Also note that the Generic Mixed 
Redset Theorem and its proof remain valid if we let 
$\rho=|\text{uniset}(f,w)|$. In this manner we get a possibly stronger from
of the Generic Mixed Redset Theorem. Can you show that this is indeed a 
stronger from? In other words, can you show that $|\text{uniset}(f,w)_+|\le 3$
or $|\text{uniset}(f,w)_+|\le 2$? Hint: vis-a-vis Remark 8, study all 
parametrizations of the Klein surface. You may also redo Question 3 by
putting in the above stronger form of the Generic Mixed Redset Theorem.

\centerline{}

\centerline{\bf Section 7: Complements of Hypersurfaces}

We shall now give several criteria, i.e., necessary conditions, for the ring
$R[1/f]$ to be isomorphic to the ring $R[1/f']$ where $f'$ is another member
of $R\setminus k$. In geometric terms, this amounts to giving criteria for
the complements of the hypersurfaces $f=0$ and $f'=0$ in the affine $n$-space
to be biregularly equivalent to each other. Note that $R[1/f]$ may be viewed as
the affine coordinate ring of the hypersurface $f(X_1,\dots,X_n)X_{n+1}-1=0$
in the affine $(n+1)$-space. Also note that if the hypersurfaces $f=0$ and 
$f'=0$ are automorphic, i.e., if there exists an automorphism of $R$ which 
sends $f$ to $f'$, then the rings $R[1/f]$ and $R[1/f']$ must be
isomorphic, and hence our criteria also provide necessary conditions for the
hypersurfaces $f=0$ and $f'=0$ to be automorphic. 

For any $c\in k$ let us write
$f-c=f_{c,0}f_{c,1}^{e(c,1)}\dots f_{c,q(c)}^{e(c,q(c))}$ with 
$f_{c,0}\in k^{\times}$, integers $1\le e(c,1)\le\dots\le e(c,q(c))$, and
pairwise coprime irreducible elements $f_{c,1},\dots,f_{c,q(c)}$ in
$R\setminus k$. Note that then 
$c\in\text{redset}(f)\Leftrightarrow e(c,1)+\dots+e(c,q(c))>1$. 
Let $e(c)$ denote the sequence $(e(c,1),\dots,e(c,q(c)))$. 
Let refset$(f)=$ the refined redset of $f$ be defined to be the family of
sequences $e(c)_{c\in\text{redset}(f)}$.
Let $f'-c=f'_{c,0}f'{}^{e'(c,1)}_{c,1}\dots f'{}^{e'(c,q'(c))}_{c,q'(c)}$ and
$e'(c)=(e'(c,1),\dots,e'(c,q'(c)))$ be the corresponding factorization and the
corresponding sequence for $f'-c$.
Let us write refset$(f)=\text{refset}(f')$ to mean that there exists a
bijection $\theta:\text{redset}(f)\to\text{redset}(f')$ such that for all
$c\in\text{redset}(f)$ we have $e(c)=e'(\theta(c))$, i.e., $q(c)=q'(c)$ and
$e(c,i)=e'(\theta(c),i)$ for $1\le i\le q(c)$.
Thus, geometrically speaking, refset$(f)=\text{refset}(f')$ means there 
is a multiplicities preserving bijection between 
the reducible members of the pencils $f-c=0$ and $f'-c=0$.

Let us write $R[1/f]\approx R[1/f']$ to mean that there exists a ring
isomorphism of $R[1/f]$ onto $R[1/f']$, and let us write
$R[1/f]\approx_k R[1/f']$ to mean that there exists a ring $k$-isomorphism of 
$R[1/f]$ onto $R[1/f']$. As our first criterion we want to show that
(i) if $R[1/f]\approx_k R[1/f']$ then $q(0)=q'(0)$, and (ii) if $f$ and $f'$
are irreducible in $R$ and $R[1/f]\approx_k R[1/f']$ then  
refset$(f)=\text{refset}(f')$. But since it takes only a little more effort,
we might as well prove the somewhat stronger:

\centerline{}

REFSET CRITERION. 
(i) If $R[1/f]\approx R[1/f']$ then $q(0)=q'(0)$, and (ii) if $f$ and $f'$
are irreducible in $R$ and $R[1/f]\approx R[1/f']$ then  
refset$(f)=\text{refset}(f')$. 

\centerline{}

Since $R$ is a UFD and clearly $k=\{0\}\cup U(R)$ is relatively algebraically
closed in the quotient field $k(X_1,\dots,X_n)$ of $R$, by applying the
following Corollary of Lemma with $(k',A')=(k,R[1/f])$, this follows by taking
$S=R$ and changing capital letters to lower case letter1/s in the even stronger:

\centerline{}

UFD CRITERION. 
Let $S$ be a UFD, and let $F,F'$ be in $S\setminus U(S)$.
For any $C\in\{0\}\cup U(S)$ let us write
$F-C=F_{C,0}F_{C,1}^{E(C,1)}\dots F_{C,Q(C)}^{E(C,Q(C))}$ with 
$F_{C,0}\in U(S)$, integers $1\le E(C,1)\le\dots\le E(C,Q(C))$, and
pairwise coprime irreducible elements $F_{C,1},\dots,F_{C,Q(C)}$ in
$S\setminus U(S)$. 
Let $F'-C=F'_{C,0}F'{}^{E'(C,1)}_{C,1}\dots F'{}^{E'(C,Q'(C))}_{C,Q'(C)}$
be the corresponding factorization of $F'-C$.
Assume that there exists a ring isomorphism $\Delta:S[1/F]\to S[1/F']$
such that $\Delta(U(S))=U(S)$. Then (i) $Q(0)=Q'(0)$. 
Moreover (ii) if $F$ and $F'$ are irreducible in $S$ then there 
exists a bijection $\Theta:U(S)\to U(S)$ such that for all $C\in U(S)$ we have
$Q(C)=Q'(\Theta(C))$ and $E(C,i)=E'(\Theta(C),i)$ for $1\le i\le Q(C)$.

\centerline{}

PROOF.
Clearly $S[1/F]$ is UFD in which the irreducible nonunits are associates of
the irreducible
nonunits of $S$ except $F_{0,1},\dots,F_{0,Q(0)}$ which have become units.
Also clearly $U(S[1/F])/U(S)$ is a free abelian group of rank $Q(0)$
generated by $F_{0,1},\dots,F_{0,Q(0)}$. Similarly
$U(S[1/F'])/U(S)$ is a free abelian group of rank $Q'(0)$
generated by $F'_{0,1},\dots,F'_{0,Q(0)}$. For any ring isomorphism
$\Delta:S[1/F]\to S[1/F']$ we obviously have $\Delta(U(S[1/F])=U(S[1/F'])$.
Therefore, since we are assuming $\Delta(U(S))=U(S)$, we get an induced
isomorphism $U(S[1/F])/U(S)\to U(S[1/F'])/U(S)$ and hence $Q(0)=Q'(0)$.
This proves (i).
Now assume that $F$ and $F'$ are irreducible in $S$. Then because of the said
induced isomorphism we must have either $\Delta(F)=\alpha F'$ with 
$\alpha\in U(S)$ or $\Delta(F)=\alpha/F'$ with $\alpha\in U(S)$. 
If $\Delta(F)=\alpha F'$ then, letting $\Theta:U(S)\to U(S)$ be the bijection
given by $C\mapsto\Delta(C)/\alpha$, we have
\begin{align*}
\Delta(F_{C,0})\prod_{1\le i\le Q(C)}\Delta(F^{E(C,i)}_{C,i})
&=\Delta(F-C)=\alpha(F'-\Theta(C))\\
&=\alpha F'_{\Theta(C),0}
\prod_{1\le i\le Q'(\Theta(C))}F'{}^{E'(\Theta(C),i)}_{\Theta(C),i}
\end{align*}
and hence we get $Q(C)=Q'(\Theta(C))$ and $E(C,i)=E'(\Theta(C),i)$
for $1\le i\le Q(C)$.
If $\Delta(F)=\alpha/F'$ then, letting $\Theta:U(S)\to U(S)$ be the bijection
given by $C\mapsto \alpha/\Delta(C)$, we have
\begin{align*}
\Delta(F_{C,0})\prod_{1\le i\le Q(C)}\Delta(F^{E(C,i)}_{C,i})
&=\Delta(F-C)=(-\Delta(C)/F')(F'-\Theta(C))\\
&=(-\Delta(C)/F')F'_{\Theta(C),0}
\prod_{1\le i\le Q'(\Theta(C))}F'{}^{E'(\Theta(C),i)}_{\Theta(C),i}
\end{align*}
and hence again we get $Q(C)=Q'(\Theta(C))$ and $E(C,i)=E'(\Theta(C),i)$
for $1\le i\le Q(C)$.
In the last two sentences we have used the obvious facts that the elements
$\Delta(F_{C,1}),\dots,\Delta(F_{C,Q(C)})$ are pairwise coprime irreducible
nonunits in $S[1/F']$ and so are the elements
$F'_{\Theta(C),1},\dots,F'_{\Theta(C),Q'(C)}$;
moreover, the elements $\Delta(F_{C,0})$, $\alpha F'_{\Theta(C),0}$, and 
$(-\Delta(C)/F')F'_{\Theta(C),0}$ are units in $S[1/F']$.

\centerline{}

COROLLARY OF LEMMA. 
Let $A'$ be an affine domain over a field $k'$, let $k''$ be the algebraic
closure of $k'$ in $L'=\text{QF}(A')$, and let $\overline k=k''\cap A'$. Then 
$\overline k$ can be located only using the ring structure of $A'$ by noting 
that it is only subfield of $A'$ which equals the intersection of $A'$ with 
a finite number of DVRs of $L'$.

\centerline{}

PROOF.
Clearly $\overline k$ is a subfield of $A'$ and by the Lemma there exists a 
finite number of DVRs $V_1,\dots,V_t$ of $L'/k'$ with 
$A'\cap V_1\cap\dots\cap V_t=\overline k$. Let $V'_1,\dots,V'_{t'}$ be any
finite number of DVRs of $L'$ such that $A'\cap V'_1\cap\dots\cap V'_{t'}$ is a
subfield $\widetilde k$ of $A'$. Let $\widehat k=\overline k\cap\widetilde k$.
We want to show that then $\widehat k=\overline k$. Clearly $\widehat k$ is
a subfield of $\overline k$ and we have 
$\widehat k=\overline k\cap V'_1\cap\dots\cap V'_{t'}$.
For any $i$ the intersection $\widehat k\cap V'_i$ is either a DVR of 
$\overline k$ or equals $\overline k$. Let $i_1<\dots<i_s$ be those values of 
$i$ for which the said intersection is a DVR, and let 
$W_j=\overline k\cap V'_{i_j}$. Now $W_1,\dots,W_s$ are a finite number of
DVRs of the field $\overline k$ and their intersection is the field 
$\widehat k$. Now it is a well-known fact that if $W_1,\dots,W_s$ are any
finite number of valuation rings with a common quotient field $\overline k$
then $\overline k$ is the quotient field of their intersection, and hence in 
our situation we must have $\widehat k=\overline k$. A proof of the said fact
can be found in Theorem (11,11) on page 38 of \cite{Na1}. 
Since our valuations are real, we can
deduce the fact from the Approximation Theorem (see Theorem 18 on page 45 of
volume II of \cite{ZaS}) thus.
Since the statement is obvious when $s=0$, suppose $s>0$. Let $M(W_i)$ denote
the maximal ideal of $W_i$. By the said theorem
we can find $u\in W_1\setminus M(W_1)$ such that $u\in M(W_i)$ for all $i>1$.
Now $u\in W_1\cap\dots\cap W_s$ and, since the valuations are real, for
any $v\in W_1$ we can find an integer $m>0$ such that 
$vu^m\in W_1\cap\dots\cap W_s$ and hence 
$v\in\text{QF}(W_1\cap\dots\cap W_s)$. Therefore 
$\text{QF}(W_1\cap\dots\cap W_s)=\text{QF}(W_1)=\overline k$.

\centerline{}

EXAMPLE 6.
Take $f=X_1^{a_1}\dots X_n^{a_n}-1$ and $f'=X_1^{a'_1}\dots X_n^{a'_n}-1$
where $1\le a_1\le\dots\le a_n$ and $1\le a'_1\le\dots\le a'_n$
are integers such that gcd$(a_1,a_j)=1$ and gcd$(a'_1,a'_{j'})=1$ for some
$j\in\{2,\dots,n\}$ and $j'\in\{2,\dots,n\}$, and $a_i\ne a_{i'}$ for some
$i\in\{1,\dots,n\}$ and $i'\in\{1,\dots,n\}$.
Then clearly redset$(f)=\{-1\}$ and redset$(f')=\{-1\}$ with
refset$(f)=\{(a_1,\dots,a_n)\}$ and refset$(f')=\{(a'_1,\dots,a'_n)\}$.
Consequently refset$(f)\ne\text{refset}(f')$. 
It follows that for all $b,b'$ in $k^{\times}$ we have
$\text{refset}(bf)\ne\text{refset}(b'f')$ and therefore  by the Refset 
Criterion we get $R[1/(bf)]\not\approx R[1/(b'f')]$.
Hence in particular no automorphism of $R$ can send the ideal $fR$ to the
ideal $f'R$.

Now suppose $n=2$. 
Let $\epsilon:R\to k[T,T^{-1}]$ be the $k$-homomorphism given by 
$(X_1,X_2)\mapsto(T^{-a_2},T^{a_1})$. Then clearly ker$(\epsilon)=fR$. 
Since gcd$(a_1,a_2)=1$, for some integers $b_1,b_2$ we have $b_1a_1-b_2a_2=1$; 
by adding a positive multiple of $a_2$ to $b_1$ and adding the same
multiple of $a_1$ to $b_2$ we can arrange $b_1$ to be positive and then 
automatically $b_2$ will also become positive; clearly 
$\epsilon(X_1^{b_2}X_2^{b_1})=T$. Also we have $\beta_1a_1-\beta_2a_2=-1$
where $\beta_1=-b_1$ and $\beta_2=-b_2$; 
by adding a positive multiple of $a_2$ to $\beta_1$ and adding the same
multiple of $a_1$ to $\beta_2$ we can arrange $\beta_1$ to be positive and then 
automatically $\beta_2$ will also become positive; clearly 
$\epsilon(X_1^{\beta_2}X_2^{\beta_1})=T^{-1}$. Thus $\epsilon$ is surjective.
Similarly we find a surjective $k$-homomorphism
$\epsilon':R\to k[T,T^{-1}]$ with ker$(\epsilon')=f'R$. However, as shown 
above, $\epsilon$ and $\epsilon'$ do not differ from each other by an 
automorphism of $R$. Thus $f=0$ and $f'=0$ are ``hyperbolas'' which
are not automorphic to each other. 

\centerline{}

Next we come to the:

\centerline{}

NONRULED CRITERION. Recall that an irreducible $\overline f\in R\setminus k$ 
is said to be ruled if, after identifying $k$ with a subfield of 
$R/(\overline fR)$, there exits a subfield $\widetilde L$ of
$\overline L=\text{QF}(R/(\overline f))$ such that 
$k\subset\widetilde L\subset\widetilde L(t)=\overline L$ where $t$ is 
transcendental over $\widetilde L$. Let us relabel 
$f_{0,1},\dots,f_{0,q(0)}$ so that $f_{0,i}$ is nonruled or ruled according as
$1\le i\le r$ or $r+1\le i\le q(0)$, and let us relabel
$f'_{0,1},\dots,f_{0,q'(0)}$ so that $f'_{0,i}$ is nonruled or ruled 
according as $1\le i\le r'$ or $r'+1\le i\le q'(0)$; note that now the 
exponent sequences $e(0,1),\dots,e(0,q(0))$ and $e'(0,1),\dots,e'(0,q'(0))$ 
need not be nondecreasing. Assume that $R[1/f]\approx_k R[1/f']$. Then
$r=r'$ and, after relabelling $f_{0,1},\dots,f_{0,r}$ suitably, we have
$\text{QF}(R/(f_{0,i}R))\approx_k\text{QF}(R/(f'_{0,i}R))$ for $1\le i\le r$.

\centerline{}

PROOF. Recall that $\mathfrak V(R)=\{R_P:P\in\text{spec}(R)\}$. For 
$0\le i\le n$ let
$R_i=k[X_0/X_i,\dots,X_n/X_i]$ where $X_0=1$, and let the projective 
$n$-space $\mathfrak P$ over $k$ be defined to be the nonsingular projective 
model $\cup_{0\le i\le n}\mathfrak V(R_i)$ of $K=k(X_1,\dots,X_n)$ over $k$.
For the language of models see Abhyankar's books \cite{A02, A06, A09}.
In particular note that for any valuation ring $V$ of $K/k$, i.e., valuation
ring with quotient field $K$ and having $k$ as a subfield, the center of $V$
on $\mathfrak P$ is the unique member of $\mathfrak P$ dominated by $V$; we 
identify $k$ with a subfield of $V/M(V)$ where $M(V)$ is the maximal ideal
of $V$, and we let restrdeg$_k(V)$ denote the transcendence degree of 
$V/M(V)$ over $k$. By a prime divisor of $K/k$ we mean a DVR $V$ of $K/k$
such that restrdeg$_k(V)=n-1$; we call the prime divisor ruled if
there exits a subfield $\widetilde L$ of $V/M(V)$ such that 
$k\subset\widetilde L\subset\widetilde L(t)=V/M(V)$ where $t$ is 
transcendental over $\widetilde L$. Upon letting $V_i=R_{f_{0,i}R}$ and
$V'_i=R_{f'_{0,i}R}$ we get nonruled prime divisors $V_1,\dots,V_r$ and
$V'_1,\dots,V'_{r'}$ of $K/k$.  Note that 
$\mathfrak V(R[1/f])\subset\mathfrak P$ and 
$\mathfrak V(R[1/f'])\subset\mathfrak P$.
By Proposition 3 on page 336 of \cite{A01} we see that $V_1,\dots,V_r$ are
exactly all the nonruled prime divisors of $K/k$ whose center on $\mathfrak P$
is not in $\mathfrak V(R[1/f])$, and $V'_1,\dots,V'_{r'}$ are 
exactly all the nonruled prime divisors of $K/k$ whose center on $\mathfrak P$
is not in $\mathfrak V(R[1/f'])$. Since $R[1/f]\approx_k R[1/f']$, it follows
that $r=r'$ and, after relabelling $f_{0,1},\dots,f_{0,r}$ suitably, we have
$\text{QF}(R/(f_{0,i}R))\approx_k\text{QF}(R/(f'_{0,i}R))$ for $1\le i\le r$.

\centerline{}

EXAMPLE 7.
Take $f=X_1^a+\dots+X_n^a-1$ and $f'=X_1^{a'}+\dots+X_n^{a'}-1$ where
$a$ and $a'$ are positive integers with $a>a'$ and $a>n>1$; assume that $a$ and
$a'$ are
nondivisible by the characteristic of $k$ in case the latter is nonzero.
The polynomials $f$ and $f'$ are clearly irreducible in $R$.
It is expected that: 

(7*) $f$ is nonruled and $\text{QF}(R/(fR))\not\approx_k\text{QF}(R/(f'R))$.
\newline
Assuming (7*), by the Nonruled Criterion we would get
$R[1/f]\not\approx_k R[1/f']$.

\centerline{}

QUESTION 5. Can you prove the above statement (7*)? 
Hint for $n=2$: show that the genus of the nonsingular plane curve 
$f=0$ is $(a-1)(a-2)/2$ and that of $f'=0$ is $(a'-1)(a'-2)/2$; see \cite{A06}.
Hint for $n>2$: show that the arithmetic genus of the nonsingular
hypersurface $f=0$ is $(a-1)\dots(a-n)/n!$ and that of $f'=0$ is 
$(a'-1)\dots(a'-n)/n!$; see \cite{A06}.
Further hint for $n>2$: by using the domination part of the
desingularization theory, show that the arithmetic genus of a 
nonsingular projective model is a birational invariant; see \cite{A09}.

\centerline{}

To avoid the problem of showing that $f$ is nonruled in the above Question 5,
let us establish the:

\centerline{}

GENERIC CRITERION. Assume that $f$ and $f'$ are irreducible in $R$. Also assume
that $R[1/f]\approx_k R[1/f']$. Then the generic members of the pencils 
$R^{\sharp}((f-Z)R^{\sharp})\approx_{k(Z)}R^{\sharp}((f'-Z)R^{\sharp})$
where we recall that $Z$ is an indeterminate over $R$ and
$R^{\sharp}=k(Z)[X_1,\dots,X_n]$. Hence in particular the said generic members
are birationally equivalent to each other, i.e.,
$\text{QF}(R^{\sharp}((f-Z)R^{\sharp}))
\approx_{k(Z)}\text{QF}(R^{\sharp}((f'-Z)R^{\sharp}))$.

\centerline{}

PROOF. Given a $k$-isomorphism $\delta:R[1/f]\to R[1/f']$, 
as in the proof of the UFD Criterion, for some $\alpha\in k^{\times}$
we have either $\delta(f)=\alpha f'$ or $\delta(f)=\alpha/f'$.
Consequently we have either $\delta(k[f]^{\times})=\delta(k[f']^{\times})$ or
$\delta(k[f]^{\times})=\delta(k[1/f']^{\times})$ respectively.
Recall that $k^{\sharp}=k(f)$ and $A^{\sharp}=k(f)[X_1,\dots,X_n]$. Let
$k'{}^{\sharp}=k(f')$ and $A'{}^{\sharp}=k(f')[X_1,\dots,X_n]$.
Clearly $A^{\sharp}$ is the localization of $R[1/f]$ at $k[f]^{\sharp}$,
and $A'{}^{\sharp}$ is the localization of $R[1/f']$ at $k[f']^{\sharp}$
as well as at $k[1/f']^{\times}$. Therefore $\delta$ has a unique extension
to an isomorphism $\delta^{\sharp}:A^{\sharp}\to A'{}^{\sharp}$.
As noted before, the pair $(k(Z),R^{\sharp}((f-Z)R^{\sharp}))$ is isomorphic
to the pair $(k^{\sharp},A^{\sharp})$, and similarly the pair
$(k(Z),R^{\sharp}((f'-Z)R^{\sharp}))$ is isomorphic
to the pair $(k'{}^{\sharp},A'{}^{\sharp})$. It follows that
$R^{\sharp}((f-Z)R^{\sharp})\approx_{k(Z)}R^{\sharp}((f'-Z)R^{\sharp})$ and
hence $\text{QF}(R^{\sharp}((f-Z)R^{\sharp}))
\approx_{k(Z)}\text{QF}(R^{\sharp}((f'-Z)R^{\sharp}))$.

\centerline{}

QUESTION 6. In connection with Example 7, can you show that
$\text{QF}(R^{\sharp}((f-Z)R^{\sharp}))
\not\approx_{k(Z)}\text{QF}(R^{\sharp}((f'-Z)R^{\sharp}))$?
Assuming this, by the Generic Criterion, we would get
$R[1/f]\not\approx_k R[1/f']$.

\centerline{}

REMARK 12. In geometric terms, considering the birational equivalence of the 
complements of two hypersurfaces $f=0$ and $f'=0$ in the affine $n$-space, and 
relating it to their refined redsets as well as to the birational
equivalence of their irreducible components and 
the biregular equivalence of the generic members of their associated pencils
$f-c=0$ and $f'-c=0$, we have the following:

(I) The first part of the Redset Criterion says that if the affine 
complements of two hypersurfaces are biregularly equivalent then they have 
the same number of irreducible components.

(II) The second part of the Redset Theorem says that if the affine 
complements of two irreducible hypersurfaces are biregularly equivalent then 
they have the same refined redsets.

(III) the Nonruled Criterion says that if the affine complements of two 
hypersurfaces are biregularly equivalent then 
their nonruled irreducible components are birationally equivalent.

(IV) the Generic Criterion says that if the affine complements of two 
irreducible hypersurfaces are biregularly equivalent then 
the generic members of their associated pencils are biregularly equivalent and 
hence birationally equivalent.

\centerline{}

QUESTION 7. In the complex case, can you topologize the conclusions of
Remark 12 by replacing biregular equivalence by some kind of topological type?
Can you somehow relate this to the topology of complements as exemplified
by the work of Zariski, Fan, Teicher, and others as discussed in \cite{GaT}?  
Can you also tie it to Abhyankar's algebraization of the tame fundamental 
groups of complements as described in \cite{A03}?

\centerline{}

QUESTION 8. We have already noted the fact that the hypersurfaces $f=0$ and
$f'=0$ being automorphic implies the biregular equivalence of their
complements. Can you exploit this fact to link-up the results of Remark 12 with
the epimorphism theorems and problems discussed by Abhyankar 
in his Kyoto Notes \cite{A04}

\centerline{}

\centerline{\bf Section 8: Redset of a Plane Curve and Zariski's Lemma}

We shall now show how, in case of characteristic zero, 
the finiteness of the redset of a hypersurface can be deduced
from that of a plane curve via the famous Lemma 5 of Zariski's Bertini II
paper \cite{Za1}. 
It may be noted that, Abhyankar \cite{A05, A06, A08} reduced 
the Galois case of the Jacobian
Problem to the birational case by means of Zariski's Lemma and then settled
the birational case by using Zariski's Main Theorem for which reference may
be made to \cite{A09}.

\centerline{}

REMARK 13.
As pointed out in Remark 1, for $n=2$, the Lemma next to the Redset Theorem
follows by taking $V_1,\dots,V_t$ to be the valuation rings of the places at
infinity of the irreducible plane curve $f=0$, and then the argument in the
proof of the Redset Theorem shows that if $k$ is relatively algebraically
closed in $L$ then redset$(f)$ is finite. To see how, for
characteristic zero, the $n=2$ case of the Redset Theorem implies the $n>2$ 
case we can proceed thus. Let $\iota$ be the smallest positive integer $\le n$
such that $f\in k[X_1,\dots,X_{\iota}]$. If $\iota=1$ then by doing nothing,
and if $\iota>1$ then by applying a $k-$automorphism to $R$ of the form
$X_j \mapsto X_j$ for all $j\in\{1,\dots,n\}\setminus\{\iota\}$ and 
$X_\iota \mapsto X_\iota + X^l_1$ with $l>$ twice the total degree of $f$, 
we can arrange matters
so that $f$ is essentially monic in $X_1$, i.e., for some $e'\in k^{\times}$ 
and integer $e>0$
we have $f=e'X^e_1 + $ terms of $X_1$-degree $<e$.
For any positive integer $i\leq n$ let $R_i=k_i[X_1,\dots,X_i]$ with 
$k_i=k(X_{i+1},\dots,X_n)$.  Now clearly $f \in R_\iota\setminus k_\iota$ 
and by Gauss Lemma, $f$ is irreducible in $R_\iota$.  
Let us identify $k_\iota$ with a subfield
of $L_\iota=\text{QF}(R_\iota/(fR_\iota))$.  
If $\iota=2$ and $k_\iota$ were relatively
algebraically closed in $L_\iota$ then by the 
$n=2$ case of the Redset Theorem we would see that 
$\{c\in\widetilde k:f-c=gh$ with $g,h$ in $R_\iota\setminus k_\iota\}$ 
is finite and hence so is redset$(f)$.  In the next Example $8$ we shall
show that $f$ irreducible in 
$R_\iota$ does not imply $k_\iota$ relatively 
algebraically closed in $L_\iota$.  
However, in case of characteristic zero, by 
taking $(L',Z_1,\dots,Z_s)=(L,\phi(X_2),\dots,\phi(X_n))$ in the Second 
Version of Zariski's Lemma given below, we can arrange $k_2$ to be relatively 
algebraically closed in $L_2$.

\centerline{}

EXAMPLE 8.  The geometric significance of the condition that the ground
field $k$ of the irreducible polynomial $f$ be relatively algebraically
closed in its function field $L$ is due to the well-known fact that,
in case of characteristic zero, it is equivalent to assuming $f$ to be
absolutely irreducible, i.e., irreducible in $R^*$.  To illustrate
this take $n=2$ and $f=X^2_1+X^u_2$ with integer $u>0$.  Assume that
$-1$ is not a square in $k$; for instance $k$ could be the field of
real numbers.  Let $R_1=k_1[X_1]$ and $R^*_1=k^*_1[X_1]$ where $k^*_1$
is an algebraic closure of $k_1=k(X_2)$. Now clearly $f$ is always irreducible 
in $R_1$, but it is irreducible in $R^*_1\Leftrightarrow u$ is odd.
Let $\phi_1:R_1 \to R_1/(fR_1)$ be the canonical epimorphism, and 
let us identify $k_1$ with a subfield of $L_1=\text{QF}(R_1/(fR_1))$.
Note that now $L_1=k(\overline X_1,X_2)$ with $\overline X_1 =\phi_1(X_1)$.
If $u=2v$ is even then $(\overline X_1 /X^v_2)^2=-1$ and hence
$k_1$ is not relatively algebraically closed in $L_1$.  If $u=2v+1$ is odd
then $(\overline X_1 /X^v_2)^2=-X_2$ and hence 
$L_1=k(\overline X_1/X^v_2)$ and therefore $k_1$ is relatively algebraically 
closed in $L_1$.

\centerline{}

REMARK 14.  Before coming to Zariski's Lemma, let us note that he calls
maximally algebraic (m.a. for short) what we have called relatively
algebraically closed., i.e., a field is m.a. in an overfield if every element
of the overfield which is algebraic over the field belongs to the field;
likewise he calls a field quasi-maximally algebraic (q.m.a. for short)
in an overfield if every element of the overfield which is separable
algebraic over the field belongs to the field; this is sometimes
called relatively separably closed.  Recall that an overfield is said to
be regular over (or a regular extension of) a field if the overfield is
a finitely generated extension of the field such that the field is m.a. in
the overfield and the overfield is separably generated over the field.

The well-known fact about an irreducible polynomial mentioned in Example 8,
also applies to any prime ideal $P$ in $R$, after identifying $k$ with a
subfield of the function field QF$(R/P)$ of the variety $\mathcal V(P)$.
Namely, $P$ is absolutely prime $\Leftrightarrow$ QF$(R/P)$ is regular
over $k$.  Recall that $P$ is absolutely prime means $PR^*$ is prime, and
note that then: $P$ is absolutely prime $\Leftrightarrow P\overline k[X_1,
\dots ,X_n]$ is prime for every field extension $\overline k$ of $k
\Leftrightarrow P\overline k[X_1,\dots,X_n]$ is prime for every algebraic
field extension $\overline k$ of $k \Leftrightarrow P \overline k[X_1,
\dots, X_n]$ is prime for every finite algebraic field extension 
$\overline k$ of $K$.  Let us call an ideal $Q$ in $R$ quasiprime if
it is primary;  note that this implies $Q \neq R$;  also note $fR$
is quasiprime $\Leftrightarrow f=gh^\mu$ for some $g \in k^{\times}$
and irreducible $h\in R \setminus k$ with integer $\mu > 0$.  Let us call $P$ 
absolutely quasiprime to mean that $PR^*$ is quasiprime, and note that then: 
$P$ is absolutely quasiprime $\Leftrightarrow P\overline k[X_1,\dots,X_n]$ is
quasiprime for every field extension 
$\overline k$ of $k\Leftrightarrow P\overline k[X_1,\dots,X_n]$ is quasiprime 
for every algebraic field extension 
$\overline k$ of $k\Leftrightarrow P \overline k[X_1,\dots,X_n]$ is quasiprime 
for every finite algebraic field extension $\overline k$ of $k$. 
As a well-known variation of the above well-known fact we have that $P$ is
absolutely quasiprime $\Leftrightarrow k$ is q.m.a. in QF$(R/P)$. Proofs of
all these assertions can be found in \cite{ZaS}.

Out of the following four versions of Zariski's Lemma, the first two
constitute Lemma 5 of \cite{Za1}, the third is Proposition I.6.1 of
\cite{Za3}, and the fourth is Theorem 2.4 of \cite{Mat} or Proposition 9.31 
of \cite{FrJ}.

\centerline{}

ZARISKI'S LEMMA. For any finitely generated field extension $L'$ of $k$
we have the following.

FIRST VERSION. Assume that $k$ is of characteristic zero and m.a. in $L'$. Let
elements $Z_1,Z_2$ in $L'$ be algebraically independent over $k$. Then for all
except a finite number of $c$ in $k$ we have that $k(Z_1+cZ_2)$ is m.a. 
in $L'$.

SECOND VERSION. Assume that $k$ is of characteristic zero and m.a. in $L'$. 
Let elements $Z_1,\dots,Z_s$ in $L'$, with $s>1$, be algebraically independent 
over $k$. Then by applying a $k$-linear automorphism to $k[Z_1,\dots,Z_s]$
it can be arranged that $k(Z_1,\dots,Z_{s-1})$ is m.a. in $L'$.
In other words, there exists a nonsingular $n\times n$ matrix $C=(C_{ij}$
over $k$ such that upon letting $Z_i^C=\sum_{1\le j\le s}C_{ij}Z_j$
we have that $k(Z_1^C,\dots,Z_{s-1}^C)$ is m.a. in $L'$. Moreover, the
constants $C_{ij}$ are nonspecial in the sense that for every $\Lambda\in k$
there is $H_i(\Lambda)\subset k$ with $|k\setminus H_i(\Lambda)|<\infty$
for $1\le i<s$ such that: if $C_{ij}=1$ for $1\le i\le s$, $C_{ij}=0$ for
$1<i+1<j\le s$, $C_{12}\in H_1(1)$, and $C_{i.i+1}\in H_i(C_{i-1,i})$ for
$1<i<s$. then $k(Z_1^C,\dots,Z_{s-1}^C)$ is m.a. in $L'$. 

THIRD VERSION. Assume that $k$ is q.m.a. in $L'$. 
Let $z_1,\dots,z_r$ be a finite number of elements in $L'$ such that
trdeg$_k k(z_1,\dots,z_r)>1$. Let $u_1,\dots,u_r$ be indeterminates over $L'$.
Then $k(u_1,\dots,u_r,u_1z_1+\dots+u_rz_r)$ is q.m.a. in $L'(u_1,\dots,u_r)$.

FOURTH VERSION. Assume that $L'$ is regular over $k$. Let
elements $Z_1,Z_2$ in $L'$ be algebraically independent over $k$, and assume
that $D(Z_2)\ne 0$ for some derivation $D$ of $L'/k$. Then for all
except a finite number of $c$ in $k$ we have that $L'$ is regular over
$k(Z_1+cZ_2)$. 

\centerline{}

\centerline{\bf Section 9: Singset of a Plane Curve
and the Zeuthen-Segre Invariant}

We have found bounds on redset and
primset for the special pencil $(f-c)_{c\in k}$ but in case of the singset
we have only stated that it is finite under appropriate conditions.
Actually, from the proof of the Singset Theorem it does follow that, in case of
characteristic zero, $|\text{singset}(f)|\le$ the number of irreducible
components of the variety of partials $\mathcal V(f_{X_1},\dots,f_{X_n})^*$.

Assuming $n=2$ with $(X_1,X_2)=(X,Y)$, and $k$ is of characteristic zero with
$k=$ its algebraic closure $k^*$, it is possible to give a more quantitative
estimate of the singset using the rank $\rho(f)$ of $f$ as introduced in
Section 11 of \cite{AbA}.  As we shall see, this rank $\rho(f)$
is related to the Zeuthen-Segre invariant.
For the convenience of the reader let us review
the definition of $\rho(f)$.

Let $g,h$ in $R$.  For $Q=(u,v)$ in the affine plane $\mathcal A=k^2$ we 
define the intersection multiplicity $I(g,h;Q)$ to be the $k$-vector-space 
dimension of $S/(g,h)S$ where $S$ is the localization of $R$ at the maximal 
ideal $(X-u,Y-v)R$.  Note that: if $(g,h)R\not\subset(X-u,Y-v)R$ then 
$I(g,h;Q)=0$; if $(g,h)R\subset qR$ for some $q\in R\setminus k$ with 
$q(u,v)=0$ then $I(g,h;Q)=\infty$; otherwise $I(g,h;Q)$ is a positive 
integer.  We define further intersection multiplicities by putting
$$
I(g,h;\mathcal A) = \sum_{Q\in \mathcal A}I(g,h;Q) 
\quad\text{and}\quad
I(g,h;f)=\sum_{\{Q\in \mathcal A:f(u,v)=0\}}I(g,h;Q)
$$
and
$$
I(g,h;\mathcal A\setminus f)=\sum_{\{Q \in \mathcal A:f(u,v)\neq 0\}}I(g,h;Q)
$$
with the usual conventions about infinity, and we note that these are 
nonnegative integers or $\infty$.  We also put
$$
\widehat I(g,h;\mathcal A)=\text{max}_{\mu \in k}I(g,h-\mu;\mathcal A)
$$
and we note that $\widehat I(g,h;\mathcal A)$ is a nonnegative integer or 
$\infty$, and: $\widehat I(g,h;\mathcal A)
=\infty\Leftrightarrow\text{gcd}(g,h-c)\ne 1$ for some $c\in k$. Next we put
$$
\alpha(g,h;\mathcal A)=\{\lambda \in k:I(g,h-\lambda; \mathcal A)
<\widehat I(g,h;\mathcal A)\}\
$$
and
$$
\beta(g,h;\mathcal A)=\sum_{0\neq\lambda \in \alpha(g,h;\mathcal A)}
[\widehat I(g,h;\mathcal A)-I(g,h-\lambda; \mathcal A)]
$$
and we note that then $\beta(g,h;\mathcal A)$ is a nonnegative integer or
$\infty$.  Finally we define
$$
\rho(f)=I(f_X,f_Y;\mathcal A\setminus f)+\beta(f_Y,f;\mathcal A)
$$
and we note that this is a nonnegative integer or $\infty$.

By augmenting $\mathcal A$ by points on the line at infinity we get the
projective plane $\mathcal P$ over $k$. For any 
$Q\in\mathcal P\setminus\mathcal A$ we define $I(g,h;Q)$ in an obvious manner
and we put
$$
I(g,h;\mathcal P)=\sum_{Q\in\mathcal P}I(g,h;Q).
$$
For any $Q=(u,v)\in \mathcal A$ we put
$$
\overline\chi(f;Q)=(\text{number of branches of $f$ at $Q$})-1
$$
and we note that if $f(u,v)=0$ then $f(X+u,Y+v)$ is a product of 
$\overline\chi(f;Q)+1$ irreducible nonconstant power series in
$k[[X,Y]]$ in case $f(u,v)=0$, and if $f(u,v)\ne 0$ then
$\overline\chi(f;Q)=-1$. We let
$$
\overline\chi(f;\mathcal A)
=\sum_{Q=(u,v)\in \mathcal A\text{ with }f(u,v)=0}\overline\chi(f;Q)
$$
and we note that this is a nonnegative integer. We define the
integer $\overline\chi(f;\mathcal P)\ge\overline\chi(f;\mathcal A)$ in an 
analogous manner and we put 
$$
\overline\chi(f;\infty)
=\overline\chi(f;\mathcal P)-\overline\chi(f;\mathcal A).
$$

If $f$ is irreducible then by $\gamma(f)$ we denote its genus, and 
we note that by the genus formula for $f$, given any $g\in R$ with
$g-c\not\in fR$ for all $c\in k$, we have
$$
2\gamma(f)-2=\text{deg}(d\phi(g))=\sum_{V\in\mathfrak R(f,\mathcal P)}
\text{ord}_V(d(\phi(g))
$$
where deg$(\phi(g))$ is the degree of the divisor of the differential
of $\phi(g)$ in the function field $L/k$ of $f$, with canonical
epimorphism $\phi:R\to A=R/(fR)$ and $L=\text{QF}(A)$, and where
$\mathfrak R(f,\mathcal P)$ is the set of all DVRs of $L/k$; we also put
$\mathfrak R(f,\mathcal A)=\{V\in\mathfrak R(f,\mathcal P): A\subset V\}$ and
$\mathfrak R(f,\infty)
=\mathfrak R(f,\mathcal P)\setminus \mathfrak R(f,\mathcal A)$;
moreover, for any $Q\in\mathcal P$ we let $\mathfrak R(f,Q)$ denote the set of
$V\in\mathfrak R(f,\mathcal P)$ having center $Q$ on $f$, and we note that
then $\mathfrak R(f,Q)$ is a finite set which is empty if and only if
``$Q$ does not lie on $f$.''

In the general case, by writing $f=f_1\dots f_s$ with irreducible 
$f_1,\dots,f_s$, we generalize the definition of $\gamma(f)$ by putting
$$
\gamma(f)=1+\sum_{1\le i\le s}(\gamma(f_i)-1).
$$
Relabelling $f_1,\dots,f_s$ suitably we can arrange that the first $r$ of them
are pairwise nonassociates, and every $f_i$ is an associate of $f_j$ for some
$j\le r$; now we put rad$(f)=f_0f_1\dots f_r$ where $0\ne f_0\in k$ is chosen
so that the coefficients of the highest lexicographic terms of $f$ and rad$(f)$
are equal; the lexicographic order is such that 
deg$(X^aY^b)\ge\text{deg}(X^{a'}Y^{b'})\Leftrightarrow$ either $b=b'$ and 
$a\ge a'$ or $b>b'$. Now we define the algebraic rank of $f$ by putting
$$
\rho_a(f)=2\gamma(\text{rad}(f))+\overline\chi(\text{rad}(f);\mathcal P).
$$
If $k=\mathbb C$ then by $\rho_t(f)$ we denote the rank of the first 
homology group of $f$, i.e., of the point-set 
$\{(u,v)\in\mathbb C^2:f(u,v)=0\}$.

Consider the condition:
\begin{equation*}
\text{$f$ is $Y$-monic of $Y$-degree $N>0$,} 
\tag{*}
\end{equation*}
i.e., deg$_{X,Y}(f-f_0Y^N)<N$ with $f_0\in k^{\times}$ and integer $N>0$. 
Also consider the conditions: 
\begin{equation*}
\text{gcd$(f_Y,f-c)=1$ for all $c\in k$,}
\tag{1*}
\end{equation*}
and
\begin{equation*}
\text{gcd$(f_Y,f-c)=1$ for all $c\in k$, and $f$ is irreducible,}
\tag{2*}
\end{equation*}
and 
\begin{equation*}
\text{gcd$(f_Y,f-c)=1$ for all $c\in k$, $f$ is irreducible, 
and $k=\mathbb C$,}
\tag{3*}
\end{equation*}
and
\begin{equation*}
\text{rad}(f)=f.
\tag{4*}
\end{equation*}
In (11.2) of \cite{AbA} it is shown that:
\begin{equation*}
(\text{*})+(\text{1*})\Rightarrow\rho(f)=(1-N)+I(f,f_Y;\mathcal A)-I(f_X,f_Y;f)
\tag{9.1}
\end{equation*}
where all the terms are integers.  In (11.5) of \cite{AbA} it is shown that:
\begin{equation*}
(\text{*})+(\text{2*})\Rightarrow\rho(f)=\rho_a(f)
\tag{9.2}
\end{equation*}
and just after (11.5) it is asserted that:
\begin{equation*}
(\text{*})+(\text{3*})\Rightarrow\rho(f)=\rho_t(f).
\tag{9.3}
\end{equation*}
In a moment we shall generalize (9.1) by showing that:
\begin{equation*}
(\text{*})+(\text{4*})
\Rightarrow\rho_a(f)=(1-N)+I(f,f_Y;\mathcal A)-I(f_X,f_Y;f)
\tag{9.4}
\end{equation*}
where obviously all the terms are integers. From this we shall deduce that:
\begin{equation*}
(\text{*})
\Rightarrow\rho_a(f)=(1-N)+\text{deg}_Y[f]
+I(f,f_Y/\widehat{[f]};\mathcal A)-I(f_X,f_Y/\widehat{[f]};f)
\tag{9.5}
\end{equation*}
with all terms integers, where
$$
[f]=\text{gcd}(f,f_X,f_Y)
\quad\text{ and }\quad
\widehat{[f]}=\text{gcd}(f_X,f_Y)
$$
with the gcds made unique by requiring them to be $Y$-monic.
While proving (9.5) we shall also show that:
\begin{equation*}
(*)\Rightarrow\begin{cases}
f=[f]\text{rad}(f)\\
\text{and }\widehat{[f]}=[f]\widetilde f
\text{ where }\widetilde f\in R\text{ with }
\mathcal V(f,\widetilde f)=\emptyset\\
\text{and }\widehat{[f]}=\prod_{c\in\text{multset}(f)^*}[f-c]\\
\text{and hence }|\text{multset}(f)^*|\le\text{deg}_Y\widehat{[f]}.
\end{cases}
\tag{9.6}
\end{equation*}
As a consequence of (9.5) we shall show that:
\begin{equation*}
\begin{cases}
\text{there exists a unique integer $\rho_{\pi}(f)$ together with}\\
\text{a unique finite subset defset$(f)$ of $k$ such that}\\
\text{$\rho_a(f-c)=\rho_\pi(f)$ 
for all $c\in k\setminus\text{defset}(f)$ and}\\
\text{$\rho_a(f-c)\ne\rho_\pi(f)$ for all $c\in\text{defset}(f)$.}\\
\end{cases}
\tag{9.7}
\end{equation*}
In reference to (9.7) we put
$$
\rho_\pi(f)=\text{the pencil-rank of the pencil }(f-c)_{c\in k}
$$
and
$$
\text{defset}(f)=\text{the deficiency set of }f.
$$
 From (9.7) we shall deduce that:
\begin{equation*}
(\text{*})\Rightarrow\rho_\pi(f)=(1-N)
+\widehat I(f,f_Y/\widehat{[f]};\mathcal A).
\tag{9.8}
\end{equation*}
 From (9.8) we shall deduce that:
\begin{equation*}
(\text{*})\Rightarrow
\rho_\pi(f)\ge\rho_a(f-c)\text{ for all }c\in k\setminus\text{multset}(f)^*.
\tag{9.9}
\end{equation*}
 From (9.8) we shall deduce the jungian formula for the pencil-rank saying that:
\begin{equation*}
\rho_\pi(f)=1-|\mathcal V_\infty(f)|
+\sum_{c\in\text{defset}(f)}[\rho_\pi(f)-\rho_a(f-c)]
\tag{9.10}
\end{equation*}
with the points at infinity of $f$ defined by
$$
\mathcal V_\infty(f)=\text{set of height-one members of }\mathcal V(f^+)
$$
where $f^+$ is the degree form of $f$ consisting of its highest degree terms.
We are using the adjective jungian in view of the fundamental contribution of
Jung \cite{Jun} to the theory of algebraic rank.
The Zeuthen-Segre invariant of $f$ is the integer $\zeta(f)$ defined by putting
$$
\zeta(f)=-|\mathcal V_\infty(f)|-\rho_\pi(f)
+\sum_{c\in\text{defset}(f)}[\rho_\pi(f)-\rho_a(f-c)]
$$
and from (9.10) it immediately follows that:
\begin{equation*}
\zeta(f)=-1.
\tag{9.11}
\end{equation*}

\centerline{}

Finally from (9.11) we shall deduce the:

\centerline{}

(9.12) DEFSET THEOREM. If (*) then 
$\text{singset}(f)\setminus\text{multset}(f)^*\subset\text{defset}(f)$, and
$|\text{defset}(f)|\le 1+\rho_a(f)+\text{deg}_Y\widehat{[f]}$ with
$|\text{singset}(f)|\le 1+\rho_a(f)+2\text{deg}_Y\widehat{[f]}$.

\centerline{}

Before turning to the proof of items (9.4) to (9.12), let us establish 
some common

\centerline{}

NOTATION AND CALCULATION.
Given any $f'\in R^{\times}$, write
$f'=\nz f'_1\dots f'_{s'}$ where $f'_1,\dots,f'_{s'}$ are
irreducible members of $R\setminus k$ and $\nz\in k^{\times}$, 
and let $\phi'_i:R\to A'_i=R/(f'_iR)$
be the canonical epimorphism and identify $k$ with a subfield of
$L'_i=\text{QF}(A'_i)$. Let
$Z_i=\{V\in\mathfrak R(f'_i,\infty):\text{ord}_V\phi'_i(f)\ge 0\}$ and
$P_i=\{V\in\mathfrak R(f'_i,\infty):\text{ord}_V\phi'_i(f)<0\}$, 
where the letters $Z$ and $P$ are meant to suggest zeros and poles,
and note that these are obviously finite sets. 
Let
$$
Z(f,f')=\sum_{1\le i\le s'}\sum_{V\in Z_i}\text{ord}_V\phi'_i(f)
\quad\text{ and }\quad
P(f,f')=\sum_{1\le i\le s'}\sum_{V\in P_i}\text{ord}_V\phi'_i(f)
$$
and note that $Z(f,f')$ is a nonnegative integer or $\infty$ according as
gcd$(f,f')=1$ or gcd$(f,f')\ne 1$, and $P(f,f')$ is always a nonpositive
integer.  Clearly for each $V\in Z_i$ there is a unique $c_i(V)\in k$ 
such that ord$_V\phi'_i(f-c_i(V))>0$.
Let $D(f,f')=\cup_{1\le i\le t}\{c_i(V):V\in S_i\}$ and
Let $E(f,f')=D(f,f')\cup\text{singset}(f)$. 
Note that clearly $D(f,f')$ is a finite subset of $k$ and hence
by the Singset Theorem so is $E(f,f')$. Assuming (*),
write $f=f_1\dots f_s$ where $f_1,\dots,f_s$ are irreducible $Y$-monic 
members of $R\setminus k$, and let $\phi_i:R\to A_i=R/(f_iR)$
be the canonical epimorphism and identify $k$ with a subfield of
$L_i=\text{QF}(A_i)$. 
Without assuming (*), consider the conditions:
\begin{equation*}
\begin{cases}
\text{for a given $Q=(u,v)\in\mathcal A$ with $f(u,v)=0$, and}\\
\text{for every irreducible factor $g$ of $f'$ in $R\setminus k$ with
$g(u,v)=0$}\\
\text{we have $f_Y\in gR$ with $f\not\in gR$ and $f_X\not\in gR$,}\\
\end{cases}
\tag{$1'$}
\end{equation*}
\begin{equation*}
f'=f_Y\text{ and }Q=(u,v)\in\mathcal A\text{ with }f(u,v)=0
\tag{$2'$}
\end{equation*}
and
\begin{equation*}
\text{gcd}(f,f')=1.
\tag{$3'$}
\end{equation*}
Then we have (I) to (III) stated below.

\centerline{}

(I) If (*)+$(1')$ then
$$
I(f,f';Q)-I(f_X,f';Q)=I(X-u,f';Q)
$$
where all the terms are integers.

\centerline{}

(II) If (*)+$(2')$ then
$$
I(X-u,f';Q)
=\overline\chi(f;Q)+\sum_{1\le i\le s}\sum_{V\in\mathfrak R(f_i,Q)}
\text{ord}_V(d\phi_i(X))\\
$$
where all the terms are integers.

\centerline{}

(III) If (*)+$(3')$ then for all $\lambda\in k$ we have
$$
I(f,f';\mathcal A)=-Z(f,f')-P(f-\lambda,f')
$$
where all the terms are integers.

\centerline{}

PROOF OF (I). If (*)+$(1')$ then we have, with all terms integers,
\begin{align*}
\text{LHS of (I)}
&=\sum_{1\le i\le s'}[I(f,f'_i;Q)-I(f_X,f'_i;Q)]\\
&=\sum_{1\le i\le s'}\sum_{V\in\mathfrak R(f'_i,Q)}
[\text{ord}_V(\phi'_i(f))-\text{ord}_V(\phi'_i(f_X))]\\
&=\sum_{1\le i\le s'}\sum_{V\in\mathfrak R(f'_i,Q)}
\text{ord}_V(\phi'_i(X-u))\\
&=\sum_{1\le i\le s'}I(X-u,f'_i;Q)\\
&=\text{RHS of (I).}\\
\end{align*}

\centerline{}

PROOF OF (II). If (*)+$(2')$ then we have, with all terms integers,
\begin{align*}
\text{LHS of (II)}
&=-1+I(X-u,f;Q)\\
&=-1+\sum_{1\le i\le s}\sum_{V\in\mathfrak R(f_i,Q)}
\text{ord}_V(\phi_i(X-u))\\
&=-1+\sum_{1\le i\le s}\sum_{V\in\mathfrak R(f_i,Q)}
[\text{ord}_V(d\phi_i(X))+1]\\
&=\text{RHS of (II).}\\
\end{align*}

\centerline{}

PROOF OF (III). If (*)+$(3')$ then we have, with all terms integers,
\begin{align*}
\text{LHS of (III)}
&=-Z(f,f')-P(f,f')\\
&\quad(\text{since number of zeros of a function equals number of its poles})\\
&=\text{RHS of (III)}\\
&\quad(\text{since $\text{ord}_V\phi'_i(f)=\text{ord}_V\phi'_i(f-\lambda)$
for all $V\in P_i$ and $\lambda\in k$}.)\\
\end{align*}

\centerline{}

PROOF OF (9.4). Assume (*)+(4*). Write $f=f_1\dots f_r$ with pairwise distinct
irreducible $f_1,\dots,f_r$ in $R\setminus k$. 
Let $\phi_i:R\to A_i=R/(f_i)R$ be the canonical epimorphism and identify $k$ 
with a subfield of $L_i=\text{QF}(A_i)$. 
Letting $\overline\sum$ stand for summation over 
$\{Q=(u,v)\in\mathcal A:f(u,v)=0\}$ we get, with all terms integers,
\begin{align*}
\text{RHS of (9.4)}
&=1-N+\overline\sum[I(f,f_Y;Q)-I(f_X,f_Y;Q)]\\
&=1-N+\overline\sum
\left[\overline\chi(f;Q)+\sum_{1\le i\le s}\sum_{V\in\mathfrak R(f_i,Q)}
\text{ord}_V(d\phi_i(X))\right]\\
&\qquad\qquad\qquad\qquad\qquad
(\text{by taking $f'=f_Y$ in (I) and (II)})\\
&=1-N+\overline\chi(f;\mathcal A)+\sum_{1\le i\le s}
\sum_{V\in\mathfrak R(f_i,\mathcal A)}\text{ord}_V(d\phi_i(X))\\
&=1-N+\overline\chi(f;\mathcal A)+(2\gamma(f)-2)-\sum_{1\le i\le s}
\sum_{V\in\mathfrak R(f_i,\infty)}\text{ord}_V(d\phi_i(X))\\
&\qquad\qquad\qquad\qquad
(\text{by genus formula for $f_i$ and definition of $\gamma(f)$})\\
&=2\gamma(f)+\overline\chi(f;\mathcal A)-N-1-\sum_{1\le i\le s}
\sum_{V\in\mathfrak R(f_i,\infty)}[\text{ord}_V(\phi_i(X))-1]\\
&\quad\qquad\qquad\qquad\qquad\qquad\qquad\qquad\qquad
(\text{because ord}_V\phi_i(X)\ne 0)\\
&=2\gamma(f)+\overline\chi(f;\mathcal A)-N-1
-[-N-(\overline\chi(f;\infty)+1)]\\
&=2\gamma(f)+\overline\chi(f;\mathcal P)\\
&=\text{LHS of (9.4).}\\
\end{align*}

\centerline{}

PROOF OF (9.5) AND (9.6).
Assume (*). 
Let $\overline f=\text{rad}(f)$.  Then by the argument in the proof of the
Singset Theorem we see that
\begin{equation*}
[f]=f/\overline f
\quad\text{ and }\quad
\widehat{[f]}=[f]h\text{ where $h\in R$ with }\mathcal V(f,h)=\emptyset.
\tag{i}
\end{equation*}
and
$$
\widehat{[f]}=\prod_{c\in\text{multset}(f)^*}[f-c]
\quad\text{and hence }\quad |\text{multset}(f)^*|\le\text{deg}_Y\widehat{[f]}.
$$
which proves (9.6).
Applying (9.4) to $\overline f$ we get
$$
\rho_a(f)=1-\text{deg}_Y\overline f
+I(\overline f,\overline f_Y;\mathcal A)
-I(\overline f_X,\overline f_Y;\overline f).
$$ 
The first two terms of the above RHS combine to give
$1-N+\text{deg}_Y[f]$ which are the first three terms of the RHS of (9.5).
It remains to compare the intersection multiplicity terms, i.e., the proof of
(9.5) will be completed by showing that 
\begin{equation*}
I(\overline f,\overline f_Y;\mathcal A)
-I(\overline f_X,\overline f_Y;\overline f)
=I(f,f_Y/\widehat{[f]};\mathcal A)-I(f_X,f_Y/\widehat{[f]};f).
\tag{i*}
\end{equation*}
For any $Q=(u,v)\in\mathcal A$ with $f(u,v)=0$, by 
taking $(\overline f,\overline f_Y)$ for $(f,f')$ in (I) we get
\begin{equation*}
I(\overline f,\overline f_Y;Q)
-I(\overline f_X,\overline f_Y;Q)
=I(X-u,\overline f_Y;Q),
\tag{ii}
\end{equation*}
and, in view of (i), by taking $(f,f_Y/\widehat{[f]})$ for $(f,f')$ 
in (I) we get
\begin{equation*}
I(f,f_Y/\widehat{[f]};Q)-I(f_X,f_Y/\widehat{[f]};Q)
=I(X-u,f_Y/\widehat{[f]};Q).
\tag{iii}
\end{equation*}
By summing over $Q=(u,v)\in\mathcal A$ with $f(u,v)=0$, the LHS of (ii) gives
the LHS of (i*), and the RHS of (iii) gives the RHS of (i*). Therefore the
proof of (9.5) will be complete by proving that, for any $Q=(u,v)\in\mathcal A$ 
with $f(u,v)=0$, the RHS of (ii) equals the RHS of (iii). 
Let $\psi:R\to R/(X-u)R$ be the canonical epimorphism, and let $W$ be the 
localization of $\psi(R)$ at $\psi((Y-v)R)$. Now
\begin{align*}
\text{RHS of (ii)}
&=\text{ord}_W\psi(\overline f_Y)\\
&=-1+\text{ord}_W\psi(\overline f)
\quad(\text{because }\overline f(u,v)=0)\\
&=-1+\text{ord}_W\psi(f/\widehat{[f]})\quad(\text{by (i)})\\
&=\left[-1+\text{ord}_W\psi(f)\right]-\text{ord}_W\psi(\widehat{[f]})\\
&=\text{ord}_W\psi(f_Y)-\text{ord}_W\psi(\widehat{[f]})
\quad(\text{because }f(u,v)=0)\\
&=\text{ord}_W\psi(f_Y/\widehat{[f]})\\
&=\text{RHS of (iii)}
\end{align*}
and this completes the proof.

\centerline{}

PROOF OF (9.7) TO (9.12). 
For any $k$-automorphism $\sigma$ of $R$ we clearly have
$\rho_a(f-c)-\rho_a(\sigma(f)-c)
=|\mathcal V_{\infty}(\sigma(f)-c)|-|\mathcal V_{\infty}(f-c)|$ 
for all $c\in k$, and 
hence: (9.7) is true for $f\Leftrightarrow$ it is true for $\sigma(f)$.
After having defined $\rho_\pi(f)$ and defset$(f)$,
for any $k$-automorphism $\sigma$ of $R$ we clearly have 
$\jmath(\sigma(f))=\jmath(f)$, and
hence: (9.9) is true for $f\Leftrightarrow$ it is true for $\sigma(f)$.
Also it is well-known that $\sigma(f)$ satisfies (*) for some 
$k$-automorphism $\sigma$ of $R$. Therefore in proving (9.7) to (9.12), without 
loss of generality we may and {\bf we shall assume that $f$ satisfies (*).}
Let $f'=f_Y/\widehat{[f]}$. Clearly $f_X$ and $f_Y$ are unchanged if we 
replace $f$ by $f-c$ with $c\in k$, and hence so are $\widehat{[f]}$ and $f'$.
Given any $c\in k$, by taking $f-c$ for $f$ in (9.5) we get
\begin{equation*}
\rho_a(f-c)=1-N+\text{deg}_Y[f-c]+I(f-c,f';\mathcal A)-I(f_X,f';f-c)
\tag{1}
\end{equation*}
with all terms integers. 
By (9.6) we have gcd$(f-c,f')=1$ and hence by taking $(0,f-c)$ for 
$(\lambda,f)$ in (III) we get
\begin{equation*}
I(f-c,f';\mathcal A)=-Z(f-c,f')-P(f,f')
\tag{2}
\end{equation*}
where all the terms are integers.
In view of (9.6), by the definitions of the sets $D(f,f')$,
multset$(f)^*$, singset$(f)$ and $E(f,f')$ we see that
\begin{equation*}
\begin{cases}
\text{$Z(f-c,f')>0$ or $Z(f-c,f')=0$}\\
\text{according as $c\in D(f,f')$ or }c\not\in D(f,f'),
\end{cases}
\tag{3}
\end{equation*}
and
\begin{equation*}
\text{deg$_Y[f-c]>0$ or $=0$ according as $c\in\text{multset}(f)^*$ or
$c\not\in\text{multset}(f)^*$,}
\tag{4}
\end{equation*}
and
\begin{equation*}
I(f_X,f';f-c)\ge 0\text{ and if $c\not\in\text{singset}(f)$ then }
I(f_X,f';f-c)=0,
\tag{5}
\end{equation*}
and
\begin{equation*}
E(f,f')=D(f,f')\cup\text{singset}(f)\text{ with }
\text{multset}(f)^*\subset\text{singset}(f).
\tag{6}
\end{equation*}
By (1) to (6) we get
\begin{equation*}
c\not\in E(f,f')\Rightarrow
\rho_a(f-c)=1-N-P(f,f')
\tag{7}
\end{equation*}
with all terms integers. Since $E(f,f')$ is a finite set and the above RHS is 
independent of $c$, this proves (9.7) and establishes the existence of 
$\rho_\pi(f)$ and defset$(f)$. Now by (1) to (7) we see that
\begin{equation*}
\rho_\pi(f)=1-N-P(f,f')
\tag{8}
\end{equation*}
and
\begin{equation*}
\text{defset}(f)\subset E(f,f')
\tag{9}
\end{equation*}
and
\begin{equation*}
\rho_\pi(f)=(1-N)+\widehat I(f,f';\mathcal A)
\tag{10}
\end{equation*}
with all terms integers, which proves (9.8).
By (9.5) and (9.8) we get
\begin{equation*}
\begin{cases}
\rho_\pi(f)-\rho_a(f-c)=&-\text{deg}_Y[f-c] +I(f_X,f';f-c)\\
&+\left(\widehat I(f,f';\mathcal A)-I(f-c,f';\mathcal A)\right).
\end{cases}
\tag{11}
\end{equation*}
Since the second parenthesis term and the third big parenthesis term
in the above RHS are obviously nonnegative, in view of (4) we see that 
\begin{equation*}
\rho_\pi(f)-\rho_a(f-c)\ge -\text{deg}_Y[f-c]
\tag{12}
\end{equation*}
and
\begin{equation*}
\rho_\pi(f)\ge\rho_a(f-c)\text{ for all }c\in k\setminus\text{multset}(f)^*
\tag{13}
\end{equation*}
which proves (9.9). 
By (1), (2) and (8) we get
\begin{equation*}
\rho_\pi(f)-\rho_a(f-c)=-\text{deg}_Y[f-c]+I(f_X,f';f-c)+Z(f-c,f')
\tag{14}
\end{equation*}
with all terms integers.
Now, with all terms integers, we have
\begin{align*}
\text{RHS of (9.10)}
&=\sum_{c\in E(f,f')}[\rho_\pi(f)-\rho_a(f-c)]\\
&\qquad\qquad\qquad\qquad\qquad(\text{by (9)})\\
&=-\text{deg}_Y\widehat{[f]}
+\sum_{c\in E(f,f')}[I(f_X,f';f-c)+Z(f-c,f')]\\
&\quad\qquad\qquad\qquad\qquad\qquad(\text{by (9.6), (4), (6) and (14)})\\
&=-\text{deg}_Y\widehat{[f]}+I(f_X,f';\mathcal A)
+\sum_{c\in E(f,f')}Z(f-c,f')\\
&\qquad\qquad\qquad\qquad\qquad\qquad\qquad(\text{by (5) and (6)})\\
&=-\text{deg}_Y\widehat{[f]}+I(f_X,f';\mathcal A)
+\sum_{1\le i\le s'}\sum_{V\in Z_i}\text{ord}_V\phi'_i(f-c_i(V))\\
&\qquad\qquad\qquad\qquad\qquad\qquad\qquad\qquad\qquad(\text{by (7) and (9)})\\
&=-\text{deg}_Y\widehat{[f]}+I(f_X,f';\mathcal A)
+\sum_{1\le i\le s'}
\sum_{V\in Z_i}[\text{ord}_V\phi'_i(f_X)+\text{ord}_V\phi'_i(X)]\\
&\qquad
(\text{because $f_Y\in f'_iR$ and 
ord}_V\phi'_i(f-c_i(V))\ne 0\ne\text{ord}_V\phi'_i(X))\\
&=-\text{deg}_Y\widehat{[f]}+I(f_X,f';\mathcal A)+Z(f_X,f')+Z(X,f')\\
&=-\text{deg}_Y\widehat{[f]}+[-Z(f_X,f')-P(f_X,f')]+Z(f_X,f')+Z(X,f')\\
&\qquad\qquad\qquad\qquad
(\text{by taking $(0,f_X)$ for $(\lambda,f)$ in (III)})\\
&=-\text{deg}_Y\widehat{[f]}+Z(X,f')-P(f_X,f'))\\
\end{align*}
and we have
\begin{align*}
&-\text{deg}_Y\widehat{[f]}+Z(X,f')-P(f_X,f'))\\
&=-\text{deg}_Y\widehat{[f]}+Z(X,f')
-\sum_{1\le i\le s'}\sum_{V\in P_i}\text{ord}_V\phi'_i(f_X)\\
&=-\text{deg}_Y\widehat{[f]}+Z(X,f')
-\sum_{1\le i\le s'}\sum_{V\in P_i}
[-\text{ord}_V\phi'_i(X)+\text{ord}_V\phi'_i(f)]\\
&\qquad\qquad
(\text{because $f_Y\in f'_iR$ and 
ord}_V\phi'_i(f)\ne 0\ne\text{ord}_V\phi'_i(X))\\
&=-\text{deg}_Y\widehat{[f]}+Z(X,f')+P(X,f')-P(f,f')\\
&=-\text{deg}_Y\widehat{[f]}+Z(X,f')+P(X,f')-P(f,f')\\
&=-\text{deg}_Y\widehat{[f]}-\text{deg}_Yf'-P(f,f')\\
&=-\text{deg}_Yf_Y-P(f,f')\\
&=1-N-P(f,f')\\
&=\rho_\pi(f)\\
&\quad(\text{by (8)})\\
\end{align*}
and so we get
\begin{equation*}
\rho_\pi(f)=1-|\mathcal V_\infty(f)|
+\sum_{c\in\text{defset}(f)}[\rho_\pi(f)-\rho_a(f-c)]
\tag{15}
\end{equation*}
and hence
\begin{equation*}
\zeta(f)=-1
\tag{16}
\end{equation*}
which proves (9.10) and (9.11).
By (9.6) we see that
$$
c\not\in\text{multset}(f)^*\Rightarrow I(f_X,f';f-c)=I(f_X,f_Y;f-c)
$$
and hence
$$
c\in\text{singset}(f)\setminus\text{multset}(f)^*\Rightarrow
I(f_X,f';f-c)>0
$$
and therefore, because the third big parenthesis term in the RHS of (11) is 
nonnegative, in view of (4) and (11) we conclude that
\begin{equation*}
\text{singset}(f)\setminus\text{multset}(f)^*\subset\text{defset}(f).
\tag{17}
\end{equation*}
Rewriting the jungian excess formula (9.9) we get
\begin{equation*}
\rho_\pi(f)=\sum_{c\in\text{defset}(f)}[\rho_\pi(f)-\rho_a(f-c)]
\tag{18}
\end{equation*}
and adding $\rho_a(f)-\rho_\pi(f)$ to both sides we obtain
\begin{equation*}
\rho_a(f)=\sum_{c\in\text{defset}(f)\setminus\{0\}}[\rho_\pi(f)-\rho_a(f-c)].
\tag{19}
\end{equation*}
By (4), (12) and (19) we get
$$
\rho_a(f)+\sum_{0\ne c\in\text{multset}(f)^*}\text{deg}_Y[f-c]
\ge\sum_{0\ne c\in\text{defset}(f)\setminus\text{multset}(f)^*}
[\rho_\pi(f)-\rho_a(f-c)]
$$
and hence by (9.6) we have
$$
\rho_a(f)+\text{deg}_Y\widehat{[f]}
\ge\sum_{0\ne c\in\text{defset}(f)\setminus\text{multset}(f)^*}
[\rho_\pi(f)-\rho_a(f-c)]
$$
and therefore, since every term in the above RHS is nonnegative, we obtain
$$
|\text{defset}(f)\setminus\{0\}|\le \rho_a(f)+\text{deg}_Y\widehat{[f]}
$$
and hence
\begin{equation*}
|\text{defset}(f)|\le 1+\rho_a(f)+\text{deg}_Y\widehat{[f]}
\tag{20}
\end{equation*}
By (17) we have
$$
|\text{singset}(f)|\le|\text{defset}(f)|+|\text{multset}(f)^*|
$$
and hence, in view of (9.6), by (20) we get 
\begin{equation*}
|\text{singset}(f)|\le 1+\rho_a(f)+2\text{deg}_Y\widehat{[f]}.
\tag{21}
\end{equation*}
The proof of (9.12) is completed by (17), (20) and (21). 

\centerline{}

QUESTION 9. Is it always true that
$\text{redset}(f)\subset\text{defset}(f)$?

\centerline{}

QUESTION 10.  More explicitly, is it possible to find a nonsingular reducible 
$f\in R\setminus k$ satisfying (*) such that $[f]=1$ and
$\rho_a(f)=\rho_\pi(f)$.  If such an $f$ is found then
$0\in\text{redset}(f)\setminus\text{defset}(f)$ 
and we have a negative answer to Question 9.
In this connection it would be worthwhile to study the famous examples of
Artal-Bartolo, Cassou-Nogues, and Luengo-Velasco \cite{ACL} in which they
construct curves whose redsets and singsets are empty and whose ranks are
arbitrarily high.

\centerline{}

QUESTION 11. In connection with the Defset Theorem (9.12), consider the 
Example given by the pencil
$(h(Y)-c)_{c\in k}$ where $f=h(Y)\in k[Y]$ is monic of degree $N>1$;
using (9.5) we can easily show that: (i) $\rho_\pi(f)=1-N$,
(ii) $\rho_a(f-c)=1-N+\text{deg}_Y[f-c]$ for all $c\in k$, and
(iii) $\rho_\pi(f)=\sum_{c\in k}[\rho_\pi(f)-\rho_a(f-c)]$;
thus $\rho_a(f-c)\ge\rho_\pi(f)$ for all $c\in k$, and
$\rho_a(f-c)>\rho_\pi(f)$ for precisely those $c$ for which $c=h(c')$ with
$h_Y(c')=0$.
It is easy to make similar examples for pencils obtained by replacing $Y$
by any nonconstant polynomial in $R$.
However we ask: is it true that if the pencil 
$(f-c)_{c\in k}$ is noncomposite, i.e., if $k(f)$ is m.a. in $k(X,Y)$, then 
$\rho_\pi(f)\ge\rho_a(f-c)$ for all $c\in k$?
If this has a positive answer then, for noncomposite 
pencils, the term $2\text{deg}_Y\widehat{[f]}$ may be dropped from the
estimate given in (9.12). 

\centerline{}

QUESTION 12.  Can rank and defset be generalized to $n>2$?

\centerline{}

{\bf MORAL.} Thus it is all a matter of measuring change in a 
quantity relative to the corresponding change in another quantity on which
the previous quantity depends. This after all is the Newton-Leibnitz idea
of derivative. By avoiding limits, which are algebraically awkward to take,
it also gives rise to the twisted derivative introduced in \cite{A07} 
and exploited in numerous succeeding papers summarized in \cite{A10} 
for calculating Galois groups
and hence fundamental groups. The same principle of measuring change gives
rise to the invariants $I,\beta,\gamma,\rho,\rho_a,\rho_\pi,\jmath,\zeta$.
Amongst these, $\beta$ measures the change in $I$, the quantities
$\beta,\gamma,\rho,\rho_a,\rho_\pi$ are different incarnations of the same
underlying reality, and the almost identical quantities $\jmath$ and $\zeta$ 
measure the change in $\rho$.
Or, as may be easier to remember, the genus $\gamma$ measures the change in the
intersection multiplicity $I$, while the Zeuthen-Segre invariant $\zeta$ 
measures the change in the genus $\gamma$.

\centerline{}

\centerline{\bf Section 10: Defset of a General Pencil}

In this section we continue to assume $n=2$ with $(X_1,X_2)=(X,Y)$, and $k$ is 
of characteristic zero with $k=$ its algebraic closure $k^*$. 
We shall now extend our study of rank and defset to a general pencil 
$(f-cw)_{c\in k\cup\{\infty\}}$ where 
$w\in R^{\times}=R\setminus\{0\}$ with gcd$(f,w)=1$ and
$d=\text{max}(\text{deg}(f),\text{deg}(w))$.
Recall that by convention $f-\infty w=w$. 

Let $\mathcal L_{\infty}=\mathcal P\setminus\mathcal A
=\{(\infty,v):v\in k\cup\{\infty\}\}$ and call this the line at infinity.
In homogeneous coordinates $(X,Y,Z)$ we think of $\mathcal L_{\infty}$
as given by $Z=0$, and also as the $(Y/Z)$-axis. As abbreviations, 
for the two special
points of this line, i.e., for the two infinite points of the fundamental
triangle, we put $(\infty)=(\infty,0)$ and $((\infty))=(\infty,\infty)$
and note that in homogeneous coordinates they are $(1,0,0)$ and $(0,1,0)$.

Let $R_0=k[[X,Y]]$ and note that $k(X,Y)\cap R_0$ is the localization of
$R$ at the maximal ideal $(X,Y)R$ and $R_0$ is its completion.
Let $e\ge 0$ be an integer.  Define the $e$-th homogenization of any $g\in R$ 
by putting $g_{[e]}=g_{[e]}(X,Y,Z)=Z^{e+\text{deg}(g)}g(X/Z,Y/Z)$ in case 
$g\ne 0$, and $g_{[e]}=g_{[e]}(X,Y,Z)=0$ in case $g=0$.
For any $Q=(u,v)\in\mathcal P$, define the $k$-monomorphism 
$T^{[Q,e]}:R\to R_0$ with $g\mapsto g^{[Q,e]}=g^{[Q,e]}(X,Y)$ by putting
\begin{equation*}
g^{[Q,e]}(X,Y)=\begin{cases}
g(X+u,Y+v)&\text{if }Q\in\mathcal A\\
g_{[e]}(1,Y+v,X)&\text{if }((\infty))\ne Q\in\mathcal L_{\infty}\\
g_{[e]}(X,1,Y)&\text{if }Q=((\infty))\\
\end{cases}
\tag{10.1}
\end{equation*}
and call this the Taylor map at $[Q,e]$. 

Maclaurin and Taylor were two disciples of Newton of calculus-fame, the 
former expanding things around the origin and the latter around other points! 
Newton having accomplished many other things, the adjective 
``of calculus-fame'' is clearly directed more towards the disciples!!.

For any $Q\in\mathcal P$ let $T^{[Q]}:R\to R_0$ be given by 
$g\mapsto g^{[Q]}=g^{[Q,0]}$. 
Note that if $g^{[Q]}(0,0)\ne 0$ then $\overline\chi(g;Q)=-1$, and
if $g^{[Q]}(0,0)=0\ne g$ then $g^{[Q]}$ is product of $1+\overline\chi(g;Q)$
irreducibles in $R_0$. Also note that for any $g,h$ in $R$ we have
$I(g,h;Q)=[R_0/(g^{[Q]},h^{[Q]})R_0:k]$. 
For any $g,h,\overline g,\overline h$ in $R$ we put
\begin{equation*}
I(g,h;\overline g\setminus\overline h)
=\sum_{\{Q\in\mathcal A:\overline g^{[Q]}(0,0)=0\ne\overline h^{[Q]}(0,0)\}}
I(g,h;Q).
\tag{10.2}
\end{equation*}
Let
$$
t=|\mathcal V(f,w)|
$$
and let $Q_1=(u_1,v_1),\dots,Q_t=(u_t,v_t)$ in $\mathcal A$ be the $t$ points 
of $\mathcal V(f,w)$ called the finite base points of the pencil. 
For any $g,h,\overline g$ in $R$ we put
\begin{equation*}
I(g,h;\overline g\setminus\mathcal V(f,w))
=\sum_{\{Q\in\mathcal A\setminus\{Q_1,\dots,Q_t\}:
\overline g^{[Q]}(0,0)=0\}}
I(g,h;Q).
\tag{10.3}
\end{equation*}

To avoid repeating the considerations of the previous section and for 
simplicity of calculation, most of the time we shall suppose that:
\begin{equation*}
\text{the elements }f,w,1\text{ are linearly independent over }k,
\tag{r*}
\end{equation*}
i.e., by replacing $f,w$ by suitable $k$-linear combinations of them, our pencil
cannot be converted into the special pencil considered in the previous section.

 From the previous section recall condition
\begin{equation*}
\text{$f$ is $Y$-monic of $Y$-degree $N>0$}
\tag{*}
\end{equation*}
and consider condition
\begin{equation*}
\text{(*) and $w$ is $Y$-monic of $Y$-degree }0<M<N.
\tag{**}
\end{equation*}
We claim that:
\begin{equation*}
\text{(r*)}\Rightarrow\begin{cases}
\text{there exists a $k$-automorphism $\sigma$ of $R$}\\
\text{together with $c_i^*\ne c_2^*$ in $k\cup\{\infty\}$ such that}\\
\text{if we let the
pair $(\sigma(f-c_1^*w),\sigma(f-c_2^*w))$ be called $(f,w)$}\\ 
\text{then condition (**) is satisfied.} 
\end{cases}
\tag{10.4}
\end{equation*}
Namely, by a well-known argument there exists a $k$-automorphism $\sigma_1$
of $R$ such that $\sigma_1(f), \sigma_1(w)$ are both Y-monic. If their degrees 
are distinct, then we already have (**). If not, then we  can find $c_1^*,
c_2^* \in k$ such that $\sigma_1(f-c_1^*w)$ has bigger degree than
$\sigma_1(f - c_2^*w)$, which must be positive because of (r*). By applying 
a suitable $k$-automorphism $\sigma_2$ of $R$ to
these two polynomials, we can arrange that condition (**) is satisfied. Now 
take $\sigma$ to be the composition of $\sigma_2$ with $\sigma_1$.

\centerline{}

Analogous to (9.5) we shall prove that:
\begin{equation*}
({}^{**})\Rightarrow\begin{cases}
\text{for any $c\in k$ we have:}\\
\rho_a(f-cw)\\
=(1-N)\\
\quad+\text{deg}_Y[f-cw]\\
\quad+I\left(f-cw,\frac{f_Yw-fw_Y}{\widehat{[f,w]}};(f-cw)\setminus w\right)\\
\quad-I\left(f_Xw-fw_X,\frac{f_Yw-fw_Y}{\widehat{[f,w]}};
(f-cw)\setminus w\right)\\
\quad+\sum_{1\le i\le t}\left(I(X-u_i,\frac{f-cw}{[f-cw]};Q_i)-1\right)\\
\end{cases}
\tag{10.5}
\end{equation*}
with all terms integers, where
$$
\widehat{[f,w]}=\text{gcd}(f_Xw-fw_X,f_Yw-fw_Y)
$$
with, as before, the gcd made unique by requiring it to be $Y$-monic. While
proving (10.5) we shall show that:
\begin{equation*}
({}^{**})\Rightarrow
\text{for any $c\in k\cup\{\infty\}$ we have }
[f-cw]=\text{gcd}(f-cw,\widehat{[f,w]}).
\tag{10.6}
\end{equation*}
As a consequence of (10.5) we shall show that:
\begin{equation*}
\text{(r*)}\Rightarrow\begin{cases}
\text{there exists a unique integer $\rho_{\pi}(f,w)$ together with}\\
\text{a unique finite subset defset$(f,w)$ of $k\cup\{\infty\}$ such that}\\
\text{$\rho_a(f-cw)=\rho_\pi(f,w)$ 
for all $c\in (k\cup\{\infty\})\setminus\text{defset}(f,w)$ and}\\
\text{$\rho_a(f-cw)\ne\rho_\pi(f,w)$ for all $c\in\text{defset}(f,w)$.}\\
\end{cases}
\tag{10.7}
\end{equation*}
In reference to (10.7), we put
$$
\rho_\pi(f,w)=\text{the pencil-rank of the pencil }
(f-cw)_{c\in k\cup\{\infty\}}
$$
and
$$
\text{defset}(f,w)=\text{the deficiency set of }(f,w).
$$
 From (10.7) we shall deduce that:
\begin{equation*}
(\text{**})\Rightarrow\begin{cases}
\rho_\pi(f,w)\\
=(1-N)\\
\quad+\text{max}_{c\in k}
I\left(f-cw,\frac{f_Yw-fw_Y}{\widehat{[f,w]}};(f-cw)\setminus w\right)\\
\quad+\sum_{1\le i\le t}\text{min}_{c\in k\setminus\text{multset}(f,w)^*}
\left(I(X-u_i,\frac{f-cw}{[f-cw]};Q_i)
-1\right)\\
\end{cases}
\tag{10.8}
\end{equation*}
with all terms integers.
 From (10.8) we shall deduce that:
\begin{equation*}
(\text{**})\Rightarrow\begin{cases}
\rho_\pi(f,w)\ge\rho_a(f-cw)\\
\text{for all }
c\in k\setminus(\text{multset}(f,w)^*\cup\text{conset}(f,w)^*)
\end{cases}
\tag{10.9}
\end{equation*}
where, without assuming $k=k^*$, the contact set of $(f,w)$ is defined by 
putting
$$
\text{conset}(f,w)^*=\cup_{Q=(u,v)\in\mathcal V(f,w)^*}\text{conset}(f,w;Q)^*
$$
with the set conset$(f,w;Q)^*$ of size at most one defined by
$$
\begin{cases}
\text{conset}(f,w;Q)^*\\
=\{c\in k^*:I(X-u,f-cw;Q)\\
\quad\qquad\qquad>I(X-u,f-c'w;Q)\text{ for some }
c'\in k^*\setminus\text{multset}(f,w)^*\}.
\end{cases}
$$
 From (10.8) we shall deduce that:
\begin{equation*}
\text{(**)}\Rightarrow\begin{cases}
\rho_a(f)+\rho_a(w)\\
=1-t
+\sum_{c\in\text{defset}(f,w)\setminus\{0,\infty\}}[\rho_\pi(f,w)-\rho_a(f-cw)].
\end{cases}
\tag{10.10}
\end{equation*}
Clearly (10.10) is equivalent to saying that:
\begin{equation*}
\text{(**)}\Rightarrow\zeta(f,w)=-1
\tag{10.11}
\end{equation*}
where
$$
\zeta(f,w)=-t-2\rho_\pi(f,w)
+\sum_{c\in\text{defset}(f,w)}[\rho_\pi(f,w)-\rho_a(f-cw)].
$$

 Finally, from (10.11) we shall deduce the:

\centerline{}

(10.12) GENERAL DEFSET THEOREM. If (**) then we have:

(i) $(a(f,w)\setminus b(f,w))\subset\text{defset}(f,w)$ where
$a(f,w)=k\cap\text{singset}(f,w)$ and
$b(f,w)=\text{multset}(f,w)^*\cup\text{conset}(f,w)^*$.

(ii)
$|\text{defset}(f,w)|\le \rho_a(f)+\rho_a(w)+\text{deg}_Y\widehat{[f,w]}+2t+1.$

(iii) $|\text{singset}(f,w)|
\le\rho_a(f)+\rho_a(w)+2\text{deg}_Y\widehat{[f,w]}+3t+2$.

\centerline{}

Before turning to the proof of items (10.5) to (10.12), let us establish 
some common

\centerline{}

NOTATION AND CALCULATION.
Given any $f'\in R^{\times}$, write
$f'=\nz f'_1\dots f'_{s'}$ where $f'_1,\dots,f'_{s'}$ are
irreducible members of $R\setminus k$ and $\nz\in k^{\times}$, 
and let $\phi'_i:R\to A'_i=R/(f'_iR)$
be the canonical epimorphism and identify $k$ with a subfield of
$L'_i=\text{QF}(A'_i)$. Let
$Z_i=\{V\in\mathfrak R(f'_i,\infty):
\text{ord}_V\phi'_i(f)\ge \text{ord}_V\phi'_i(w)\}$ and
$P_i=\{V\in\mathfrak R(f'_i,\infty):
\text{ord}_V\phi'_i(f)<\text{ord}_V\phi'_i(w)\}$,
where the letters $Z$ and $P$ are meant to suggest zeros and poles,
and note that these are obviously finite sets. 

Consider  the conditions:
\begin{equation*}
\text{gcd}(w,f')=1,
\tag{$0'$}
\end{equation*}
\begin{equation*}
\begin{cases}
\text{for a given $Q=(u,v)\in\mathcal A$ with $f(u,v)=0$, and}\\
\text{for every irreducible factor $g$ of $f'$ in $R\setminus k$ with
$g(u,v)=0\ne w(u,v)$}\\
\text{we have $w^2(f/w)_Y\in gR$ with $f\not\in gR$ and 
$w^2(f/w)_X\not\in gR$,}\\
\end{cases}
\tag{$1'$}
\end{equation*}
\begin{equation*}
f'=w^2(f/w)_Y\text{ and }Q=(u,v)\in\mathcal A\text{ with }f(u,v)=0\ne w(u,v)
\tag{$2'$}
\end{equation*}
and
\begin{equation*}
\text{gcd}(f,f')=1.
\tag{$3'$}
\end{equation*}

Assuming $(0')$:
Let
$$
Z((f,w),f')=\sum_{1\le i\le s'}\sum_{V\in Z_i}
[\text{ord}_V\phi'_i(f)-\text{ord}_V\phi'_i(w)]
$$
and
$$
P((f,w),f')=\sum_{1\le i\le s'}\sum_{V\in P_i}
[\text{ord}_V\phi'_i(f)-\text{ord}_V\phi'_i(w)]
$$
and note that $Z((f,w),f')$ is a nonnegative integer or $\infty$ according as
gcd$(f,f')=1$ or gcd$(f,f')\ne 1$, 
and $P((f,w),f')$ is always a nonpositive integer.  
Clearly for each $V\in Z_i$ there is a unique $c_i(V)\in k$ 
such that ord$_V\phi'_i(f-c_i(V)w)>\text{ord}_V\phi'_i(w)$.
Let $D((f,w),f')=\cup_{1\le i\le t}\{c_i(V):V\in Z_i\}$ and
let $E((f,w),f')=D((f,w),f')\cup\text{singset}((f,w))$. 
Note that clearly $D((f,w),f')$ is a finite subset of $k$ and hence
by the General Singset Theorem so is $E((f,w),f')$. 

Without assuming $(0')$ but assuming (**):
Write $f=f_1\dots f_s$ where $f_1,\dots,f_s$ are irreducible $Y$-monic 
members of $R\setminus k$, and let $\phi_i:R\to A_i=R/(f_iR)$
be the canonical epimorphism and identify $k$ with a subfield of
$L_i=\text{QF}(A_i)$. 

With these conditions in mind, we have (I) to (IV) stated below.

\centerline{}

(I) If (**)+$(1')$ then
$$
I(f,f';Q)-I(w^2(f/w)_X,f';Q)=I(X-u,f';Q)
$$
where all the terms are integers.

\centerline{}

(II) If (**)+$(2')$ then
$$
I(X-u,f';Q)
=\overline\chi(f;Q)+\sum_{1\le i\le s}\sum_{V\in\mathfrak R(f_i,Q)}
\text{ord}_Vd\phi_i(X)\\
$$
where all the terms are integers.

\centerline{}

(III) If (**)+$(0')+(3')$ then for all $\lambda\in k$ we have
$$
I(f,f';\mathcal A)-I(w,f';\mathcal A)=-Z((f,w),f')
-P((f-\lambda w,w),f')
$$
where all the terms are integers.

\centerline{}

(IV) If (*) then
$$
1-\text{deg}_Yf=2+\overline\chi(f;\infty)
+\sum_{1\le i\le s}\sum_{V\in\mathfrak R(f_i,\infty)}\text{ord}_Vd\phi_i(X)
$$
with all terms integers.

\centerline{}

PROOF OF (I). If (**)+$(1')$ then we have, with all terms integers,
\begin{align*}
\text{LHS of (I)}
&=\sum_{1\le i\le s'}\left[(I(f,f'_i;Q)-I(w^2(f/w)_X,f'_i;Q)\right]\\
&=\sum_{1\le i\le s'}\sum_{V\in\mathfrak R(f'_i,Q)}
\left[\text{ord}_V\phi'_i(f)-\text{ord}_V\phi'_i(w^2(f/w)_X)\right]\\
&=\sum_{1\le i\le s'}\sum_{V\in\mathfrak R(f'_i,Q)}
\text{ord}_V\phi'_i(X-u)\\
&=\sum_{1\le i\le s'}I(X-u,f'_i;Q)\\
&=\text{RHS of (I).}\\
\end{align*}

\centerline{}

PROOF OF (II). If (**)+$(2')$ then we have, with all terms integers,
\begin{align*}
\text{LHS of (II)}
&=-1+I(X-u,f;Q)\\
&=-1+\sum_{1\le i\le s}\sum_{V\in\mathfrak R(f_i,Q)}
\text{ord}_V\phi_i(X-u)\\
&=-1+\sum_{1\le i\le s}\sum_{V\in\mathfrak R(f_i,Q)}
[\text{ord}_Vd\phi_i(X)+1]\\
&=\text{RHS of (II).}\\
\end{align*}

\centerline{}

PROOF OF (III). If (**)+$(0')+(3')$ then we have, with all terms integers,
\begin{align*}
\text{LHS of (III)}
&=-Z(f,f')-P(f,f')+Z(w,f')+P(w,f')\\
&\quad(\text{since number of zeros of a function equals number of its poles})\\
&=\text{RHS of (III)}\\
&\quad(\text{since $\text{ord}_V\phi'_i(f)=\text{ord}_V\phi'_i(f-\lambda w)$
for all $V\in P_i$ and $\lambda\in k$}.)\\
\end{align*}

\centerline{}

PROOF OF (IV). If (*) then we have, with all terms integers,
\begin{align*}
\text{RHS of (IV)}
&=2+\overline\chi(f;\infty)+\sum_{1\le i\le s}
\sum_{V\in\mathfrak R(f_i,\infty)}\left(\text{ord}_V\phi_i(X)-1\right)\\
&\qquad(\text{because ord}_V\phi_i(X)<0)\\
&=2+\overline\chi(f;\infty)-\text{deg}_Yf
-\sum_{1\le i\le s}|\mathfrak R(f_i,\infty)|\\
&=1-\text{deg}_Yf\\
&=\text{LHS of (IV).}
\end{align*}

\centerline{}

PROOF OF (10.5) AND (10.6).
Assume (**). 
Then by the argument in the proof of the
Singset Theorem we see that
\begin{equation*}
\text{for any }c\in k\cup\{\infty\}\text{ we have }
f-cw=[f-cw]\overline f\text{ where }\overline f=\text{rad}(f-cw)
\tag{i}
\end{equation*}
and
$$
\text{for any $c\in k\cup\{\infty\}$ we have }
[f-cw]=\text{gcd}(f-cw,\widehat{[f,w]})
$$
which proves (10.6).
Take $f'=\frac{f_Yw-fw_Y}{\widehat{[f,w]}}$.
For any $c\in k$ and $Q=(u,v)\in\mathcal A$ with $f(u,v)-cw(u,v)=0\ne w(u,v)$,
by using (I) we get
\begin{equation*}
I(f-cw,f';Q)-I(f_Xw-fw_X,f';Q)=I(X-u,f';Q)
\tag{ii}
\end{equation*}
with all terms integers. Let $\overline f=\text{rad}(f-cw)$.  
Let $\psi:R\to R/(X-u)R$ be the canonical epimorphism, and let $W$ be the 
localization of $\psi(R)$ at $\psi((Y-v)R)$.
Then we have, with all terms integers,
\begin{align*}
I(X-u,w^2(\overline f/w)_Y;Q)
&=\text{ord}_W\psi(w)+\text{ord}_W\psi(\overline f)-1\\
&=\text{ord}_W\psi(\overline f)-1\\
&\qquad(\text{since }w(u,v)\ne 0)\\
&=\text{ord}_W\psi(f-cw)-\text{ord}_W\psi([f-cw])-1\\
&=\text{ord}_W\psi(w^2((f-cw)/w)_Y)-\text{ord}_W\psi([f-cw])\\
&=\text{ord}_W\psi(w^2((f-cw)/w)_Y)-\text{ord}_W\psi(\widehat{[f-w]})\\
&\qquad(\text{since }\mathcal V(f-cw)\cap\mathcal V(\widehat{[f,w]}/[f-cw])
=\emptyset\\
&=\text{ord}_W\psi(f')\\
&=I(X-u,f';Q)
\end{align*}
and hence, with all terms integers, we have
\begin{equation*}
I(X-u,w^2(\overline f/w)_Y;Q)=I(X-u,f';Q).
\tag{iii}
\end{equation*}
Now, upon letting $\widetilde\sum$ and $\widehat\sum$ stand for summations over
$\{Q=(u,v)\in\mathcal A:f(u,v)-cw(u,v)=0\ne w(u,v)\}$ and
$\{Q=(u,v)\in\mathcal A:f(u,v)-cw(u,v)=0\}$ respectively, we have, with all
terms integers,
\begin{align*}
\text{RHS of (10.5)}
&=1-N+\text{deg}_Y[f-cw]
+\sum_{1\le i\le t}\left(I(X-u_i,\overline f;Q_i)-1\right)\\
&\qquad+\widetilde\sum I(X-u,w^2(\overline f/w)_Y;Q)\\
&\qquad\qquad\qquad\qquad(\text{by (i), (ii) and (iii)})\\
&=1-\text{deg}_Y\overline f
+\widehat\sum
\left[\overline\chi(\overline f;Q)+\sum_{1\le i\le\overline s}
\sum_{V\in\mathfrak R(\overline f_i,Q)}
\text{ord}_Vd\overline\phi_i(X)\right]\\
&\qquad(\text{by taking $\overline f$ for $f$ in (II) and writing 
$\overline f=\overline f_1\dots\overline f_{\overline s}$ with}\\ 
&\quad\qquad\text{irreducible $Y$-monic members 
$\overline f_1,\dots,\overline f_{\overline s}$ of $R\setminus k$ and}\\
&\quad\qquad\text{letting $\overline\phi_i:R\to R/(\overline f_iR)$
be the canonical epimorphism})\\
&=\rho_a(\overline f)\\
&\quad(\text{by (IV)})\\
&=\text{LHS of (10.5).}
\end{align*}

\centerline{}

PROOF OF (10.7) TO (10.12). 
In view of (10.4), as in the proof of (9.7) to (9.12), while proving
(10.7) to (10.12) we may and {\bf we shall assume that $(f,w)$ satisfies (**).}
Take 
$$
f'=\frac{f_Yw-fw_Y}{\widehat{[f,w]}}
\quad\text{ and }\quad
f''=f_X w - f w_X.
$$
Let
$$
a(f,w)=k\cap\text{singset}(f,w)\quad\text{ and }\quad
b(f,w)=\text{multset}(f,w)^*\cup\text{conset}(f,w)^*.
$$
Now in the RHS of (10.5), for all $c\in k\setminus\text{multset}(f,w)^*$
the second line is zero, for $c\in k\setminus D((f,w),f')$ the third line
attains a maximum, for $c\in\text{singset}(f,w)$ the fourth line is zero, and 
for $c\in k\setminus b(f,w)$ the fifth line attains a minimum.
This proves that
\begin{equation*}
\text{defset}(f,w)\subset\left(D((f,w),f')\cup\text{singset}(f,w)\cup
\text{conset}(f,w)^*\right)
\tag{1}
\end{equation*}
and establishes (10.7) and (10.8). Clearly (10.9) is evident from (10.8).
Also note that if $c\in\text{singset}(f,w)\setminus b(f,w)$ then the above 
consideration of the RHS lines of (10.5) establishes that 
$c\in\text{defset}(f,w)$, and hence
\begin{equation*}
(a(f,w)\setminus b(f,w))\subset\text{defset}(f,w)
\tag{2}
\end{equation*}
which proves part (i) of (10.12). Clearly $|\text{conset}(f,w)^*|\le t$
and hence by (10.9) and assuming (10.10) we get
\begin{equation*}
|\text{defset}(f,w)|\le\rho_a(f)+\rho_a(w)+t-1
+|\text{multset}(f,w)^*|+|\text{conset}(f.w)^*|+2
\tag{3}
\end{equation*}
where the last number $2$ is added for possible $0,\infty$ in defset$(f,w)$
not accounted by (10.10). This gives
\begin{equation*}
|\text{defset}(f,w)|\le\rho_a(f)+\rho_a(w)+\text{deg}_Y\widehat{[f,w]}
+2t+1
\tag{4}
\end{equation*}
showing that (10.10) $\Rightarrow$ part (ii) of (10.12). By (2) we see that
we get
$$
|\text{singset}(f,w)|\le|\text{defset}(f,w)|+\text{deg}_Y\widehat{[f,w]}
+t+1
$$
and hence by (4) we get
\begin{equation*}
|\text{singset}(f,w)|\le\rho_a(f)+\rho_a(w)+2\text{deg}_Y\widehat{[f,w]}
+3t+2
\tag{5}
\end{equation*}
showing that (10.10) $\Rightarrow$ part (iii) of (10.12). 
Thus, it only remains to prove (10.10).
We shall do this in STEPS (6) to (12).

STEP (6). Using the proof of (10.5) we shall now show that:
\begin{equation*}
\text{(***)}\Rightarrow\begin{cases}
\rho_a(w) &= (1-M)+\text{deg}_Y[w]\\
&\quad+I\left(w,f';w\setminus f\right)-I\left(f'',f';w\setminus f\right)\\
&\quad+\sum_{1\le i\le t}\left(I\left(X-u_i,\frac{w}{[w]};Q_i\right)-1\right)\\
\end{cases}
\tag{10.5*}
\end{equation*}
where condition
\begin{equation*}
\text{$f$ and $w$ are $Y$-monic of positive $Y$-degree }N\ne M.
\tag{***}
\end{equation*}
is obviously weaker than condition (**). So  let $\overline w=\text{rad}(w)$.  
Then by the argument in the proof of the Singset Theorem we see that:
$[w]=\text{gcd}(w,\widehat{[f,w]})$ and $w=[w]\overline{w}$.

The idea of the proof is to redo the calculations in (10.5) reversing
the role of $f,w$ and for this purpose, let us note that under our
current notation,
we have:
\begin{equation*}
w^2(f/w)_Y = -f^2(w/f)_Y \text{ and } w^2(f/w)_X = -f^2(w/f)_X
\tag{i}
\end{equation*}

Note that our arguments in (I) thru (IV) remain valid under exchange of
$f,w$, if we change the degree condition (**) by the weaker condition (***).
Indeed, the degrees, $N,M$ never enter the
calculations in (I) to (IV).

For  $Q=(u,v)\in \mathcal A$ with $w(u,v)=0\ne f(u,v)$,
by using calculations of (I) we get
\begin{equation*}
I(w,f';Q)-I(f'',f';Q)=I(X-u,f';Q)
\tag{ii}
\end{equation*}
with all terms integers.
Let $\psi:R\to R/(X-u)R$ be the canonical epimorphism, and let $W$ be the
localization of $\psi(R)$ at $\psi((Y-v)R)$.
Then we have, in view of (i) and with all terms integers,
\begin{align*}
I(X-u,f^2(\overline w/f)_Y;Q)
&=\text{ord}_W\psi(f)+\text{ord}_W\psi(\overline w)-1\\
&=\text{ord}_W\psi(\overline w)-1\\
&\qquad(\text{since }f(u,v)\ne 0)\\
&=\text{ord}_W\psi(w)-\text{ord}_W\psi([w])-1\\
&=\text{ord}_W\psi(f^2(w/f)_Y)-\text{ord}_W\psi([w])\\
&=\text{ord}_W\psi(f^2(w/f)_Y)-\text{ord}_W\psi(\widehat{[f,w]})\\
&\qquad(\text{since }\mathcal V(w)\cap\mathcal V(\widehat{[f,w]}/[w])
=\emptyset\\
&=\text{ord}_W\psi(f')\\
&=I(X-u,f';Q)
\end{align*}
and hence, with all terms integers, we have
\begin{equation*}
I(X-u,f^2(\overline w/f)_Y;Q)=I(X-u,f';Q).
\tag{iii}
\end{equation*}
Now, upon letting $\widetilde\sum$ and $\widehat\sum$ stand for summations over
$$\{Q=(u,v)\in\mathcal A:w(u,v)=0\ne f(u,v)\}
\quad\text{ and }\quad
\{Q=(u,v)\in\mathcal A:w(u,v)=0\}
$$ 
respectively, we have, with all terms integers,
\begin{align*}
\text{RHS of (10.5*)}
&=1-M+\text{deg}_Y[w]+\sum_{1\le i\le t}\left(I(X-u_i,\overline
w;Q_i)-1\right)\\
&\qquad+\widetilde\sum I(X-u,f^2(\overline w/f)_Y;Q)\\
&\qquad\qquad\qquad\qquad(\text{by (ii) and (iii)})\\
&=1-\text{deg}_Y\overline w
+\widehat\sum
\left[\overline\chi(\overline w;Q)+\sum_{1\le i\le\overline s}
\sum_{V\in\mathfrak R(\overline w_i,Q)}
\text{ord}_Vd\overline\phi_i(X)\right]\\
&\qquad(\text{by taking $\overline w$ for $w$ in (II) and writing
$\overline w=\overline w_1\dots\overline w_{\overline s}$ with}\\
&\quad\qquad\text{irreducible $Y$-monic members
$\overline w_1,\dots,\overline w_{\overline s}$ of $R\setminus k$ and}\\
&\quad\qquad\text{letting $\overline\phi_i:R\to R/(\overline w_iR)$
be the canonical epimorphism})\\
&=\rho_a(\overline w)\\
&\quad(\text{by (IV)})\\
&=\text{RHS of (10.5)*.}
\end{align*}

STEP (7).
Combining (10.5*) with (10.5) we see that for any $c\in k$ we have:
\begin{equation*}
\rho_a(f-cw)+\rho_a(w)+t-1=\sum_{1\le j\le 4}F_j(c)
\tag{7*}
\end{equation*}
where
$$
F_1(c)=(1-N-M)+\text{deg}_Y[f-cw]+\text{deg}_Y[w]
$$
and
$$
F_2(c)=I\left((f-cw)w,f';(f-cw)w\right)
-I\left(f'',f';(f-cw)w\right)
$$
and
$$
F_3(c)=\sum_{1\le l\le t}
\left(I(f'',f';Q_l)-I((f-cw)w,f';Q_l)\right)
$$
and
$$
F_4(c)=\sum_{1\le l\le t}
\left(I\left(X-u_l,\frac{(f-cw)w}{[f-cw][w]};Q_l\right)-1\right).
$$
Fix some $c_\pi \in k$ such that $\rho_a(f-c_\pi w) =
\rho_\pi(f,w)$ and such that $c_\pi$ gives the various extremal values
as described in (10.8). 
Explicitly, we assume that $c_\pi$ is chosen so that it satisfies the
following additional conditions (i) to (iv).

(i) For each $V\in \mathfrak R(f_i',\infty)$, with $1 \le i \le s' $, we have
$\text{ord}_V(\phi_i'(f-c_\pi w))$ = 
$\text{min}(\text{ord}_V(\phi_i'(f),\text{ord}_V(\phi_i'(w))$.

(ii) For each $V\in \mathfrak R(f_i',Q_l)$, with $1\le i\le s'$ and 
$1\le l\le t$, we have
$\text{ord}_V(\phi_i'(f-c_\pi w) = \text{min}(\text{ord}_V(\phi_i'(f),
\text{ord}_V(\phi_i'(w))$.

(iii) For $1 \le l \le t$ we have $I(X-u_i,f-c_\pi w;Q_i)=
\text{min}(I(X-u_i,f ; Q_i), I(X-u_i,w; Q_i))$.

(iv) $c_\pi \not \in \text{defset}(f,w)$ and hence in particular
$I\left(f'',f';(f-c_\pi w) \setminus w \right) =0$ by the General Singset 
Theorem.

For various numerical functions $F(c)$ to be considered, with $c$ varying
in $k$, let $H(F(c))$ denote the variation
$\sum_{c\in D}(F(c_\pi)-F(c))$ where $D$ is a finite subset of $k$ which is
defined below and which is a large enough ``defset'' to be applicable to all
the relevant $F$'s.

(i*) For each $V\in \mathfrak R(f_i',\infty)$, with $1 \le i \le s' $, 
we define $c_V\in k$ thus.
In case $\text{ord}_V\phi'_i(f)<\text{ord}_V\phi'_i(w)$,
we take $c_V = c_\pi$. In case
$\text{ord}_V\phi'_i(f) \ge \text{ord}_V\phi'_i(w)$, we take $c_V$ to be
the unique element of $k$ such that
$\text{ord}_V\phi'_i(f - c_V w)> \text{ord}_V\phi'_i(f - c_\pi w)$.
Let $D_{\text{i}}$ be the set of all these elements $c_V$.

(ii*) For each $V\in \mathfrak R(f_i',Q_l)$, with $1\le i\le s'$ and 
$1\le l\le t$, we define $c_V\in k$ thus.
In case $\text{ord}_V\phi'_i(f)<\text{ord}_V\phi'_i(w)$,
we take $c_V = c_\pi$. In case
$\text{ord}_V\phi'_i(f) \ge \text{ord}_V\phi'_i(w)$, we take $c_V$ to be
the unique element of $k$ such that
$\text{ord}_V\phi'_i(f - c_V w)> \text{ord}_V\phi'_i(f - c_\pi w)$.
Let $D_{\text{ii}}$ be the set of all these elements $c_V$.

(iii*) For $1 \le l \le t$  we define $c_l\in k$ thus. In case 
$I(X-u_l,f-c_\pi w;Q_l)=\text{max}\{I(X-u_l,f-cw;Q_l):c\in k\}$, we take 
$c_l=c_{\pi}$. In case 
$I(X-u_l,f-c_\pi w;Q_l)<\text{max}\{I(X-u_l,f-cw;Q_l):c\in k\}$,
we take $c_l$ to be the unique element of $k$ such that
$I(X-u_l,f-c_l w;Q_l)=\text{max}\{I(X-u_l,f-cw;Q_l):c\in k\}$.
Let $D_{\text{iii}}=\{c_1,\dots,c_t\}$.

(iv*) Let $D_{\text{iv}}$ 
be the union of defset$(f,w)\cap k$ and singset$(f,w)$.

Let $D=D_{\text{i}}\cup D_{\text{ii}}\cup D_{\text{iii}}\cup D_{\text{iv}}$.

Now, in view of (7*), equation (10.10) is equivalent to the equation 
$$
\rho_a(f - c_\pi w)+\rho_a(w)+t-1=\sum_{1\le j\le 4}H(F_j(c))
$$
and hence to the equation
\begin{equation*}
\sum_{1\le j\le 4}H(F_j(c))
=\sum_{1\le j\le 4}F_j(c_\pi).
\tag{10.10*}
\end{equation*}
{\it In words, (10.10*) says that the function $\sum_{1\le j\le 4}F_j(c)$, or
equivalently the function $\rho_a(f - cw)+\rho_a(w)+t-1$, replicates itself,
i.e., it has a constant value at most points and that value equals the total
variation of the function.}

STEP (8). In view of (9.6) and (10.6) we see that
\begin{equation*}
H(F_1(c)))=-\text{deg}_Y\widehat{[f,w]}+\text{deg}_Y[w].
\tag{8*}
\end{equation*}

STEP (9).
For each $V\in \mathfrak R(f_i',\infty)$, with $1\le i\le s'$, we clearly have
\begin{equation*}
-\text{ord}_V\phi'_i(f'')
=\text{ord}_V\phi'_i(X)-\text{ord}_V\phi'_i((f-c_V w)w).
\tag{9a}
\end{equation*}
Now
$$
\begin{cases}
\text{ord}_V\phi'_i(f)<\text{ord}_V\phi'_i(w)\\
\Rightarrow
\text{ord}_V\phi'_i((f-cw)w)
=\text{ord}_V\phi'_i(fw)\text{ for all }c\in k
\end{cases}
$$
and thus in this case, using (9a), we get
\begin{align*}
-H\left(\text{ord}_V\phi'_i((f- cw)w)\right)-\text{ord}_V\phi'_i(f'')
= & ~0 -\text{ord}_V\phi'_i(f'')\\
= & ~\text{ord}_V\phi'_i(X)-\text{ord}_V\phi'_i((f-c_\pi w)w).
\end{align*}
Likewise
$$
\begin{cases}
\text{ord}_V\phi'_i(f)\ge\text{ord}_V\phi'_i(w)\\
\Rightarrow\text{ord}_V\phi'_i((f-cw)w)
=\text{ord}_V\phi'_i(fw)\text{ except for exactly one }c=c_V\in k
\end{cases}
$$
and thus in this case, again using (9a), we get
\begin{align*}
& -H\left(\text{ord}_V\phi'_i((f-
cw)w)\right)-\text{ord}_V\phi'_i(f'')\\
&= -\text{ord}_V (\phi'_i((f- c_\pi w)w))
+\text{ord}_V (\phi'_i((f- c_V w)w))
-\text{ord}_V\phi'_i(f'')\\
&= ~\text{ord}_V\phi'_i(X)-\text{ord}_V\phi'_i((f-c_\pi w)w).
\end{align*}
Consequently we always have
\begin{equation*}
\begin{cases}
-H\left(\text{ord}_V\phi'_i((f- cw)w)\right)-\text{ord}_V\phi'_i(f'')\\
=\text{ord}_V\phi'_i(X)-\text{ord}_V\phi'_i((f-c_\pi w)w).
\end{cases}
\tag{9b}
\end{equation*}

Clearly
$$
H(F_2(c)) = H(I\left((f-cw)w,f';(f-cw)w\right))
-H(I\left(f'',f';(f-cw)w\right)).
\leqno(9c)
$$
Since $I\left(f'',f';(f-c_\pi w)w
\setminus w \right))=0 $ by our choice of $c_\pi$, we also have
$$
I\left(f'',f';(f-c_\pi w)w\right)=I\left(f'',f';w\right).
\leqno\text{(9d)}
$$

Now 
\begin{align*}
H(I(f'',f';(f-cw)w))
&=\sum_{c\in D}(I(f'',f';(f-c_\pi w)w)-I(f'',f';(f-cw)w))\\
&=-\sum_{c\in D}I(f'',f';(f-cw)w \setminus w)\quad\text{ by (9d)}\\
&=-I((f'',f';\mathcal A\setminus w))\\
&=I(f'',f'; w)-I(f'',f';\mathcal A)\\
&=I(f'',f'; w)+\sum_{1\le i\le s'}
\sum_{V\in \mathfrak R(f'_i,\infty)}\text{ord}_V\phi'_i(f'')
\end{align*}
and
\begin{align*}
&H(I\left((f-cw)w,f';(f-cw)w\right))\\
&=-\sum_{1\le i\le s'}
\sum_{V\in \mathfrak R(f'_i,\infty)}H(\text{ord}_V\phi'_i((f-cw)w)\\
&=\sum_{1\le i\le s'}\sum_{V\in \mathfrak R(f'_i,\infty)}
[\text{ord}_V\phi'_i(f'')+\text{ord}_V\phi'_i(X)
-\text{ord}_V\phi'_i((f-c_\pi w)w)]
&\quad\text{ by (9b)}
\end{align*}
and hence by (9c) we get
$$
H(F_2(c))=-\text{deg}_Y(f')-I(f'',f'; w)
-\sum_{1\le i\le s'}
\sum_{V\in \mathfrak R(f_i',\infty)} \text{ord}_V\phi'_i((f-c_\pi w)w).
$$
Clearly have
$\text{deg}_Y(f')=\text{deg}_Y(f_Yw-fw_Y)-\text{deg}_Y(\widehat{[f,w]})$
and, since $N \ne M$, we also have
$$
\text{deg}_Y(f')=N+M-1-\text{deg}_Y(\widehat{[f,w]}).
$$ 
Combining the above two displayed equations we conclude that
\begin{align*}
H(F_2(c))
=&1-N-M+\text{deg}_Y\widehat{[f,w]}-I(f'',f'; w)\\
&-\sum_{1\le i\le s'}
\sum_{V\in \mathfrak R(f'_i,\infty)}\text{ord}_V\phi'_i((f-c_\pi w)w)\\
\end{align*}
where the last line is clearly equal to $I((f-c_\pi w)w,f',\mathcal A)$ 
and hence by using (8*) we get
$$
\begin{cases}
H(F_1(c))+H(F_2(c))\\
=1-N-M+\text{deg}_Y[w]-I(f'',f'; w)
+I((f-c_\pi w)w,f';\mathcal A).\\
\end{cases}
\leqno\text{(9e)}
$$
Clearly deg$_Y[f-c_\pi w]=0$ and hence by (9d) and (9e) we conclude that
$$
H(F_1(c))+H(F_2(c))=F_1(c_\pi)+F_2(c_\pi).
\leqno(9^*)
$$

STEP (10). Upon letting
$$F_{3,l}(c)= I(f'',f';Q_i) - I((f-cw)w,f';Q_l)$$
we get
\begin{align*}
&H(F_{3,l}(c))\\
&=-H(I((f-cw)w,f';Q_l))\\
&=-\sum_{1\le i\le s'}\sum_{V\in\mathfrak R(f'_i,Q_l)}
H(\text{ord}_V\phi_i'((f-cw)w)\\
&=\sum_{1\le i\le s'}\sum_{V\in\mathfrak R(f'_i,Q_l)}
\left[\text{ord}_V\phi_i'(X-u_l)-\text{ord}_V\phi_i'((f-c_\pi w)w))
+\text{ord}_V\phi_i'(f'')\right]\\
&\qquad\qquad\qquad\qquad\qquad\qquad\qquad\qquad\qquad\qquad\qquad\qquad\qquad
\qquad \text{ as in (9a)}\\
&=I(X-u_l,f';Q_l)- I((f-c_\pi w)w,f';Q_l)+I(f'',f';Q_l)\\
\end{align*}
and hence
$$
\begin{cases}
F_3(c)=\sum_{1\le l\le t}F_{3,l}(c)\text{ with}\\
H(F_{3,l}(c))=F_{3,l}(c_\pi)+I(X-u_l,f';Q_l).
\end{cases}
\leqno(10^*)
$$

STEP (11). Upon letting
$$F_{4,l}(c)
=\left(I\left(X-u_l,\frac{(f-cw)w}{[f-cw][w]};Q_l\right)-1\right)$$
we clearly have
$$
\begin{cases}
F_{4,l}(c)=&I(X-u_l,f-cw;Q_l)+I(X-u_l,w;Q_l)-1\\
&-I(X-u_l, [f-cw];Q_l)-I(X-u_l,[w];Q_l).
\end{cases}
\leqno\text{(11a)}
$$
Let $\mu_l=I(X-u_l,f-c_\pi w;Q_l)$ and 
$\theta_l=\text{max}\{I(X-u_l,f-cw;Q_l):c\in k\}$ and
$\nu_l=I(X-u_l,w;Q_l)$. Then
\begin{align*}
&H(F_{4,l}(c))\\
&=H(I(X-u_l,f-cw;Q_l))-H(I(X-u_l,[f-cw];Q_l))\\
&=\mu_l-\theta_l-H(I(X-u_l,[f-cw];Q_l))\\
&=\mu_l-\theta_l-I(X-u_l,[w];Q_l)+I(X-u_l,\widehat{[f,w]};Q_l)\\
&=\mu_l-\theta_l-I(X-u_l,[w];Q_l)+[\nu_l-1+\theta_l-I(X-u_l,f';Q_l)]\\
&=I(X-u_l,(f-c_\pi w)w;Q_l)-1-I(X-u_l,f';Q_l)-I(X-u_l,[w];Q_l)
\end{align*}
and hence, because of (11a) and the obvious fact that 
$I(X-u_l, [f-c_\pi w];Q_l)=0$, we get
$$
\begin{cases}
F_4(c)=\sum_{1\le l\le t}F_{4,l}(c)\text{ with}\\
H(F_{4,l}(c))=F_{4,l}(c_\pi)-I(X-u_l,f';Q_l).
\end{cases}
\leqno\text{(11b)}
$$
By (10*) and (11b) we see that
$$
H(F_3(c))+H(F_4(c))=F_3(c_\pi)+F_4(c_\pi).
\leqno(11^*)
$$

STEP (12). By (9*) and (11*) we get (10.10*).

\centerline{}

{\bf CONCLUSION.} 
Genus plus excess branch number is the rank of a curve. The total
variation of the rank as a curve moves thru a pencil is independent of the
pencil. A suitable modification of
this variation is the Zeuthen-Segre invariant of the surface. For the
plane it equals minus one. The defset of a polynomial is the set of
translation constants which produce nongeneral rank. The defset gives a
bound on the singset, i.e., the set of translation constants which produce 
singular curves or more generally singular hypersurfaces. The redset of a 
polynomial is the set of translation constants which produce reducible
hypersurfaces. Bounds for the redset are found in terms of the group of units 
of the affine coordinate ring.

\centerline{}

\centerline{\bf Section 11: Linear Systems and Pencils on Normal Varieties}

To say a word about the Zeuthen-Segre invariant of a
surface, let us very briefly talk about linear systems and pencils on normal 
varieties.

So assume $k=k^*$, let $E$ be an irreducible $n$-dimensional normal 
algebraic variety over $k$, for $0\le i\le n$ let $E_i$ be the set of all 
irreducible $i$-dimensional subvarieties of $E$, for any 
$C\in\cup_{0\le i\le n}E_i$ let $k(C)$ be the local ring of $C$ on $E$, let
$\widehat E_i=\{k(C):C\in E_{n-i}\}$, and finally let 
$\widehat E=\cup_{0\le i\le n}\widehat E_i$.
Then, in the language of models (see \cite{A09}), $\widehat E_i$ is the set of 
all $i$-dimensional members of the 
$n$-dimensional normal model $\widehat E$ of $k(E)/k$.

Recall that a premodel $\widehat E$ of a finitely generated field extension
$K/k$ is a collection of local domains whose quotient field is $K$ and which
have $k$ as a subring; $\widehat E$ is irredundant (resp: complete)
means any valuation ring of $K/k$ dominates at most (resp: at least) one of 
its members; $\widehat E$ is a model if it is an irredundant premodel which
can be expressed as a finite union 
$\widehat E=\cup_{0\le j\le m}\mathfrak V(B_j)$ where 
$B_j=k[x_{j0},\dots,x_{jm}]=$ an affine domain over $k$ and  
where $\mathfrak V(B_j)=$ the set of all localizations $(B_j)_P$ with 
$P$ varying over spec$(B_j)$.
The normality assumption says that there is an injection $R\to B_j$ such that
$B_j$ is the integral closure of the image in a finite algebraic field 
extension of the quotient field of the image, or equivalently that every member
of $\widehat E$ is normal.  The normality assumption 
implies that $\widehat E_1$ is a set of DVRs of $k(E)/k=K/k$.
Recall that $E$ or $\widehat E$ is nonsingular means every member of
$\widehat E$ is a regular local ring, and so nonsingular $\Rightarrow$ normal.
If $m$ can be taken to be $0$ then we call $\widehat E$ (resp: $E$)
to be an affine model (resp: affine variety).
If we can find nonzero elements $z_0,\dots,z_m$ in $k(E)$ such that 
$x_{ji}=z_i/z_j$ for all $i,j$ in $\{0,\dots,m\}$ then we call $\widehat E$ 
(resp: $E$) to be a projective model (resp: projective variety).

Now any $C\in E_{n-i}$ can be recovered from $k(C)\in\widehat E_i$ 
by observing that (closed) points $P$ in $C$ are characterized by saying that
$k(P)$ are those members of $\widehat E_n$ for which $k(C)$ belongs to
$\mathfrak V(k(P))$. Thus we may dispense with the geometric beginning of
commencing with an algebraic variety, and start (and end) with a model. 
This economy of thought is the beauty of the language of models.

Let $\mathcal D$ be the group of all divisors on $E$, i.e., the set of all
maps $E_{n-1}\to \mathbb Z$ with finite support. Let $\mathcal D_+$ be the
set of all effective divisors, i.e.,  nonzero divisors $D$ with
$D(E_{n-1})\subset\mathbb N=$ the set of all nonnegative integers.
The degree deg$(D)$ of any divisor $D$ is defined to be $\sum D(C)$ taken over
all $C$ in $E_{n-1}$. The divisor $(z)$ of any $z\in k(S)^\times$ is defined
by the equation $(z)(C)=\text{ord}_{k(C)}z$.
For any $k$-vector subspace $H$ of $k(E)$, by $\mathbb P(H)$ we denote 
the associated projective space, i.e., the set of all $1$-dimensional 
subspaces of $H$; for any nonzero $z,z'$ in any $y\in\mathbb P(H)$ we clearly 
have $(z)=(z')$ and this common divisor is denoted by $(y)$.
By a linear system on $E$ we mean a subset $\mathcal C$ of $\mathcal D_+$
for which there exists a finite dimensional $k$-vector-subspace $H$
of $k(E)$ together with $D'\in\mathcal D$ such that $y\mapsto (y)+D'$ gives
a bijection $\mathbb P(H)\to\mathcal C$, and for which there does not exist
$D_0\in\mathcal D_+$ such that for all $D\in\mathcal C$ we have $D\ge D_0$;
the second proviso means that we assume our linear systems to be devoid of
fixed components.  It is easily seen that the dimension of $H$ depends only on 
$\mathcal C$, and we call this dimension decreased by one to be the 
dimension of $\mathcal C$. 
If the dimension of $\mathcal C$ is one then we call $\mathcal C$ a pencil.

Now assume that $k$ is of characteristic zero with $n=2$, and $E$ is an 
irreducible  nonsingular projective algebraic surface. For any
$C\in E_1$ let $\gamma(C)$ be its genus, and let us generalize this to any
$D\in\mathcal D_+$ by putting
$$
\gamma(D)=1+\sum_{C\in\ E_1}D(C)(\gamma(C)-1).
$$
For any $Q\in E_0$ the completion $\widetilde Q$ of $k(Q)$ is clearly 
isomorphic to $k[[X,Y]]$. For any $D\in\mathcal D_+$ we let
$\chi(D;Q)$ denote the number of branches of $D$ at $Q$, i.e., upon letting
$M$ stand for maximal ideal,
$$
\left(\prod_{\{C\in E_1:D(C)>0\text{ and }k(Q)\subset k(C)\}}
[k(Q)\cap M(k(C))]^{D(C)}\right)\widetilde Q
=U_0U_1\dots U_{\chi(D;Q)}\widetilde Q
$$
where $U_0$ is a unit in $\widetilde Q$, and
$U_1,\dots,U_{\chi(D;Q)}$ are irreducible nonunits in $\widetilde Q^\times$.
Let
$$
\overline\chi(D;Q)=\chi(D;Q)-1
$$
and put
$$
\overline\chi(D;E)=\sum_{\{Q\in E_0:\chi(D;Q)>1\}}\overline\chi(D:Q).
$$
Let rad$(D)\in\mathcal D_+$ be defined by
$$
(\text{rad}(D))(C)=\begin{cases}
1&\text{if }D(C)>0\\
0&\text{if }D(C)=0.\\
\end{cases}
$$
Also put
$$
\rho_a(D)=2\gamma(\text{rad}(D))+\overline\chi(\text{rad}(D);E).
$$
Finally define the base point set of a pencil $\mathcal C$ on $E$ by putting
$$
B(\mathcal C)=\cap_{D\in\mathcal C}S_0(D)
$$
where the curve-support and point-support of $D$ are given by
$$
S_1(D)=\{C\in E_1:D(C)\ne 0\}\text{ and }
S_0(D)=\{Q\in E_0:k(Q)\subset\cup_{C\in S_1(D)}k(C)\}.
$$

If there exists a finite subset defset$(\mathcal C)$ and an integer 
$\rho_\pi(\mathcal C)$ such that for any $D\in\mathcal C$ we have:
$\rho_\pi(\mathcal C)=\rho_a(D)\Leftrightarrow 
D\in\mathcal C\setminus\text{defset}(\mathcal C)$, then these two objects are
clearly unique and we call them the {\bf deficiency set} and the 
{\bf pencil-rank} of the pencil $\mathcal C$ on $E$.  When this is so, 

we define the Zeuthen-Segre invariant of $\mathcal C$ to be the
integer $\zeta(\mathcal C)$ given by
$$
\zeta(\mathcal C)=-|B(\mathcal C)|-2\rho_\pi(\mathcal C)
+\sum_{D\in\text{defset}(\mathcal C)}[\rho_\pi(\mathcal C)-\rho_a(D)].
$$

\centerline{}

{\bf EPILOGUE.} Assume it has been shown that defset$(\mathcal C)$ and 
$\rho_\pi(\mathcal C)$ exist for every pencil on $E$, and $\zeta(\mathcal C)$ 
depends only on $E$ and not on $\mathcal C$. Let 
$\zeta(E)$ be the Zeuthen-Segre invariant of the surface $E$, i.e.,
the common value of $\zeta(\mathcal C)$ for all pencils $\mathcal C$ on $E$.
Taking any two distinct members $F$ and $G$ of $\mathcal C$, and 
adding $\rho_a(F)+\rho_a(G)-\zeta(E)$ to both sides of the above 
equation we get the jungian formula
$$
\rho_a(F)+\rho_a(G)=-\zeta(E)-|B(\mathcal C)|
+\sum_{D\in\text{defset}(\mathcal C)
\setminus\{F,G\}}[\rho_\pi(\mathcal C)-\rho_a(D)].
$$
Thinking of $F$ and $G$ as ``curves on the surface'' $E$, their  
``common points'' or ``set-theoretic intersection'' is given by
$I^*(F,G)=S_0(F)\cap S_0(G)$, and clearly
we have $B(\mathcal C)=I^*(F,G)$.
Moreover, $F$ and $G$ have no common component, i.e., there is no $C\in E_1$
with $F(C)\ne 0\ne G(C)$. Also $\mathcal C$ is ``generated'' by
$F$ and $G$; so every member of $\mathcal C$ can symbolically be written as
$F-cG$ with $c\in k\cup\{\infty\}$ with $F-\infty G=G$. 
In the situation of the previous section, this symbolism becomes more real
by taking $E$ to be the projective plane $\mathcal P$ over $k$, and taking
$F=Z^df(X/Z,Y/Z)$ and $G=Z^dw(X/Z,Y/Z)$.
As a thought for the future, going back to the section on More General
Pencils, $\rho_\pi$ of the pencil
could be defined as $\rho_a$ of the generic member $(f,w)^{\flat}$ with
affine coordinate ring $k(f/w)[X,Y]$ over ground field $k(f/w)$.
Similar trick could be played when $E$ is any nonsingular surface.

\end{document}